\documentclass[10pt]{article}

\usepackage[backend=bibtex,doi=false,isbn=false,url=false,eprint=false]{biblatex}
\usepackage{latexsym}
\usepackage{amssymb}
\usepackage{amsthm}
\usepackage{amscd}
\usepackage{amsmath}
\usepackage{mathrsfs}
\usepackage{graphicx}
\usepackage{hyperref}
\usepackage{shuffle}
\usepackage{stmaryrd}
\usepackage{enumitem}
\usepackage[noabbrev]{cleveref}
\usepackage{ytableau}

\usepackage[all]{xy}
\input xy \xyoption{frame}
\xyoption{dvips}

\usepackage[colorinlistoftodos]{todonotes}
\usetikzlibrary{chains,scopes,decorations.markings}

\usepackage{tikz-cd}

\theoremstyle{definition}
\newtheorem* {theorem*}{Theorem}
\newtheorem* {conjecture*}{Conjecture}
\newtheorem{theorem}{Theorem}[section]

\theoremstyle{definition}

\newtheorem {claim}{Claim}
\newtheorem* {example*}{Example}

\newtheorem{lemma}[theorem]{Lemma}
\theoremstyle{definition}
\newtheorem{definition}[theorem]{Definition}
\theoremstyle{definition}

\newtheorem{conjecture}[theorem]{Conjecture}
\newtheorem{proposition}[theorem]{Proposition}
\newtheorem{corollary}[theorem]{Corollary}

\newtheorem* {remark*}{Remark}
\newtheorem{remark}[theorem]{Remark}
\theoremstyle{definition}
\newtheorem {example}[theorem]{Example}
\theoremstyle{definition}

\theoremstyle{definition}

\theoremstyle{definition}

\theoremstyle{definition}

\newenvironment{claimproof}{%
  \proof[Proof of the claim]}{\endproof}

\newcommand{\Des}{\operatorname{Des}}
\newcommand{\DesR}{\operatorname{Des}_R}
\newcommand{\DesL}{\operatorname{Des}_L}

\newcommand{\supp}{\operatorname{supp}}

\newcommand{\Ess}{\operatorname{Ess}}
\newcommand{\rank}{\operatorname{rank}}

\newcommand{\DesV}{\operatorname{Des}_V}

\newcommand{\last}{\operatorname{end}}

\newcommand{\upshift}{\mathsf{shift}_{[\uparrow]}}
\newcommand{\downshift}{\mathsf{shift}_{[\downarrow]}}

\def\topRows{\mathsf{topRows}}
\def\cols{\mathsf{cols}}

\def\({\left(}
\def\){\right)}

\newcommand{\cP}{\mathcal{P}}
\newcommand{\cQ}{\mathcal{Q}}

\newcommand{\cC}{\mathcal{C}}
\newcommand{\cD}{\mathcal{D}}
\newcommand{\cZ}{\mathcal{Z}}

\def\cX{\mathcal{X}}
\def\cW{\mathcal{W}}
\def\cY{\mathcal{Y}}

\def\ZZ{\mathbb{Z}}

\def\spanning{\textnormal{-span}}

\def\cyc{\textsf{cyc}}

\def\fG{\mathcal{G}^{(\beta)}}
\def\sh{\mathrm{sh}}

\def\fk{\mathfrak}

\def\barr{\begin{array}}
\def\earr{\end{array}}
\def\ba{\begin{aligned}}
\def\ea{\end{aligned}}
\def\be{\begin{equation}}
\def\ee{\end{equation}}

\def\qquand{\qquad\text{and}\qquad}
\def\quand{\quad\text{and}\quad}

\def\quord{\quad\text{or}\quad}

\def\inv{{\mathsf{inv}}}

\def\cH{\mathcal H}

\def\hs{\hspace{0.5mm}}

\def\ds{\displaystyle}

\def\PP{\mathbb{Z}_{>0}}
\def\NN{\mathbb{Z}_{\geq 0}}

\def\fkS{\fk S}

\def\ben{\begin{enumerate}}
\def\een{\end{enumerate}}

\def\bei{\begin{itemize}}
\def\eei{\end{itemize}}

\def\cE{\mathcal E}

\def\hs{\hspace{0.5mm}}

\def\fpf{{\mathsf {FPF}}}

\def\D{\mathsf{D}}

\def\ellhat{\hat\ell}

\def\Sfpf{\hat {\fk S}^\fpf}

\def\Ifpf{I^\fpf}

\def\e{\textbf{e}}

\newcommand{\xRightarrow}[2][]{\ext@arrow 0359\Rightarrowfill@{#1}{#2}}

\renewcommand{\O}{\mathsf{O}}
\newcommand{\Sp}{\mathsf{Sp}}

\newcommand{\cA}{\mathcal{A}}
\newcommand{\cB}{\mathcal{B}}

\newcommand{\arc}[2]{ \ar @/^#1pc/ @{-} [#2] }
\def\arcstop{\endxy\ }
\def\arcstart{\ \xy<0cm,-.06cm>\xymatrix@R=.1cm@C=.10cm }
\newcommand{\arcstartc}[1]{\ \xy<0cm,-.15cm>\xymatrix@R=.1cm@C=#1cm}

\def\ellhat{\hat\ell}
\def\Sfpf{\widehat {\fk S}^\fpf}
\def\iS{\widehat \fkS}

\def\cG{\mathcal{G}}

\def\cV{\mathcal{V}}

\def\O{\mathrm{O}}

\def\fpf{\mathsf{fpf}}

\definecolor{darkred}{rgb}{0.7,0,0} 
\newcommand{\defn}[1]{{\color{darkred}\emph{#1}}} 

\usepackage{pbox}
\tikzset{every loop/.style={min distance=10 mm, in=60, out=120, looseness=10}}
\tikzset{
dot/.style = {circle, fill, inner sep=2.2,outer sep=0},
}

\numberwithin{equation}{section}
\allowdisplaybreaks[1]
\UseCrayolaColors

\def\bpartial{\partial^{(\beta)}}
\def\bpi{\pi^{(\beta)}}
\def\fkG{\mathcal{G}^{(\beta)}}
\def\fkGO{\mathcal{G}^{\mathsf{O}}}
\def\fkGSp{\mathcal{G}^{\mathsf{Sp}}}

\def\MX{M\hspace{-0.5mm}X}
\def\Mat{\mathsf{Mat}}
\def\codim{\mathrm{codim}}

\def\cF{\mathcal{F}}

 \def\crb{\mathsf{crb}}
\def\shiftable{\mathsf{shiftable}}
\def\dom{\mathsf{dom}}
\def\iG{\widehat{\cG}}
\def\ellfpf{\ell_{\mathsf{fpf}}}
\def\ellhat{\ell_{\mathsf{inv}}}

\def\ODes{\mathrm{ODes}}

\def\SD{\mathsf{SD}}

\def\forced{\mathsf{F}}
\def\prohib{\mathsf{P}}

\def\GC{\operatorname{GC}^{\mathsf{O}}}
\def\GCSp{\operatorname{GC}^{\mathsf{Sp}}}

\def\GQ{GQ}
\def\GP{GP}

\def\Sh{\mathsf{Sh}}

\def\cB{\mathcal{B}_\inv}

\def\SubDiag{\mathsf{SubDiag}}
\def\word{\mathsf{word}}

\newcommand{\Ivex}{I^{\mathsf{vex}}_\infty}
\def\iD{\hat D}

\definecolor{lightskyblue}{RGB}{176,226,255}
\newcommand{\customvarpi}{\varpi}

\begin{document}
\title{On some Grothendieck expansions}
\author{
    Eric MARBERG\\
    Department of Mathematics \\
    Hong Kong University of Science and Technology \\
    {\tt emarberg@ust.hk}
    \and
        Jiayi WEN\\
    Department of Mathematics \\
    University of Southern California \\
    {\tt jiayiwen@usc.edu}
}

\date{}

\maketitle

\begin{abstract}
The orthogonal and symplectic groups act on the complete flag variety with finitely many orbits. The orthogonal Grothendieck polynomials $\mathcal{G}^{\mathsf{O}}_z$ and symplectic Grothendieck polynomials $\mathcal{G}^{\mathsf{Sp}}_z$ are distinguished representatives for the $K$-theory classes of the corresponding orbit closures. There is a simple formula to expand $\mathcal{G}^{\mathsf{Sp}}_z$ as a linear combination of Grothendieck polynomials $\mathcal{G}^{(\beta)}_w$, which represent the $K$-theory classes of Schubert varieties. Although the constructions of $\mathcal{G}^{\mathsf{Sp}}_z$ and $\mathcal{G}^{\mathsf{O}}_z$ are similar, finding the $\mathcal{G}^{(\beta)}$-expansion of $\mathcal{G}^{\mathsf{O}}_z$ or even computing $\mathcal{G}^{\mathsf{O}}_z$ is much harder. If $z$ is vexillary then $\mathcal{G}^{\mathsf{O}}_z$ has a nonnegative $\mathcal{G}^{(\beta)}$-expansion, but the associated coefficients are mostly unknown. This paper derives several new formulas for $\mathcal{G}^{\mathsf{O}}_z$ and its $\mathcal{G}^{(\beta)}$-expansion when $z$ is vexillary. Among other applications, we prove that the latter expansion has a nontrivial stability property.
\end{abstract}

\tableofcontents

\section{Introduction}

This article is concerned with combinatorial formulas for expansions of three different families of 
\defn{Grothendieck polynomials} related to $K$-theory classes of Schubert varieties.
We first briefly introduce these polynomials in the context of \defn{matrix Schubert varieties}.
We will then discuss the central problem of interest and some of our partial solutions.

\subsection{Grothendieck polynomials}

Let $n$ be a positive integer. Write $S_n$ for the group of permutations of the integers $\ZZ$ with support in $[n]:=\{1,2,\dots,n\}$. Let $\Mat_{n\times n}$ be the affine variety of complex $n\times n$ matrices.

For $w \in S_n$,
the \defn{matrix Schubert variety} $\MX_w$ is the set of matrices $M \in \Mat_{n\times n}$
whose upper $i\times j$ sub-matrices $M_{[i][j]}$  satisfy the rank conditions
\[\rank(M_{[i][j]}) \leq |\{ t \in [i]  : w(t) \leq j\}|\quad\text{for all $i,j\in [n]$.}\]
Results in \cite{KnutsonMiller} identify an equivariant $K$-theory class 
\[[\MX_w] \in \ZZ[a_1^{\pm1},a_2^{\pm1},\dots,a_n^{\pm1}] \cong K_T(\Mat_{n\times n }). \]
The \defn{Grothendieck polynomial} $\fkG_w$ introduced in \cite{FominKirillov} can be formed from
$[\MX_w]$
 by substituting $a_i \mapsto 1+\beta x_i$ for all $i \in [n] $
and dividing by $(-\beta)^{\codim(\MX_w)}$. 

\begin{example}\label{G-ex}
The Grothendieck polynomial of the reverse permutation $w_0 = n\cdots 321 \in S_n$
is the monomial $\fkG_{w_0} = x_1^{n-1}x_2^{n-2}\cdots x_{n-1}$.
As will be explained in Section~\ref{formulas-sect}, all other instances of $\fkG_w$ for $w \in S_n$ can be computed by applying certain 
\defn{divided difference operators} to $\fkG_{w_0}$.
When $n=3$ we have 
\[
\ba
 \fkG_{123} &=1, \\
 \fkG_{231} &= x_1,
 \ea
 \qquad
 \ba
 \fkG_{132} &=x_1+x_2 + \beta x_1x_2, \\
 \fkG_{231} &= x_1x_2,
 \ea
 \qquad
 \ba
 \fkG_{312} &=x_1^2, \\
 \fkG_{321} &= x_1^2x_2.
 \ea
 \]
\end{example}
The Grothendieck polynomial $\fkG_w$ always belongs
to
 $\NN[\beta][x_1,x_2,\dots]$ and is homogeneous of degree $\codim(\MX_w)$
 if we consider $\beta $ to have $\deg \beta=-1$. This degree coincides with the value of the usual Coxeter length function
 \be\label{ell-def}
  \ell(w) := |\{ (i,j) \in \PP\times \PP : i<j\text{ and }w(i)>w(j)\}|.
  \ee
 In turn, if $w \in S_\infty := \bigcup_{m\geq 1} S_m$
 then the value of $\fkG_w$ does 
 not 
depend on the choice of $n$ with $w \in S_n$. 
Finally,
it is  known \cite[Cor.~3.3]{MP2021Sp} that the set of polynomials  $\{\fkG_w : w \in S_\infty \}$ is a $\ZZ[\beta]$-basis for $ \ZZ[\beta][x_1,x_2,\dots]$.

\subsection{Orthogonal Grothendieck polynomials}

The main subject of this article is the following analogue of Grothendieck polynomials studied in \cite{MP2020,WyserYong}.
Let 
\[\textstyle I_n = \{ z \in S_n : z=z^{-1}\} \quand I_\infty = \bigcup_{m\geq 1} I_m = \{ z \in S_\infty : z=z^{-1}\}\]
be the sets of involutions in $S_n$ and $S_\infty$.
For $z \in I_n$ define
\[ \Mat^\O_{n\times n} := \left\{ X \in \Mat_{n\times n} : X^\top = X\right\}
\quand  \MX^\O_z = \MX_z \cap \Mat^\O_{n\times n}.\]
Results in \cite{MP2020} identify an equivariant $K$-theory class 
\[[\MX^\O_z] \in \ZZ[a_1^{\pm1},a_2^{\pm1},\dots,a_n^{\pm1}] \cong K_T(\Mat^\O_{n\times n }). \]
Similar to $\fkG_w$, the \defn{orthogonal Grothendieck polynomial} $\fkGO_z$  is formed from 
$[\MX^\O_z]$
 by substituting $a_i \mapsto 1+\beta x_i$ for all $i $
and dividing by $(-\beta)^{\codim(\MX^\O_z)}$. 
These polynomials have another geometric interpretation \cite[Thm.~2.19]{MP2020}
as $K$-theory representatives for the orbit closures of the orthogonal group acting on the complete flag variety.

Each $\fkGO_z$ belongs to  $\NN[\beta][x_1,x_2,\dots,x_n]$ and 
 does not 
depend on the choice of $n$ with $z \in I_n$.  
If we consider $\beta$ to have $\deg \beta=-1$  then  $\fkGO_z$ is homogeneous
of degree $\codim(\MX^\O_z)$. This degree can be computed by 
setting 
 \be\label{ellhat-eq}
\ba 
\cyc(z) := |\{ a \in \PP : a<z(a)\}|
\quand
\ellhat(z) := \tfrac{\ell(z) + \cyc(z)}{2} .
\ea
\ee
Then  by \cite[Lem.~5.4]{Pawlowski-universal} we have 
$ \codim(\MX^\O_z)=\ellhat(z) \in \NN$.

\begin{example}
The reverse permutation $w_0 = n\cdots 321 \in S_n$ is an involution.
For this index, Wyser and Yong derived the product formula \cite[Prop.~5]{WyserYong}
\be
\label{w0-eq}
\textstyle \fkGO_{w_0} = \prod_{1 \leq i \leq j \leq n-i} (x_i+x_j+\beta x_ix_j)\ee
but did not identify any other general expressions for $\fkGO_z$.
If $z \in I_n$ is \defn{vexillary} in the sense of being $2143$-avoiding, 
then $\fkGO_z$  can be computed using 
divided difference operators,
as will be explained in more detail in Section~\ref{formulas-sect}.
For $n=3$ such calculations give the formulas
\[
\ba
\fkGO_{1} & = 1 = \fkG_{123}, \\[-10pt]\\
\fkGO_{(1,2)} &= 2x_1 + \beta x_1^2  = 2\fkG_{213} + \beta \fkG_{312}, \\[-10pt]\\
\fkGO_{(2,3)} &= 2x_2 + 2x_1 + \beta x_2^2 + 4\beta x_1 x_2 + \beta x_1^2 + 2 \beta^2 x_1 x_2^2 + 2 \beta^2 x_1^2 x_2 + \beta^3 x_1^2 x_2^2 \\ &= 2\fkG_{132}  +\beta \fkG_{231}+ \beta \fkG_{1423} + \beta^2 \fkG_{2413}, \\[-10pt]\\
\fkGO_{(1,3)} &= 2x_1 x_2 + 2x_1^2 + 3 \beta x_1^2 x_2 +\beta x_1^3 + \beta^2 x_1^3 x_2 
\\&=2\fkG_{231} + 2\fkG_{312} + 3\beta \fkG_{321} + \beta \fkG_{4123} + \beta^2\fkG_{4213}.
\ea
\]
\end{example}

 For involutions $z \in I_\infty$ that are not vexillary,
no simple algebraic method is available for computing $\fkGO_z$. 
However, at least in the vexillary case, it is known \cite[Prop.~3.29]{MS2023} that
$ \fkGO_z $ is a linear combination of ordinary Grothendieck polynomials $\fkG_w$ with coefficients in $\NN[\beta]$. 
The nontrivial part of this fact is the non-negativity of the coefficients, which follows from a general geometric result of Brion \cite{Brion2002}.
It is an open problem to describe
the $\fkG$-expansion of $\fkGO_z$ explicitly. This problem is our  focus here.
%

 \subsection{Symplectic Grothendieck polynomials}

 It is instructive to contrast this problem with what is known about the  following \defn{symplectic Grothendieck polynomials}. 
 Let $\Ifpf_n$ be the set of fixed-point-free involutions of $\ZZ$
 sending $i \mapsto i-1$ for all even $i \notin [n]$.
 This set is in bijection with the fixed-point-free involutions in $S_n$ when $n$ is even.
 %
 For $z \in\Ifpf_n$ define
 \[ \Mat^\Sp_{n\times n} := \left\{ X \in \Mat_{n\times n} : X^\top = -X\right\}
 \quand  \MX^\Sp_z = \MX_z \cap \Mat^\Sp_{n\times n}.\]
As in the previous two cases, there is an equivariant $K$-theory class \cite{MP2020}
 \[ [\MX^\Sp_z] \in \ZZ[a_1^{\pm1},a_2^{\pm1},\dots,a_n^{\pm1}] \cong K_T(\Mat^\Sp_{n\times n }). \]
 The \defn{symplectic Grothendieck polynomial} $\fkGSp_z$  is formed from 
 the class $[\MX^\Sp_z]$ 
  by substituting $a_i\mapsto 1+\beta x_i$ for all $i \in [n]$
 and dividing by $(-\beta)^{\codim(\MX^\Sp_z)}$. 
Like $\fkGO_z$,
each
  $\fkGSp_z$ belongs to $\NN[\beta][x_1,x_2,\dots]$,  
  is homogeneous of degree $\codim(\MX^\Sp_z)$ if we consider $\deg\beta=-1$,
and  does not 
 depend on the choice of $n$ with $z \in \Ifpf_n$. 
 
 The family of symplectic Grothendieck polynomials first appeared \cite{WyserYong},
 although in a slightly different form and with another geometric interpretation.
 The symplectic group acts on the complete flag variety with orbits indexed by $\Ifpf_n$.
 In \cite[\S2.1]{WyserYong}, Wyser and Yong constructed  polynomials representing
the $K$-theory classes of the closures of these orbits.
By \cite[Thm.~2.19]{MP2020},  
these polynomials coincide  with $\fkGSp_z$ after  setting $\beta=-1$ and replacing $x_i$ with $1-x_i$.

The polynomials $\fkGSp_z$ are more tractable than $\fkGO_z$ in several respects. 
Each $\fkGSp_z$ can be determined inductively by applying divided difference operators to a product formula similar to \eqref{w0-eq} corresponding to the reverse permutation \cite[Thms.~3 and 4]{WyserYong}.
There is also a positive $\fkG$-expansion of each $\fkGSp_z$, which can be described explicitly in the following way.

  Write $w_i=w(i)$ for $w \in S_\infty$. Then let $\approx$ be the transitive closure of the relation
  on $S_\infty$ with $v^{-1} \approx w^{-1}$ if some $i \in 2\NN$
  and $a<b<c<d$ have
 \[ \{v_{i+1}v_{i+2}v_{i+3}v_{i+4},\ w_{i+1}w_{i+2}w_{i+3}w_{i+4}\} \subseteq \{adbc,\ bcad,\ bdac\},\] while $v_j=w_j$ for all $j \notin \{i+1,i+2,i+3,i+4\}$.
For example, 
 \[
 (\hs15{2634}\hs)^{-1} \approx (\hs15{3624}\hs)^{-1} \approx (\hs{15{34}}{26}\hs)^{-1} \approx (\hs{3415}26\hs)^{-1} \approx (\hs{3514}26\hs)^{-1}.
 \]
 %
 Given $z \in \Ifpf_n$ let $a_1<a_2<\dots$ be the integers with $0<a_i<b_i:=z(a_i ) $
 and define 
 $\alpha_\fpf(z)$ to be inverse of  $a_1 b_1a_2 b_2 \cdots \in S_\infty.$
Then 
we have \cite[Thm.~3.12]{MP2020}
\be\textstyle\label{fkgsp-eq}
 \fkGSp_z = \sum_{w\approx \alpha_\fpf(z)} \beta^{\ell(w) - \ellfpf(z)} \fkG_w
 \quad\text{where }\ellfpf(z) := \ell(\alpha_\fpf(z)).
 \ee
 In particular, for each $z \in \Ifpf_n$ there is a map $\GCSp_z : S_\infty \to \{0,1\}$ with 
\be\textstyle
\fkGSp_z = \sum_{w \in S_\infty} \beta^{\ell(w) - \ellfpf(z)} \cdot  \GCSp_z(w) \cdot  \fkG_w.
\ee

The formula \eqref{fkgsp-eq} is a special case of \cite[Thm.~3]{Knutson2009},
which gives a general method of decomposing $K$-theory classes of \defn{multiplicity-free subvarieties}
of generalized flag manifolds. The expansion in \cite{Knutson2009} is somewhat inexplicit as it involves coefficients that may be zero,
and in any case the method 
cannot be applied to decompose $\fkGO_z$, which does not represent a multiplicity-free variety.

\subsection{Summary of results}\label{GC-eq-sect}

Both $\fkG_w$ and $\fkGO_z$ are homogeneous if we set $\deg \beta =-1$.
As $ \left\{ \fkG_w : w \in S_\infty \right\}$ is a $\ZZ[\beta]$-basis for $\ZZ[\beta][x_1,x_2,\dots]$,
it follows that 
for each $z \in I_\infty$ 
there is a unique  coefficient function
$\GC_{z} : S_\infty \to \ZZ$
with 
\be\textstyle
\fkGO_z = \sum_{w \in S_\infty}   \beta^{\ell(w) - \ellhat(z)} \cdot \GC_z(w) \cdot \fkG_w.
\ee
%
The \defn{support} of this map, given by
\[
\supp(\GC_z) := \left\{ w \in S_\infty : \GC_z(w) \neq 0\right\}
,\]
is always a finite set of permutations. 
In contrast to $\GCSp_z$,
this set seems difficult to compute in an efficient way.

When $z \in I_\infty$ is vexillary, as noted earlier, the function $\GC_z$ takes all values in $\NN$.
(It is unknown whether this remains true when $z$ is not vexillary.)
More precisely,
define the \defn{visible descent set} of $z \in I_\infty$ to be
\be\label{vdes-def} \DesV(z) := \{i \in \PP : z(i) >  z(i+1) \leq i\},\ee
and define  the \defn{right descent set} of a permutation $w\in S_\infty$ to be
\be\DesR(w) := \{ i \in \PP  : w(i)>w(i+1)\}.\ee
Then for any vexillary $z \in I_\infty$  with $\DesV(z) \subseteq [n]$, 
one has \cite[Prop.~3.29]{MS2023}  
\be\label{vdes-eq}
 \fkGO_z \in\NN[\beta]\spanning\left\{\fkG_w : \DesR(w) \subseteq[n]\right\} \subseteq \ZZ[\beta][x_1,x_2,\dots,x_n]
.\ee
Whereas $\GCSp_z$ takes all values in $\{0,1\}$, 
  the maximum value of $\GC_z$ is unbounded as $z$ varies over all vexillary  involutions.

We mention that if one sets $\beta=0$ then $\fkG_w$, $\fkGO_y$, and $\fkGSp_z$
turn into the \defn{(involution) Schubert polynomials} $\fkS_w$, $\iS_y$, and $\Sfpf_z$ studied in \cite{HMP1,Manivel,WyserYong},
which represent cohomology classes of matrix Schubert varieties rather than $K$-theory classes.
For these polynomials, the relevant expansions are all much simpler.
Both $\iS_y$ and $\Sfpf_z$ are equal to a constant 
times a multiplicity-free sum of Schubert polynomials $\fkS_w$. Moreover, the terms that appear are predicted by a general formula of Brion \cite{Brion}
and are described combinatorially in \cite{CJW}.

This work contains the first concrete results about the more general $\fkG$-expansion of $\fkGO_z$,
or equivalently the Schubert class expansion of the $K$-theory classes of symmetric matrix Schubert varieties.
Our main theorems can be summarized as follows.

Results in \cite{MP2020} give a product formula for $\fkGO_z$ when $z$ is \defn{dominant}
and a divided difference recurrence for orthogonal Grothendieck polynomials   indexed by vexillary involutions.
We prove in Section~\ref{prelim-sect}  that these two identities determine
all instances of $\fkGO_z $ when $z$ is vexillary. This provides a (partial) orthogonal counterpart to results of Wyser and Yong \cite{WyserYong} 
showing that all symplectic Grothendieck polynomials $\fkGSp_z$ for $z \in \Ifpf_n$ are obtained by applying divided difference operators to $\fkGSp_{w_0}$.
By contrast, one cannot use divided difference operators in this way to compute $\fkGO_z $ for arbitrary (non-vexillary) involutions $z$; see the discussion in \cite[\S2.2]{WyserYong}.

 Section~\ref{exp1-sect} then derives two exact formulas for $\GC_z $ when $z$ is any \defn{quasi-dominant} involution (see Definition~\ref{qd-def} and Theorems~\ref{qd-thm} and \ref{dom-thm}). 
Several ingredients are needed to prove these theorems,
including the \defn{transition formulas} in \cite{Lenart,LenartSottile} and new identities for certain \defn{involution Grothendieck polynomials} introduced in \cite{MP2021}.
Section~\ref{ex-conj-sect} also presents a general conjecture about the support of $\GC_z$ for any vexillary $z \in I_\infty$,
and a specific prediction for this set when $z$ is a reverse permutation.

 In Section~\ref{vex-sect} we explain a new formula for the monomial expansion of $\fkGO_z$ when $z $ is any vexillary involution (see Theorem~\ref{ivex-thm}).
This somewhat complicated result is complementary to an even more complicated Pfaffian formula  for the same polynomials
that was derived in \cite{MP2020}.
Our new theorem can be viewed as the first instance of a $K$-theoretic extension of Brion's cohomology formula in \cite{Brion} beyond the multiplicity-free case considered in \cite{Knutson2009}.
This formula may be useful in future work for
determining the \defn{Castelnuovo--Mumford regularity} of symmetric matrix Schubert varieties.

As an application, we show that the
$\fkG$-expansion of $\fkGO_z$ has a nontrivial \defn{shift invariance} property when $z$ is vexillary (see Theorem~\ref{supp-thm}).
This property leads to a new proof of the existence of \defn{stable limits} of orthogonal Grothendieck polynomials,
and affirms \cite[Conj.~3.39]{MS2023} in the vexillary case.

We also derive a number of interesting identities in special cases.
 For example,   Section~\ref{gr-sect} presents a more explicit formula for $\fkGO_z$ when $z$ is \defn{I-Grassmannian}.
 This formula is a non-symmetric generalization of the main theorem in \cite{Chiu.Marberg2023},
 which expands Ikeda and Naruse's \defn{$K$-theoretic Schur $Q$-functions} $\GQ_\lambda$ in terms of their \defn{$K$-theoretic Schur $P$-functions} $\GP_\mu$.

\subsection*{Acknowledgments}

This work was partially supported by Hong Kong RGC grants  16306120 and 16304122. 
We thank Brendan Pawlowski  
for many useful conversations.

\section{Preliminaries}\label{prelim-sect}

This section reviews some facts about divided difference operators and related formulas for $\fkG_w$ and $\fkGO_z$. Then we   
develop some properties of a natural directed graph on  
vexillary involutions. Finally, we discuss a variant of $\fkGO_z$ from \cite{MP2021}.


\subsection{Divided difference operators}

Continue to let $\beta$ and $x_i$ be commuting indeterminates.
The group $S_\infty$ acts on $\ZZ[\beta][x_1,x_2,\dots]$ by permuting the $x_i$ variables.
For each $i \in \PP$ the \defn{divided difference operators}  $\partial_i$ and  $\bpartial_i$ act on $ \ZZ[\beta][x_1,x_2,\dots]$
by the formulas
\be
\partial_i f =  \tfrac{f-s_i f}{x_i-x_{i+1}} \quand \bpartial_i f = \partial_i((1+\beta x_{i+1}) f) = -\beta f + (1+\beta x_i) \partial_i f.
\ee
These operators satisfy the Coxeter braid relations for $S_\infty$ along with $\partial_i\partial_i = 0$ and $\bpartial_i \bpartial_i = -\beta \bpartial_i$.
There is a useful product formula
\be\label{fg-eq}
\bpartial_i(fg) = (s_i f) \cdot \(\bpartial_i g+\beta g\) + (\bpartial_i f)\cdot  g
\ee
which holds for any $f,g\in \ZZ[\beta][x_1,x_2,\dots]$. This implies that if $s_if =f$ then 
\be\label{fg-eq2} \bpartial_i(fg) = f \cdot \bpartial_i g.\ee

\subsection{Product formulas and recurrences}\label{formulas-sect}

The \defn{Rothe diagram} of a permutation $w \in S_\infty$ is the set
$ D(w)$ of positive integer pairs $ (i, j)$ satisfying both $ i<w^{-1}(j)$ and $j< w(i)$.
A \defn{partition} is a weakly decreasing sequence of nonnegative integers $\lambda = (\lambda_1\geq \lambda_2\geq \dots \geq 0)$ 
with a finite sum.

A permutation $w \in S_\infty$ is \defn{dominant} 
if
there is a partition $\lambda$ such that 
$D(w) =\D_\lambda$ where $\D_\lambda := \{(i,j)  : 1 \leq j \leq \lambda_i\}.$
In this case we say that $w$ is dominant of shape $\lambda$. 
The dominant permutations are exactly the elements of $S_\infty$ that are $132$-avoiding \cite[Ex.~2.2.2]{Manivel},
so each dominant permutation is also vexillary.
 
There is a unique dominant $w \in S_\infty$ of each partition shape $\lambda$,
and results in \cite{FominKirillov} show that the Grothendieck polynomial of this permutation is
\be\label{dom-fk-eq}
\textstyle \fkG_w = \prod_{(i,j) \in \D_\lambda} x_i . \ee 
For any $i \in \PP$ and $w \in S_\infty$ we have \cite{FominKirillov}
\be\label{Grothendieck-recursion}
\bpartial_i \fkG_w = \begin{cases} \fkG_{ws_i} &\text{if }i \in\DesR(w) \\ -\beta \fkG_w&\text{if }i \notin \DesR(w).\end{cases}
\ee
Combined with \eqref{dom-fk-eq}, this recursive formula  determines $\fkG_w$ for all $w \in S_\infty$.

The dominant permutation $z$ of shape $\lambda$ belongs to $I_\infty$
if and only if $\lambda=\lambda^\top$ is a \defn{symmetric} partition, since $D(z) = D(z^{-1}) = D(z)^\top$.
In this case
\be\label{idom-eq} \fkGO_z = \prod_{\substack{(i,j) \in \D_\lambda \\ i\leq j}} x_i\oplus x_j 
\quad\text{where $x\oplus y := x+y+\beta xy$}\ee by \cite[Thm.~3.8]{MP2020}.
Moreover, if $z \in I_\infty$ is vexillary then 
\be\label{orthogonal-recursion}\bpartial_i \fkGO_z = \begin{cases} \fkGO_{s_izs_i} &\text{if }i \in\DesR(z) \\ -\beta \fkGO_z&\text{if }i\notin\DesR(z)\end{cases}
\ee
for all $i \in \PP$  such that $s_i z s_i $ is also vexillary and not equal to $z$ \cite[Prop.~3.23]{MP2020}.
The recurrence \eqref{orthogonal-recursion} combined with the product formula \eqref{idom-eq}
can be used to calculate $\fkGO_z$ for any vexillary $z \in I_\infty$. 
This claim is not obvious, but will be a consequence of Theorem~\ref{vexweak-prop} in the next subsection.

 For involutions $z \in I_\infty$ that are not vexillary,
no simple algebraic method using divided difference operators is known for computing $\fkGO_z$. 
See the discussion in \cite[\S2.2]{WyserYong}, especially \cite[Examples~2.5 and 2.6]{WyserYong}.

\subsection{Vexillary involutions and weak order}

Let $\Ivex$ be the subset of vexillary elements of $I_\infty$,
that is, involutions in $S_\infty$ that are $2143$-avoiding.
 This section discusses some concrete ways of identifying the elements of this set. Both \cite{GPP} and \cite[\S4]{Pawlowski-universal} are useful references for this material.

It is sometimes useful to draw involutions $z \in I_\infty$ as \defn{arc diagrams},
which are incomplete matchings on the ordered set $\{1,2,3,\dots,n\}$ in which distinct vertices $i$ and $j$ are connected if and only if $z(i)=j$.
The value of $n$ should be chosen such that $z(i)=i$ for all $i>n$, but we do not always require $n$ to be the minimal integer with this property.
In drawing arc diagrams, we order the vertices from left to right and omit vertex labels. For example,
\[
z = (1,6)(2,4)(5,7)\quad\text{has arc diagram}\quad 
\arcstart
{
*{\bullet}    \arc{1}{rrrrr}   & *{\bullet}   \arc{0.5}{rr}    & *{\bullet}  & *{\bullet}     & *{\bullet} \arc{0.5}{rr}  & *{\bullet}    & *{\bullet}      
} 
\arcstop.
\]

\begin{remark}\label{arc-vex-rem}
The $2143$-pattern avoidance condition that characterizes $\Ivex$ is equivalent to the following property of arc diagrams \cite[Fig.~3]{GPP}:
an element of $ I_\infty$ is vexillary if and only if its arc diagram does not contain any of the following as induced vertex-ordered subgraphs:
\[
\arcstart
{
*{\bullet}  \arc{0.4}{r}  & *{\bullet}     & *{\bullet} \arc{0.4}{r}    & *{\bullet}       
} 
\arcstop
\quad
\arcstart
{
*{\bullet}  \arc{0.6}{rrr}  & *{\bullet}     & *{\bullet} \arc{0.6}{rrr}    & *{\bullet}  & *{\bullet}  & *{\bullet}       
} 
\arcstop
\quad
\arcstart
{
*{\bullet}  \arc{0.8}{rrrr}  & *{\bullet}     & *{\bullet} \arc{0.8}{rrrr}    & *{\bullet}  \arc{0.5}{rr} & *{\bullet}  & *{\bullet}  & *{\bullet}            
} 
\arcstop
\quad
\arcstart
{
*{\bullet}   \arc{0.8}{rrrr}  & *{\bullet}    \arc{0.5}{rr}   & *{\bullet} \arc{0.8}{rrrr}    & *{\bullet} & *{\bullet}  & *{\bullet}  & *{\bullet}            
} 
\arcstop
\quord
\arcstart
{
*{\bullet}  \arc{1}{rrrrr}  & *{\bullet}  \arc{0.6}{rrr}  & *{\bullet}  \arc{1}{rrrrr}  & *{\bullet}  \arc{0.6}{rrr}  & *{\bullet}  & *{\bullet}  & *{\bullet}     & *{\bullet}    
} 
\arcstop.
\]
\end{remark}

Define the \defn{vexillary weak order graph} on $\Ivex$ to be the directed graph with labeled edges $v \xrightarrow{i} w$ 
if $v,w \in \Ivex$ and $i \in \PP$ are such that 
\[w=s_i vs_i\quand \ell(v)<\ell(w).\]
One has $\ell(v)<\ell(s_i vs_i)$ precisely when $i$ and $i+1$ are not both fixed points of $v$ and   $v(i)<v(i+1)$.
In this case one always has
$\ell(s_ivs_i)  = \ell(v)+2$ but it is possible that $s_i v s_i $ is not itself vexillary.

\begin{remark}
Wyser and Yong define another \defn{weak order} relation $\prec$ on the set of all involutions $I_\infty$ in \cite[\S1.2]{WyserYong}.
If $v \xrightarrow{i} w$ is an edge in our vexillary weak order graph
then $w\prec v$ is a covering relation for this weak order.
Conversely, if $w\prec v$ is a covering relation with $v,w\in\Ivex$
then the vexillary weak order graph contains an edge  $v \xrightarrow{i} w$ for some $i$ if and only if $\cyc(v)=\cyc(w)$.
\end{remark}

For integers $p,q \in \PP$ with $q \geq 2p$, let  
\be\label{dom-def} \dom_{pq} := (1,q)(2,q-1)(3,q-2)\cdots(p,q-p+1) \in \Ivex.\ee
Also define  $ \dom_{00} = 1$.
The permutation $\dom_{pq}$ is dominant of shape 
\[
\lambda = (q-1,q-2,q-3,\dots,q-p,\underbrace{p,p,p,\dots,p}_{q-2p\text{ terms}},p-1,p-2,p-3,\dots,1).
\]

 The following terminology is new:

\begin{definition} \label{qd-def}
We say $z \in I_\infty$ is  \defn{quasi-dominant} if  $i<z(i)$ for all $i\in [\cyc(z)]$. 
\end{definition}

Since $\cyc(z)$ is the number of nontrivial cycles of $z$,
the involution $z$ is quasi-dominant precisely when it can be written in cycle notation as 
\[ z = (1,b_1)(2,b_2)\cdots (k,b_k)\]
for some numbers $b_1,b_2,\dots,b_k > k=\cyc(z)$. 
From 
this description, we see that  $z$ is quasi-dominant
if and only if $i-1<z(i-1)$ whenever $1<i<z(i)$.
It follows that if $z$ is dominant (and hence $132$-avoiding), then $z$ is also quasi-dominant.
In general, however, a quasi-dominant involution need not be dominant or vexillary.

The following property of  $\Ivex$ is possibly known to experts, but we have not been able to find a reference in the literature. 

\begin{theorem}\label{vexweak-prop}
Let $z \in \Ivex$.
Set $p=q=0$ if $z=1$ and otherwise define
\be\label{pq-eq} p =\cyc(z) = |\{ i \in \PP : i<z(i)\}| \quand
q = \max\{ i \in \PP : z(i)\neq i\}.\ee 
Then the vexillary weak order graph contains a path $z\xrightarrow{i_1} \cdots \xrightarrow{i_k} \dom_{pq}$ with \[\{i_1,\dots,i_k\} \subseteq [q-1].\]
Moreover, if $z$ is quasi-dominant, then this path 
has $p\notin\{ i_1,\dots,i_k\}.$
\end{theorem}

This result will be proved in two steps. For $i \in \PP$ and $w \in S_\infty$ define 
\be D_i(w) := \{ j : (i,j) \in D(w)\}
\quand
 \iD(w):=\{(i,j)\in D(w): i\leq j\}.\ee
Then for any finite subset $S \subset\PP\times \PP$ let
\be
\Ess(S):=\{(i,j)\in S: (i+1,j)\notin S\text{ and }(i,j+1)\notin S\}.
\ee
Note that  $z\in I_\infty$ is completely determined by $D(z)$, which by symmetry is itself determined by $\iD(z)$. 

It is useful to recall that a permutation $w \in S_\infty$ is vexillary if and only if the family of sets $\{ D_i(w): i \in \PP\}$ is totally ordered  under inclusion \cite[\S2]{Pawlowski-universal}. 
In addition, an involution
$z\in I_\infty$ is vexillary if and only if $\Ess(\iD(z))$ is totally ordered under the partial order that has $(i,j)\preceq (k,l)$ whenever $i\geq k$ and $j\leq l$ \cite[Lem.~4.18]{Pawlowski-universal}.

 For $w\in S_\infty$, define $c(w) := (c_1,c_2,c_3,\dots)$ where $c_i := |D_i(w)|$. 
It is well-known that $c(w)$ uniquely determines $w$. In particular, $w$ is dominant of shape $\lambda$ if and only if $c(w)=\lambda$. 
One has
$c_i\leq c_{i+1}$ if and only if $w(i)<w(i+1)$. 

Fix $z \in \Ivex$.
 Assume $z$ is not dominant, so that $c(z)$ is not a partition.
 Then the set of ascents $\{ i \in \PP : c_i<c_{i+1} \}$ is nonempty, and we define
the \defn{distinguished ascent} of $z$ to be 
\[
j := \min\{ i \in \PP : c_i<c_{i+1} = m\}
\]
where
\[
m:=\max\{ c_{i+1} : i\in \PP,\  c_i < c_{i+1}\}.
\]
In other words, let $j$ be the index of the ascent $c_j<c_{j+1}$ that maximizes $c_{j+1}$.

\begin{lemma}\label{vexweak-path-lem}
Suppose $z \in \Ivex$ is not dominant and has distinguished ascent $j$. Then
$j < \max\{ i \in \PP : z(i) \neq i\}$ and $z \xrightarrow{j} s_j zs_j \in \Ivex$.
\end{lemma}

\begin{proof}
 Let $q= \max\{i\in \ZZ_{>0}:z(i)\neq i\}$
 so that $z\in S_q$ and $D(z)\subset [q]\times [q]$. 
 Define $m$ as above and let $N=\max (D_{j+1}(z)).$
 
 Because  $c_{l+1}=0$ for all $l\geq q$, we must have $j<q$.
 As $c_j<c_{j+1}$, we must also have $z(j)<z(j+1)$ and there exists some $(j+1,l)\in D(z)$ such that $z(j)<l<z(j+1)$. We claim that 
 $ j\leq z(j).
 $
   To show this, we first check that 
  \be\label{claim-vex2} \varnothing\neq \{z(j)+1,z(j)+2,\cdots,N-1\}= D_{j+1}(z)\setminus D_{j}(z).\ee
  Because $z$ is vexillary, the family of sets $\{D_{i}(z): i\in \PP\}$ is totally ordered. Hence $c_j<c_{j+1}$ implies that $D_j(z)\subset D_{j+1}(z)$, and it suffices to show that 
  \[\{z(j)+1,z(j)+2,\cdots,N-1\}\subset D_{j+1}(z)\] because $z(j)>\max (D_{j}(z))$. If $i<j$ and $z(j)\in D_i(z)$, then $D_{j+1}(z)\subset D_{i}(z)$ and so $z(i)>N$. However, if $i<j$ and $z(j)\notin D_i(z)$ then $z(i)<z(j)$. Therefore \eqref{claim-vex2} holds.

  Now we prove that $j\leq z(j)$ by 
  assuming $z(j)<j$ and deriving a contradiction. By the preceding paragraph, we have 
  $(j+1,z(j)+1)\in D_{j+1}(z)$
  and so by the symmetry of $D(z)$ we also have
 $  (z(j)+1,j+1)\in D(z)$. Since $z(j)\notin D_{j+1}(z)$ and $z(j)\in D_{z(j)}(z)$, it must hold that $c_{z(j)+1}>c_{j+1}$. Because $z(j)<j$ is not the distinguished ascent, we  have 
  \[c_{z(j)}\geq c_{z(j)+1}\quad\text{and therefore}\quad D_{z(j)+1}(z)\subset D_{z(j)}(z).\] But this is impossible as $j+1\in D_{z(j)+1}(z)\subset D_{z(j)}(z)$ and $j+1>j=z(z(j))$. Therefore our claim that $j\leq z(j)$ must hold.

Finally, we prove that $s_jzs_j\in \Ivex$ by showing that $\Ess(\iD(s_jzs_j))$ is totally ordered under $\preceq$. We observe that $D(s_jzs_j)$ and $D(z)$ differ only in the rows and columns indexed by $j$ and $j+1$. Because $z(j)\geq j$ and $z(j+1)>j+1$, 
the half-diagrams $\iD(z)$ and $\iD(s_jzs_j)$ have the same  columns $j$ and $j+1$. Moreover,  
\[\iD_{j+1}(s_jzs_j)=\iD_j(z)\quand \iD_j(s_jzs_j) = \iD_{j+1}(z)\cup \{z(j)\}.\] 
  Because $(j,z(j)+1)\in \iD(s_jzs_j)$, the pair $(j,z(j))$ cannot be in $\Ess(\iD(s_jzs_j))$ and $\Ess(\iD(s_jzs_j))\setminus \Ess(\iD(z))$ has at most an element. In fact:
\bei
\item    If $(j+1,N)\notin \Ess(\iD(z))$ then $\Ess(\iD(z)) = \Ess(\iD(s_jzs_j))$. 
\item If $(j+1,N)\in \Ess(\iD(z))$ then $(j,N)$ is in $\Ess(\iD(s_jzs_j))$, but this set is still totally ordered as any $(k,l)\succeq (j+1,N)$ must have $k<j$ and $l\geq N$. 
\eei
We conclude that $\Ess(\iD(s_jzs_j))$ is totally ordered under $\preceq$ so $s_iz s_i \in \Ivex$.
\end{proof}

Now suppose $z \in \Ivex$ is dominant and has $p>0$ nontrivial cycles.
Define 
\be\label{hat-d-eq}
\iD_i(z) := \{ j : (i,j) \in \iD(z)\}
\quand \hat c(z) := (\hat c_1, \hat c_2,\hat c_3,\dots)
\ee
where $\hat c_i := |\iD_i(z)|$.
Because $z$ is dominant, we have
\be
D_i(z) = [i-1+\hat c_i]\quand\iD_i(z) = i-1+ [\hat c_i]\quad\text{for all $i\leq p$,}
\ee
along with $\hat c_i = 0$ and $z(i)\leq i$ for all $i>p$.
%
%

Continue to define $q$ as in \eqref{pq-eq}. 
Then we must have 
$D(z) \subseteq \{ (i,j) : i+j\leq q\}$ so $\hat c_i \leq q-2i+1$ for all $i \in [p]$,
with equality for all $i$ if and only if $z=\dom_{pq}$.
Hence, if $z$ is dominant but not equal to $ \dom_{pq}$ then some minimal $i \in [p]$ has $\hat c_i < q - 2i+1$;
 we refer to $j := z(i)$ as the \defn{dominant ascent} of $z$.
 
 Finally, let $\e_i$ be the tuple with $1$ in position $i$ and zeros in all other positions.
 
\begin{lemma}\label{dominant-path-lem}
Let $z \in \Ivex\setminus\{1\}$ be dominant. Define $p$ and $q$ as in \eqref{pq-eq} and assume $z \neq \dom_{pq}$.
Let $j$ be the dominant ascent of $z$ and set $i = z(j)$. 
Then $j\leq q-i $ and $s_j z s_j $ is dominant with $\ell(s_j z s_j) = \ell(z) + 2$ and $\hat c(s_j z s_j) = \hat c(z) + \e_i$.
\end{lemma}

\begin{proof}
Since $\hat c_t = q-2t+1$ for each $t \in [i-1]$,
we must have 
$z(t) = q-t+1$ for all $t \in [i-1]$.
Then, in order to have $\hat c_i < q-2i+1$, we must also have $j=z(i) < q-i+1$.

Every $t \in [p]$ must have $t <z(t)$ since $z$ has $p$ nontrivial cycles but $z(t) \leq t$ for all $t >p$.
Therefore $i<z(i)=j$ since $i \in [p]$. As $\ell(z) = |D(z)|$, 
it suffices to show that $D(s_j z s_j) = D(z) \sqcup\{(i,j)\}\sqcup\{(j,i)\}$ since the union on the right is 
 the diagram of a symmetric partition. 

The diagrams $D(z)$ and $D(s_jzs_j)$ only differ  in the rows and  columns indexed by $j$, $j+1$, $z(j)=i$, and $z(j+1)$.
Because $D(z)$ is of partition shape, we must have $z(k)<j$ for any $i<k<z(j+1)$. Hence $z(j+1)\geq j+1$ and $(i,j)\in D(s_jzs_j)$, but $D(s_jzs_j)$ has no pairs of the form $(k,j)$ with $k> i$. 
 By the symmetry of the Rothe diagram, $D(s_jzs_j)$ contains $(j,i)$ and no pairs of the form $(j,k)$ with $k>i$.
 Because $i<j$ and $j+1\leq z(j+1)$, the diagrams $\iD_j(s_jzs_j)$ and $\iD_{j}(z)$ have the same column $i$ and same column $z(j+1)$ unless $z(j+1)=j+1$.
 Thus $D(s_jzs_j)=D(z) \sqcup\{(i,j)\}\sqcup\{(j,i)\}$ as needed.
\end{proof}

We can now give a short proof of Theorem~\ref{vexweak-prop}.
\begin{proof}[Proof of Theorem~\ref{vexweak-prop}]

If $z$ is not dominant, or dominant but not equal to $\dom_{pq}$, 
then we can apply Lemma~\ref{vexweak-path-lem} or \ref{dominant-path-lem} to 
find $j \in [q-1]$ with $z \xrightarrow{j} s_j z s_j \in \Ivex$.
As the values of $p$ and $q$ do not change when $z$ is replaced by $s_j z s_j $,
this one-step path can be extended by induction to one of the desired form.

Assume $z$ is quasi-dominant.
In this case $z(i) \leq i$ for all $ i>p$,  so $\ell(s_p z s_p) \leq \ell(z)$.
Thus, 
we can only have $z \xrightarrow{j}v \in \Ivex$ if $j\neq p$, and then $v$ is quasi-dominant with $\cyc(v) = \cyc(z)=p$, so by induction any path $z\xrightarrow{i_1} \cdots \xrightarrow{i_k} \dom_{pq}$ must have $p\notin\{ i_1,\dots,i_k\}.$
\end{proof}

\subsection{Involution Grothendieck polynomials}\label{iG-sect}

There is another family of polynomials indexed by involutions $z \in I_\infty$ that share several favorable algebraic properties with $\fkGSp_z$.
We shall see that these \defn{involution Grothendieck polynomials}
are closely related to $\fkGO_z$, although they do not currently have a geometric interpretation by themselves.

Let $\leq$ denote the \defn{(strong) Bruhat order} on $S_\infty$,
which is the transitive closure of the relation with $v<w$ whenever $\ell(v) < \ell(w)$ and $w=vt$ for any transposition $t=(i, j) \in S_\infty$.
For each $u, v \in S_\infty$ there is a unique element $w \in S_\infty$ such that the set-wise product 
$
\{ x\in S_\infty : x \leq u\} \{ x \in S_\infty : x \leq v\} $ is equal to $ \{x \in S_\infty : x \leq w\}$.
If we set $u\circ v := w$ then $\circ$ gives an associative operation $ S_\infty \times S_\infty \to S_\infty$, often called the \defn{Demazure product}. 

Standard properties of the Bruhat order \cite[\S2.2]{CCG} imply that 
 $u\circ v = uv$ if and only if $\ell(uv) = \ell(u)+\ell(v)$, that $s_i\circ s_i = s_i$ for all simple transpositions $s_i:=(i,i+1) \in S_\infty$, and that $(u\circ v)^{-1} = v^{-1} \circ u^{-1}$.
 It follows 
 that
\be\label{circ-eq}
 s_i\circ w =\begin{cases} s_iw &\text{if }i \notin \DesL(w) \\ w &\text{if }i \in \DesL(w)\end{cases}
 \ \ \text{and}\ \ 
 w\circ s_i =\begin{cases} ws_i &\text{if }i \notin \DesR(w) \\ w &\text{if }i \in \DesR(w)\end{cases}\ee
for all $w \in S_\infty$ and $i \in \PP$, where $\DesL(w) :=  \DesR(w^{-1})$.
 
 \begin{remark}\label{bcs-rem}
 There is a monomial-positive formula for $\fkG_w$ that can be stated using this notation. In what follows, 
 we use the term \defn{word} to mean a finite sequence of positive integers $a=a_1a_2\cdots a_k$.
Given such a word, we set $\ell(a) =k$ and $x^a = x_{a_1}x_{a_2}\cdots x_{a_k}$.

A pair of words $(a,i)$ with $\ell(a) = \ell(i)$ is a 
\defn{compatible sequence} 
if 
  $ i $
 is weakly increasing such that
   $i_j < i_{j+1}$ whenever $a_j \leq a_{j+1}$.
Such   a compatible sequence is \defn{bounded} if $i_j \leq a_j$ for all $j$.
   
Given  $w \in S_\infty$,  let $\cH(w)$ be the set of \defn{Hecke words} $a=a_1a_2\cdots a_k$ with $w =s_{a_1}\circ s_{a_2}\circ \cdots \circ s_{a_k}$. Then by \cite[Cor.~5.4]{KnutsonMiller2004} one has
 \be\textstyle\fkG_w = \sum_{(a,i)} \beta^{\ell(a) - \ell(w)} x^i\ee
 where the sum is over all bounded compatible sequences $(a,i)$ with $a \in \cH(w)$.
 \end{remark}
 
Combining the identities in \eqref{circ-eq} shows that if
$z\in I_\infty$ and $i \in \PP$ then
\be\label{iwi-eq}
s_i\circ  z\circ s_i =\begin{cases} zs_i &\text{if }i \notin \DesR(z)\text{ and }zs_i=s_iz \\
s_izs_i &\text{if }i \notin \DesR(z)\text{ and }zs_i\neq s_iz \\
z&\text{if }i \in \DesR(z).\end{cases}
\ee
Note that $zs_i=s_iz$ if and only if $\{i,i+1\} = \{z(i),z(i+1)\}$.
The formula \eqref{iwi-eq} implies that $ w\mapsto w^{-1}\circ w$ is a surjective map $S_\infty\to I_\infty$,
so the set 
\be \cB(z) := \{w \in S_\infty : w^{-1} \circ w = z\}\ee is nonempty for all $z \in I_\infty$.


\begin{definition}\label{iG-def} The \defn{involution Grothendieck polynomial} of $z \in I_\infty$ is
\[ \textstyle \iG_z := \sum_{w \in \cB(z)} \beta^{\ell(w) - \ellhat(z)} \fkG_w \in \NN[\beta][x_1,x_2,x_3,\dots].\]
\end{definition}

The set $\cB(z)$ was extensively studied in \cite{HMP2}, which referred to its elements
as the \defn{Hecke atoms} for $z$.
The polynomials $\iG_z$ were introduced in \cite[\S4]{MP2021}.

The set $\cB(z)$ can be efficiently generated using a certain equivalence relation.
Let $\sim$ be the transitive closure of the relation
 on $S_\infty$ that has $v^{-1} \sim w^{-1}$ if there is an index $i \in \PP$
 and integers $a<b<c$ such that 
\[
\{ v_{i}v_{i+1}v_{i+2} ,\ w_{i}w_{i+1}w_{i+2} \}\subset \{cba,\ cab,\ bca\}\] while $v_j=w_j$ for all $j \notin \{i,i+1,i+2\}$.
This is the transitive relation   with 
\[ (\cdots cab \cdots)^{-1} \sim (\cdots cba \cdots)^{-1} \sim (\cdots bca \cdots)^{-1}\quad\text{for all }a<b<c,\]
where corresponding ellipses ``$\cdots$'' mask identical subsequences.

For $z \in I_\infty$ let $a_1<a_2<\dots$ be the positive integers with $a_i\leq b_i:=z(a_i ) $.
Define 
$\alpha_\inv(z)$ to be inverse of the permutation whose one-line representation is formed by removing the repeated letters from
$b_1 a_1b_2 a_2 \cdots.$
Then \cite[\S6.1]{HMP2} 
\be\label{sim-eq} \cB(z) = \{ w \in S_\infty : w \sim \alpha_\inv(z)\}.
\ee
Recall the definition of $\ellhat(z)$ from \eqref{ellhat-eq}.
It is also known \cite[\S3]{HMP1} that 
\be\label{ellhat-eq1}
\ellhat(z) =  \ell(\alpha_\inv(z))=\min\{\ell(w) : w \in \cB(z)\} = |\iD(z)|.
\ee
From \eqref{ellhat-eq} and \eqref{iwi-eq}  we get
\be\label{ellhat-eq2}
\ellhat(s_i\circ z \circ s_i) =\begin{cases} \ellhat(z) &\text{if }i \in \DesR(z) \\
\ellhat(z)+1 &\text{if }i \notin\DesR(z).\end{cases}
\ee

\begin{example}\label{cB-ex2}
If $z = (1,4)(2,5) \in I_\infty$ then $\alpha_\inv(z) = 41523^{-1} = 24513$ and
\[
\ba \cB(z) &= \{24513,\ 25413,\ 25314,\ 35214,\ 35124\} \\&= \left\{41523^{-1},\ 41532^{-1},\ 41352^{-1},\ 43152^{-1},\ 34152^{-1}\right\}
.
\ea
\]
\end{example}

\section{Expansions for quasi-dominant involutions}\label{exp1-sect}

Recall from Definition~\ref{qd-def} that an involution $z\in I_\infty$ is quasi-dominant if one has $i<z(i)$ for all integers $1 \leq i \leq \cyc(z)$,
and that every dominant involution is quasi-dominant. 
This section proves two new formulas for $\fkGO_z$ when $z$ is quasi-dominant; see Theorems~\ref{qd-thm} and \ref{dom-thm}.
We then present a general conjecture about the support of $\GC_z$ for vexillary involutions $z \in \Ivex$.

\subsection{Grothendieck expansions in the quasi-dominant case}\label{qd-sect}

For any element $z \in I_\infty$ we define
\be\label{k-eq} k(z) :=
\min \{ i \in \NN: z(j) \leq j\text{ for all } j>i\}.
\ee
Since $D(z)$ is invariant under transpose when $z\in I_\infty$, we also have
\be k(z)=\min\{ i \in \NN : (j,j) \notin D(z) \text{ for all } j>i\}.\ee

\begin{proposition}
Let $z \in I_\infty$. Then it holds that $ k(z) \geq\cyc(z)$,
with equality if and only if $z$ is quasi-dominant.
\end{proposition}

\begin{proof}
As $\{ a \in \PP : a<z(a)\} \subseteq \{1,2,\dots,k(z)\}$
 we must have $\cyc(z)\leq k(z)$, with equality only if $z$ is quasi-dominant.
Conversely, if $z$ is quasi-dominant then $z(j) \leq j$ for all $j > \cyc(z)$, so $k(z) \leq \cyc(z)$
whence $k(z) =\cyc(z)$. 
\end{proof} 

The goal of this section is to prove the following theorem.

\begin{theorem}\label{qd-thm}
If $z \in \Ivex$ is quasi-dominant 
then $\fkGO_z = \iG_z \prod_{i=1}^{k(z)} (2 + \beta x_i).$
\end{theorem}

The proof of this result appears at the end of this section.
Before giving this argument, we need to derive a few new properties of 
the polynomials $\iG_z$.

\begin{theorem}\label{iG-thm}
 Suppose $z \in I_\infty$. Then for any $i \in \PP$  it holds that
\[\bpartial_i \iG_z = \begin{cases} \iG_{zs_i} &\text{if $i \in \DesR(z)$ and $z(i) = i+1$}
\\ 
\iG_{s_i z s_i} &\text{if $i \in \DesR(z)$ and $z(i) \neq i+1$}
\\
-\beta \iG_z &\text{if $i \notin \DesR(z)$.}\end{cases}
\]
\end{theorem}
\begin{proof}
    Let $z \in I_\infty$  and  $i \in \PP$.
    If $i \notin \DesR(z)$ 
    then $z\neq s_i \circ z \circ s_i$ by \eqref{iwi-eq} and every $w \in \cB(z)$ has $w\circ s_i \in \cB(s_i \circ z \circ s_i)$.
    In this case,
    since  $\cB(z)$ and $\cB(s_i\circ z\circ s_i)$ are disjoint sets,
     every $w \in \cB(z)$ must have $w\neq w \circ s_i$, which means that $i \notin \DesR(w)$
    and $\bpartial_i \fkG_w = -\beta \fkG_w$ so $\bpartial_i \iG_z =-\beta \iG_z$.

    Instead suppose $i\in\DesR(z)$. 
    Define $y \in I_\infty$ to be   $ zs_i =s_iz$  if $z(i)=i+1$ or else $s_izs_i$ if $z(i)\neq i+1$.
    Then it follows from \eqref{iwi-eq} that 
    \ben
    \item[(a)] one has $\{y,z\} = \{ v \in I_\infty : z=s_i \circ v\circ s_i\}$ and $\ellhat(z) = \ellhat(y)+1$.
    \een
    Next, choose any $w \in \cB(z)$. If $i \in \DesR(w)$ then $ws_i \circ s_i = w$ so 
    \[
    s_i \circ (ws_i)^{-1} \circ (ws_i) \circ s_i = (ws_i\circ s_i)^{-1} \circ (ws_i \circ s_i) = w^{-1} \circ w = z.
    \]
    Alternatively, if $i \notin \DesR(w)$ then $ws_i = w\circ s_i$ so 
    \[
    s_i \circ (ws_i)^{-1} \circ (ws_i) \circ s_i= s_i\circ s_i\circ w^{-1} \circ w \circ s_i \circ s_i= s_i \circ z \circ s_i = z.
    \]
    Given these identites we conclude from property (a) that
    \ben
    \item[(b)] if $w \in \cB(z)$ then $ws_i \in \cB(y) \sqcup \cB(z)$.
    \een
    Finally, if $v \in \cB(y)$ then we certainly have $v\circ s_i \in \cB(z)$ as 
    \[ (v\circ s_i)^{-1} \circ (v\circ s_i) = s_i \circ v^{-1} \circ v \circ s_i = s_i \circ y \circ s_i = z.\]
    As $\cB(y)$ and $\cB(z)$ are disjoint, we must also have $v\neq v\circ s_i = vs_i$ so
    \ben
    \item[(c)] if $v \in \cB(y)$ then $i \notin \DesR(v)$ and $vs_i \in \cB(z)$.
    \een
Items (b) and (c) imply
    that $\cB(z) = \cC \sqcup \cD \sqcup \cE$ 
    is the union of  the disjoint sets
    \[
    \ba
    \cC &:= \{ vs_i : v \in \cB(y)\},
    \\
    \cD &:= \{ w \in \cB(z) : i \notin \DesR(w)\},
    \\
    \cE &:= \{ ws_i : w \in \cD \}.
    \ea
    \]
    We can therefore write
    \[\ba
    \bpartial_i \iG_z 
    &= \sum_{w \in \cC\sqcup \cD\sqcup \cE} \beta^{\ell(w) - \ellhat(z)} \bpartial_i \fkG_w 
    \\&= \sum_{v \in \cB(y)} \beta^{\ell(vs_i) - \ellhat(z)} \bpartial_i \fkG_{vs_i} + \sum_{w \in \cD} \beta^{\ell(w) - \ellhat(z)} \(\bpartial_i \fkG_w +\beta  \bpartial_i \fkG_{ws_i}\).
    \ea
    \]
    The last expression is equal to $\iG_y$, since  \eqref{Grothendieck-recursion} implies that if $v \in \cB(y)$ then 
    \[
    \beta^{\ell(vs_i) - \ellhat(z)} \bpartial_i \fkG_{vs_i} = \beta^{(\ell(v)+1) - (\ellhat(y)+1)}  \fkG_v = \beta^{\ell(v) -\ellhat(y)}  \fkG_v
    \]
    while if $w \in \cD$ then $\bpartial_i \fkG_w +\beta  \bpartial_i \fkG_{ws_i} = -\beta \fkG_w +\beta\fkG_{w}=0$.
    \end{proof}
    The next result is closely related to \cite[Thm.~4.5]{MP2021}.
The latter theorem
gives a general formula for $\iG_z$ as a sum of polynomial weights 
attached to certain \defn{involution Hecke pipe dreams} for $z \in I_\infty$. 
It will turn out that this sum has only one term in the case when $z \in I_\infty$ is dominant,
but to explain this we must recall a few definitions.

Let $\SubDiag = \{(i,j) \in \PP\times \PP : i\geq j\}$ and write $(i,j) \prec (i',j')$ if $i<i'$ or if $i=i'$ and $j>j'$. 
For any subset $D\subset\SubDiag$ let $\word(D)$ be the sequence formed by listing the numbers $i+j-1$ as $(i,j)$ runs over
$D$ in the total order $\prec$. For example,
$ \word(\{(1,2),(1,1),(2,3)\}) = 214.$

Let $D\subset \SubDiag$ 
and write $\word(D) = w_1w_2\cdots w_{k}$.
Then $D$ is an \defn{involution Hecke pipe dream} for $z \in I_\infty$ if 
$ s_{w_1} \circ s_{w_2}\circ \cdots \circ s_{w_k} \in \cB(z)$.
In this case we say that $D$ is \defn{reduced} if its size $|D|=k$ is equal to $\ellhat(z)$.
By \cite[Thm.~5.8]{HMP6} every $z\in I_\infty$ has at least one reduced involution Hecke pipe dream.

\begin{lemma}\label{oneplus-lem}
Suppose $D$ is an involution Hecke pipe dream for $z \in I_\infty$ and $(i,j) \in D$.
Then $z$ has a reduced involution Hecke pipe dream containing $(i,j)$.
\end{lemma}

\begin{proof}
Suppose $\word(D) = w_1w_2\cdots w_l$. Define $z_0=1$ and $z_k = s_{w_k}\circ z_{k-1} \circ s_{w_k}$ for $k \in [l]$. 
If $D$ is reduced then the desired property is immediate, so assume $D$ is not reduced.
Then it follows from \eqref{ellhat-eq2} that some $q \in [l]$ has $w_q \in \DesR(z_{q-1})$.
If $q$ is chosen to be minimal, then $q>1$ and the \defn{exchange property} stated as \cite[Prop.~3.10]{Hultman2007}
implies that there exists an index $p \in [q-1]$ such that
\[ s_{w_1}\circ \cdots \circ s_{w_p} \circ \cdots \circ \widehat{s_{w_q}} \circ \cdots \circ s_{w_l} =s_{w_1}\circ \cdots \circ \widehat{s_{w_p}} \circ \cdots \circ s_{w_q} \circ \cdots \circ s_{w_l} \in \cB(z).\]
The letters $w_p$ and $w_q$ correspond to two different positions in $D$. Only one of these positions can be $(i,j)$, so by removing the other we obtain an involution Hecke pipe dream for $z$ that still contains $(i,j)$ while having strictly fewer elements. From this point, the lemma follows by induction on the size of $D$.
\end{proof}

\begin{lemma}\label{iG-thm2-lem}
Suppose $z \in I_\infty$ is dominant of shape $\lambda$. Then $\D_\lambda \cap \SubDiag$ 
is the unique involution Hecke pipe dream for $z$.
\end{lemma}

\begin{proof}
  \cite[Thm.~3.9 and Prop.~4.15]{HMP6} show that  $\D_\lambda \cap \SubDiag$ is the unique reduced involution Hecke pipe dream for $z$.  By Lemma~\ref{oneplus-lem},
every involution Hecke pipe dream for $z$ is contained in, and hence equal to, this set.
 \end{proof}

The  $\beta=0$ case of the following product formula appears in \cite[Thm.~3.26]{HMP1}.
\begin{theorem}\label{iG-thm2}
If $z \in I_\infty$ is dominant of shape $\lambda$ then 
\[ \iG_z = \prod_{\substack{(i,j) \in \D_\lambda \\ i = j}} x_i \prod_{\substack{(i,j) \in \D_\lambda \\ i < j}} x_i\oplus x_j
\quad\text{where $x\oplus y := x+y+\beta xy$}.\]
\end{theorem}

\begin{proof}
From Lemma~\ref{iG-thm2-lem}, this result is a special case of \cite[Thm.~4.5]{MP2021}, which states that if $z \in I_\infty$ then $ \textstyle \iG_z = \sum_D \(\prod_{(i,j) \in D, i = j} x_i \prod_{(i,j) \in D, i \neq j} x_i\oplus x_j\)$ where $D$ runs over all involution Hecke pipe dreams for $z$.
\end{proof}


\begin{example}
Let  $z = (1,4)(2,5) \in I_\infty$. Then $z$ is dominant of shape $(3,3,2)$, so from Example~\ref{cB-ex2} and Theorem~\ref{iG-thm2}
we get the two expressions
\[ \ba \iG_z &=
\fkG_{24513} + \beta\fkG_{25413} + \fkG_{25314} +\beta \fkG_{35214} + \fkG_{35124}
\\
&= x_1x_2(x_1\oplus x_2)(x_1\oplus x_3)(x_2\oplus x_3)
.\ea
\]

\end{example}

We may now deliver the postponed proof of Theorem~\ref{qd-thm}.

\begin{proof}[Proof of Theorem~\ref{qd-thm}]
    Let $z\in \Ivex$ be quasi-dominant with $p=\cyc(z) = k(z)$. 
    Set $q = \max\{i \in \PP : z(i)\neq i\}$.
    By Theorem~\ref{vexweak-prop}, there is a  directed path $z\xrightarrow{i_1} \cdots \xrightarrow{i_k} \dom_{pq}$ in the vexillary weak order graph with $p\notin\{ i_1,\dots,i_k\}$,
so 
\[\textstyle \fkGO_z = \bpartial_{i_1} \cdots \bpartial_{i_k}  \fkGO_{\dom_{pq}}= \bpartial_{i_1} \cdots \bpartial_{i_k} \( \iG_{\dom_{pq}} \prod_{i=1}^{p}(2+\beta x_i)\)\]
by \eqref{idom-eq}, \eqref{orthogonal-recursion}, and Theorem~\ref{iG-thm2}.
As $\prod_{i=1}^{p}(2+\beta x_i)$ is symmetric in $x_i$ and $x_{i+1}$ when $i \neq p$,
 we can use  \eqref{fg-eq2} and Theorem~\ref{iG-thm} to  rewrite the last expression as
 the desired product
$\textstyle (  \bpartial_{i_1} \cdots \bpartial_{i_k}  \iG_{\dom_{pq}} ) 
 \prod_{i=1}^{p}(2+\beta x_i)
=  \iG_z \prod_{i=1}^{p}(2+\beta x_i)$. 
\end{proof}

\subsection{Transition formulas and applications}

Results in  \cite{Lenart,LenartSottile}
give explicit \defn{transition formulas} for certain products of Grothendieck polynomials. We review these here and then use them to derive an expression for $\GC_z$ when $z$ is quasi-dominant.
Following  \cite{LenartSottile}, we write 
\[ v\xrightarrow{(a,b)}w\qquad\text{for $v,w \in S_\infty$ and positive integers $a<b$}
\] to indicate that $w = v(a,b)$ and $\ell(w) = \ell(v)+1$, meaning that $w$ covers $v$ in the Bruhat order on $S_\infty$.
 The length condition holds precisely when $v(a) < v(b)$ and 
no $i$ with $a<i<b$ has $v(a) < v(i) <v(b)$. 

\begin{theorem}[{Lenart \cite[Thm.~3.1]{Lenart}}]
\label{lenart-transition} If $k \in \PP$ and $v \in S_\infty$ then
\[\textstyle
(1+\beta x_k)  \fkG_v = \sum_\gamma 
(-1)^{p} \cdot \beta^{p+q}
\cdot  \fkG_{\last(\gamma)}
\]
where the sum is over all saturated chains $\gamma$ in Bruhat order on $S_\infty$ of the form
{\be\label{lenart-eq}
v = v_0 \xrightarrow{(a_1,k)}
\cdots  \xrightarrow{(a_p,k)} v_p  \xrightarrow{(k,b_1)} 
\cdots 
 \xrightarrow{(k,b_q)}v_{p+q} = \last(\gamma)
 \ee}%
with $p,q\in \NN$ and
$ a_p < a_{p-1} < \dots < a_1 < k < b_q< b_{q-1} <\dots <b_1.$
\end{theorem}

A permutation $w \in S_\infty$ is \defn{$n$-Grassmannian} 
if $w(i) < w(i+1)$ for all positive integers $i\neq n$.
We say that an element of $S_\infty$ is \defn{Grassmannian} if it is $n$-Grassmannian for some $n$.
The identity permutation is the unique element of $S_\infty$ that is $n$-Grassmannian for all $n \in \NN$.

\begin{remark}\label{grass-sym-rem}
The polynomials $\fkG_w$ indexed by $n$-Grassmannian permutations are a $\ZZ[\beta]$-basis 
for the subring of  $S_n$-invariants in $\ZZ[\beta][x_1,x_2,\dots,x_n]$ by \cite[Thm.~2.2]{lenart2000}.
In particular, $\fkG_w$ is a symmetric element of $\ZZ[\beta][x_1,x_2,\dots,x_n]$ if and only if $w$ is $n$-Grassmannian. 
\end{remark}

Given a partition $\lambda$ with at most $n$ nonzero parts, 
there is a unique $n$-Grassmannian permutation $w \in S_\infty$
with 
$(w(1) -1, w(2) -2,\dots, w(n) - n) = (\lambda_n,\dots,\lambda_2, \lambda_1).$
We denote this permutation using the notation $[\lambda \mid n]$.
 
 \begin{example}\label{Gr-ex}
In one-line notation we have 
$  [1^n \mid n] = 234\cdots (n+1)1$ and 
\[ [ n \mid n] = 1,2,3,\cdots, n-1,2n,n,n+1,n+2,\dots,2n-1,\]
along with $[\emptyset\mid n] = 1$ for all $n$.
For $0<j<n$ we have 
\[[1^{n-j} \mid n] = 1,2,3,\cdots,j,j+2,j+3,j+4,\dots,n+1,j+1.\]
\end{example}

In the following lemma, we interpret $[1^0\mid n]$ as the permutation $[\emptyset\mid n] = 1$.
\begin{lemma}\label{1gr-lem}
If $0\leq j\leq n$ then 
\[ (1+\beta x_{n+1})  \fkG_{[1^{j} \mid n]} =  \sum_{i=j}^n (-\beta)^{i-j}  \fkG_{[1^i\mid n]} -  \sum_{i=j+1}^{n+1} (-\beta)^{i-j} \fkG_{[1^{i}\mid n+1]}.\]
\end{lemma} 

\begin{proof}
By inspecting Example~\ref{Gr-ex}, we see that if $v_0= [1^{j} \mid n]$ and $k=n+1$, then
the only way to construct a saturated chain in Bruhat order 
\[
v_0 \xrightarrow{(a_1,k)}v_1 \xrightarrow{(a_2,k)}\cdots  \xrightarrow{(a_p,k)} v_p  \xrightarrow{(k,b_1)}v_{p+1}  \xrightarrow{(k,b_2)}
\cdots 
 \xrightarrow{(k,b_q)} v_{p+q}
 \]
 with $a_p<a_{p-1}<\dots<a_1<k<b_q<b_{q-1}<\dots<b_1$
 is to have 
 \[ a_i = n-j+1-i\text{ for all $i \in [p]$}\quand q\in\{0,1\}\quand b_1= k+1=n+2.\]
Under these conditions, the end of the  chain is $v_{p+q}=[1^{j+p+q}\mid n+q]$.
Thus,
Theorem~\ref{lenart-transition} gives
$
 (1+\beta x_{n+1})  \fkG_{[1^{j} \mid n]} = \sum_{p=0}^{n-j} \sum_{q=0}^1 (-1)^p \beta^{p+q} \fkG_{[1^{j+p+q}\mid n+q]}$
which is equivalent to the desired formula.
\end{proof}

Fix a positive integer $k$. The \defn{$k$-Bruhat order} on $S_\infty$  is
transitive closure of the relation  with $v<_k w$ whenever $v \xrightarrow{(a,b)}w$ and $a\leq k <b$.
\begin{definition}[{cf. \cite[Def.~1.10]{LenartSottile}}]
\label{kpieri-def}
An \defn{unmarked $k$-Pieri chain} from $v\in S_\infty$ to $w \in S_\infty$ is a saturated chain in $k$-Bruhat order of the form 
\be\label{kpieri} v=v_0 \xrightarrow{(a_1,b_1)} v_1 \xrightarrow{(a_2,b_2)} \cdots \xrightarrow{(a_q,b_q)} v_q = w\ee
satisfying $1\leq a_i \leq k < b_i$ for all $i \in [q]$ along with the following two conditions:
\ben
\item[(P0)]  if $a_j =a_i > a_{i+1}$ for some $1 \leq j <i < q$ then $b_i>b_{i+1}$; and
\item[(P1)] $b_1\geq b_2\geq \dots \geq b_q$.
\een 
\end{definition}

An essential and non-obvious property of this definition is 
that for any permutations $v,w \in S_\infty$ at most one unmarked $k$-Pieri chain exists from $v$ to $w$
 \cite[Thm.~2.2]{LenartSottile}.
Following \cite{LenartSottile}, we write $v \xrightarrow{c(k)} w$ if such a chain exists. The reference \cite{LenartSottile}
gives an explicit algorithm to find the unique chain $v \xrightarrow{c(k)} w$, but we will not review this here.

Suppose $v,w \in S_\infty$  are such that $v \xrightarrow{c(k)} w$,
and assume that the unique unmarked $k$-Pieri chain from $v$ to $w$ has the form \eqref{kpieri}.
Define $\forced_k(v,w)$ to be the number of indices $i \in [q]$ such that either 
\bei
\item $b_1=b_2=\dots=b_i$ and $a_1>a_2>\dots>a_i$, or
\item $i<q$ and $b_i=b_{i+1}$ and $a_i>a_{i+1}$.
\eei
Then define $\prohib_k(v,w)$ to be the number of indices $i \in [q]$ such that 
\bei
\item $a_j=a_i$ for some $1\leq j <i$.
\eei
Notice that we have $\forced_k(v,w)\geq1$ if $v\neq w$ and that
\be\label{fp-obs}
\forced_k(v,w) + \prohib_k(v,w) \leq   \ell(w) - \ell(v)
\quand
k  \leq   \ell(w) - \ell(v) - \prohib_k(v,w).
\ee

\begin{theorem}[{Lenart--Sottile \cite[Cor.~1.16]{LenartSottile}}] 
\label{lensot-thm1}
If $0<p\leq k$ and $v \in S_\infty$ then
\[
\fkG_{[1^p\mid k]}  \fkG_{v} = \sum_{v\xrightarrow{c(k)}w}\tbinom{\ell(w) - \ell(v) - \forced_k(v,w) - \prohib_k(v,w)}{p - \forced_k(v,w)}\cdot  \beta^{\ell(w)-\ell(v) - p} \cdot   \fkG_w
\]
where the sum is over all permutations $w \in S_\infty$ such that $v\xrightarrow{c(k)} w$.
\end{theorem}

We note 
one more formula, which is a special case of the previous result:

\begin{theorem}[{Lenart \cite[Thm.~3.1]{lenart2000}}]\label{lenart-pieri}
Suppose $0<p\leq k$ and $\lambda$ is a partition with at most $k$ parts. Then 
    \[
\fkG_{[1^p\mid k]} \fkG_{[\lambda \mid k]}     = \sum_\mu \tbinom{\cols(\mu/\lambda)-1}{|\mu/\lambda|-p}
\cdot \beta^{|\mu/\lambda|-p}
\cdot \fkG_{[\mu\mid k]}
    \]
    where  the summation ranges over all partitions $\mu$ with $\ell(\mu)\leq k$ obtained by adding a vertical strip to $\lambda $ of weight at least $p$, and where $|\mu/\lambda|$ and $\cols(\mu/\lambda)$ respectively denote the numbers of cells and columns in   $\D_\mu \setminus \D_\lambda$.
\end{theorem}

We can now turn Theorem~\ref{qd-thm} into an exact, though not manifestly positive, formula for the coefficient function $\GC_z : S_\infty \to\NN$ from Section~\ref{GC-eq-sect}.

\begin{lemma}\label{prod-lem}
For all positive integers $k$ it holds that 
\[  \textstyle
 \prod_{i=1}^k (2+\beta x_i) =  2^k - \sum_{j=1}^k2^{k-j}  \cdot  (-\beta)^j  \cdot  \fkG_{[1^j \mid k]}.
   \]
\end{lemma}

\begin{proof}
Our argument is by induction on $k$. Since $\fkG_{[1^1\mid 1]} = x_1$, the base case when $k=1$ is immediate.    When $k>1$, we may assume by induction that 
\[ 
\ba
     \textstyle    \prod_{i=1}^k (2+\beta x_i) 
        &  
         = \textstyle(2+\beta x_k) \(2^{k-1} - \sum_{j=1}^{k-1}   2^{k-1-j} (-\beta)^j  \fkG_{[1^j \mid k - 1]} \)\\
        &  
         =\textstyle (2+\beta x_k) \(2^{k}  -  \sum_{j=0}^{k-1}  2^{k-1-j}(-\beta)^j  \fkG_{[1^j \mid k - 1]} \)
\ea
\]
where we interpret $[1^0\mid k-1]=1$. By Lemma~\ref{1gr-lem} 
we have 
\[
\ba
\textstyle (2+&\beta x_k)\textstyle\sum_{j=0}^{k-1}  2^{k-1-j}(-\beta)^j  \fkG_{[1^j \mid k - 1]} = 
\\ & \sum_{j=0}^{k-1} 2^{k-1-j}   \(    (-\beta)^j\fkG_{[1^j \mid k - 1]} + \sum_{i=j}^{k-1} (-\beta)^i    \fkG_{[1^i \mid k - 1]}  -  \sum_{i=j+1}^{k} (-\beta)^i \fkG_{[1^{i}\mid k]}\).
\ea
\]
The coefficient of $(-\beta)^j \fkG_{[1^j\mid k-1]}$ in the last expression is 
\[ 2^{k-1-j} + (2^{k-1} + 2^{k-2} +  \dots + 2^{k-1-j}) = 2^k\]
for each $j=0,1,\dots,k-1$, while the coefficient of $-(-\beta)^j \fkG_{[1^j\mid k]}$ is 
\[ 2^{k-1} + 2^{k-2} + \dots + 2^{k-1-(j-1)} = 2^k - 2^{k-j}\]
for each $j=1,2,\dots,k$. Thus $(2+\beta x_k)\sum_{j=0}^{k-1}  2^{k-1-j} (-\beta)^j  \fkG_{[1^j \mid k - 1]}$ is equal to
\[ \textstyle   2^k + 2^k \sum_{j=1}^{k-1} (-\beta)^j \fkG_{[1^j\mid k-1]} - \sum_{j=1}^k  (2^k-2^{k-j}) (-\beta)^j\fkG_{[1^j\mid k]}.
\]
Lemma~\ref{1gr-lem} also
  shows that  $(2+\beta x_k)   2^{k}=(2+\beta x_k)   2^{k}  \fkG_{[1^0 \mid k - 1]}$ is equal to
\[  \textstyle
2^{k+1} + 2^k \sum_{j=1}^{k-1} (-\beta)^j  \fkG_{[1^j \mid k - 1]} - 2^k \sum_{j=1}^k (-\beta)^j \fkG_{[1^j \mid k]}.
\]
These expressions have difference $ \ds \prod_{i=1}^k (2+\beta x_i) =  2^k - \sum_{j=1}^k2^{k-j}  (-\beta)^j  \fkG_{[1^j \mid k]}.$
\end{proof}

For $v,w \in S_\infty$, recall the definitions of $\forced_k(v,w)$ and $\prohib_k(v,w)$ from Theorem~\ref{lensot-thm1}.
Set $\epsilon_k(v,v) = 1$ and $\rho_k(v,v)=2^k$, and for $v \xrightarrow{c(k)} w\neq v$ define
\be
\epsilon_k(v,w) := (-1)^{1+\forced_k(v,w)}
\quand
\rho_k(v,w) := 2^{k + \ell(v) - \ell(w) + \prohib_k(v,w)}.
\ee
Set $\epsilon_k(v,w)= \rho_k(v,w)=0$ when we do not have $v \xrightarrow{c(k)} w$. 

\begin{theorem}\label{dom-thm}
Suppose $z \in \Ivex$ is quasi-dominant   and $k=k(z)$.  Then
\[\textstyle \GC_z(w) = \sum_{v \in \cB(z)} 
\epsilon_k(v,w)   \rho_k(v,w) \geq 0\quad\text{for all $w \in S_\infty$}.\]
\end{theorem}

\begin{proof}
Recall from Section~\ref{GC-eq-sect} that we always have
 $\GC_z(w)\geq 0$.
To get an given expression for this coefficient, 
we start with the formula    $ \fkGO_z = \iG_z \prod_{i=1}^{k} (2 + \beta x_i)$ from Theorem~\ref{qd-thm}. Substituting in Lemma~\ref{prod-lem} and   applying Theorem~\ref{lensot-thm1} turns this to
\[
 \fkGO_z = \sum_{v \in \cB(z)} 
\(2^{k}\cdot  \beta^{\ell(v)-\ellhat(z)} \cdot \fkG_v 
-
  \sum_{\substack{w \in S_\infty \\ v\xrightarrow{c(k)}w}} \Sigma_{v,w}^k \cdot \beta^{\ell(w)-\ellhat(z)} \cdot \fkG_w\)
 \]
 where for $v \xrightarrow{c(k) }w$  we define
\[\textstyle\Sigma_{v,w}^k :=
\sum_{p=1}^k
 (-1)^p 2^{k-p} \tbinom{\ell(w) - \ell(v) - \forced_k(v,w) - \prohib_k(v,w)}{p - \forced_k(v,w)}.
 \]
 If $v=w$ then $\ell(w) - \ell(v) =\forced_k(v,w) = \prohib_k(v,w) = 0$   so $\epsilon_k(v,w) \rho_k(v,w) = 2^k$ and $\Sigma_{v,w}^k =0$. If $v\xrightarrow{c(k)}w$ but $v\neq w$ then we must have $\forced_k(v,w)\geq 1$, so in view of  \eqref{fp-obs} we can rewrite $\Sigma_{v,w}^k$ as 
 \[\textstyle\Sigma_{v,w}^k 
 = (-1)^{\forced_k(v,w)} 2^{k-\forced_k(v,w)} \sum_{m \in \NN} (-1/2)^m \tbinom{\ell(w) - \ell(v) - \forced_k(v,w) - \prohib_k(v,w)}{m}
 \]
 which transforms via the binomial theorem to 
 $ \Sigma_{v,w}^k  = - \epsilon_k(v,w) \rho_k(v,w)$. 
 
 Using these identities, our expression above for $\fkGO_z$ becomes  the desired formula
$
 \fkGO_z = \sum_{v \in \cB(z)} \sum_{w \in S_\infty} \beta^{\ell(w)-\ellhat(z)} \cdot \epsilon_k(v,w) \rho_k(v,w) \cdot  \fkG_w.
$
\end{proof}

\subsection{Examples and conjectures}\label{ex-conj-sect}

Fix $ z \in I_\infty$. If $z\neq 1$ then define $k=k(z)$ as in \eqref{k-eq} and
 let 
 \be
 j=j(z) := \max\{ j \in \PP : z(i) = i\text{ for all }1\leq i\leq j\}.
 \ee 
 If $z=1$ is  the identity permutation 
 then set $j=j(z)=1$ and $k=k(z)=0$.
 
 \begin{definition}\label{zpieri-def}
For permutations $v,w \in S_\infty$ we write $v \xrightarrow{[z]}w$ if there exists an unmarked $k$-Pieri chain
\[ v=v_0 \xrightarrow{(a_1,b_1)} v_1 \xrightarrow{(a_2,b_2)} \cdots \xrightarrow{(a_q,b_q)} v_q = w\]
as in Definition~\ref{kpieri-def}, such that for each $i\in[q]$ it additionally holds that 
\[ j \leq a_i \leq k < b_i \qquand \text{$a_i<z(a_i)$ or $z(b_i)<b_i$} .\]
Finally, define $\cB^+(z) := \left\{ w \in S_\infty : v\xrightarrow{[z]} w \text{ for some }v \in \cB(z)\right\}$.
\end{definition}

Recall that the elements of $\cB(z)$ are referred to as  Hecke atoms for $z$.
Following this convention, we call $\cB^+(z)$ the set of \defn{extended Hecke atoms} for $z$.

\begin{example}
If $z=1 \in I_\infty$ then $j(z) = 1>0=k(z)$ so only empty chains satisfy the conditions in Definition~\ref{zpieri-def}.
Therefore  $\cB^+(1) = \cB(1) = \{1\}$.
\end{example}


Based on computations and the preceding theorem, the set of extended Hecke atoms $\cB^+(z)$ seems to be a good approximation for $\supp(\GC_z)$. In particular:



\begin{corollary}\label{qd-supp-cor}
If $z \in \Ivex$ is quasi-dominant and $k=k(z)$ then
\[
\supp(\GC_z) \subseteq \cB^+(z)
=\left\{ w \in S_\infty : v\xrightarrow{c(k)} w \text{ for some }v \in \cB(z)\right\}.
\]
\end{corollary}

 \begin{proof}
 In this case $j(z)=1$ and $k(z)=k$, so $a<z(a)$ for all $j\leq a \leq k$.
 Thus all unmarked $k$-Pieri chains from $v$ to $w$ satisfy the conditions in Definition~\ref{zpieri-def}, 
 \end{proof}

If $a$ and $b$ are integers then let $[a,b] = \{ i \in \ZZ : a\leq i \leq b\}$.
For $w \in S_\infty$ let 
\be\supp(w) := \{ i \in \PP : w(i) \neq i\} .\ee

\begin{proposition} \label{supp-prop}
If $z \in I_\infty$ has $\supp(z)  \subseteq [a,b]$ for integers $a\leq b$,
then \[\supp(v) \subseteq [a,b]\quand \supp(w) \subseteq [a-1,b+1]\]
for all $v \in \cB(z)$ and $w \in \cB^+(z)$.
\end{proposition}

\begin{proof}
If $z=1$ then $\supp(z) = \supp(v) = \supp(w) = \varnothing$
for all $v \in \cB(z)$ and $w \in \cB^+(z)$
so the result holds trivially. Assume $z\neq 1$.

Let $j=j(z)$ and $k = k(z)$ and suppose $v \in \cB(z)$.
Since $\alpha_\inv(z)$ fixes all integers $i \notin[a,b]$, 
it follows from \eqref{sim-eq} that we have $\supp(v) \subseteq [a,b]$.
Therefore, if $c \leq k$ and $ b+1< d$ then we must have $\ell(v(c,d)) \neq \ell(v)+1$,
since if $v(c) < v(d)$ then it also holds that $v(c) <  b+1 = v(b+1) < d = v(d)$.

Hence,
if there exists any unmarked $k$-Pieri chain of the form \eqref{kpieri} from $v$ to $w \in S_\infty$, then 
this chain must have $b+1 \geq b_1 \geq b_2 \geq \dots \geq b_q>k$. Since to have $v \xrightarrow{[z]} w$ we also require that $j \leq a_i$ for all $i$, 
we conclude that $v \xrightarrow{[z]} w$ implies $\supp(w) \subseteq[j,b+1]$.
The result follows since $a-1 \leq j$ as $\supp(z) \subseteq [a,b]$.
\end{proof}

We have used a computer to verify the following for all vexillary $z \in I_{11}$:

\begin{conjecture}\label{b+conj}
If $z \in \Ivex$ then $\cB(z)\subseteq \supp(\GC_z) \subseteq \cB^+(z)$.
\end{conjecture}

Both containments in this conjecture can be strict, 
but in some notable cases we have $\supp(\GC_z) = \cB^+(z)$.
 Among the dominant (respectively, vexillary) involutions in $S_8$,
 the equality $\supp(\GC_z) = \cB^+(z)$ holds in 67 of 70 cases (respectively, 179 of 323 cases).
 

 \begin{example}\label{t-ex}
 Suppose $z=t_n:=(1, n)\in I_n$.
 Then $z$ is dominant of shape $\lambda=(n-1,1^{n-2})$ and $ \cB(z)$ consists of the permutations in $S_n$
 whose {inverses} in one-line notation are the shuffles of $n1$ and $234\cdots(n-1)$
 with at most one letter between $n$ and $1$.
 The larger set $\cB^+(z)$ consists of the permutations in $S_{n+1}$
 whose {inverses} in one-line notation are the shuffles of the words $n1$ and $234\cdots(n-1)(n+1)$
 with at most two letters between $n$ and $1$, excluding the inverse of $234\cdots(n-1)(n+1)n1$.
We shall see below that $\supp(\GC_z) = \cB^+(z)$.
 \end{example}

We can generalize this example.
Let $\Sh(n)$ denote the set of words obtained by shuffling
the 2-letter word $n2$ with then $(n-1)$-letter word
\[1 \hs\widehat 2\hs 3\cdots (n-1)\hs\widehat n \hs(n+1)= 1345\cdots (n-3)(n-2)(n-1) (n+1).\]
Define $\cX(n)$ to be the subset of words in $\Sh(n)$ for which
\begin{itemize}
    \item the first letter is $1$ while the last letter is $n+1$; and
\item at most one letter appears between $n$ and $2$.
\end{itemize}
Define $\cY(n)\subset\Sh(n)$ be the set of words with exactly one or exactly two letters between $n$ and $2$.
Note that $\cX(n)$ and $\cY(n)$ have nonempty intersection.

\begin{example}
We have 
$\cX(3) = \left\{ 1\overline{3}\underline{2}4 \right\}
$, $ \cY(3) = 
\left\{
\overline{3}1\underline{2}4, \ \  1\overline{3}4\underline{2}, \ \ \overline{3}14\underline{2}\right\}
$, and 
\[\Sh(3)=\left\{{\barr{ccccccccc} \overline{3}\underline{2}14, &&1\overline{3}\underline{2}4, && 14\overline{3}\underline{2}\\
&\overline{3}1\underline{2}4, && 1\overline{3}4\underline{2},  \\
&&\overline{3}14\underline{2}
\earr}\right\}.\]
\end{example}

Define 
$c_4=(1,4,5,2) = 41352 \in S_5$ and for $n\geq 5$ let \[ c_n=(1,n,n-1,n-2,\dots,7,6,5,2) = n1342567\cdots (n-2)(n-1) \in S_n.\] 
Finally set $\cZ(n) = \varnothing$ for $n\leq 3$ and $\cZ(n) = \{c_n\}$ for $n\geq 4$.

\begin{proposition}\label{t-prop}
Let $n>2$ and $z=(2,n) \in \Ivex$. Then the following holds:
   \ben
    \item[(a)] $\cB(z)= \left\{w^{-1} : w \in \cX(n)\right\}$.
    
    \item[(b)] $\cB^+(z)=\left\{w^{-1} : w \in \cX(n)\cup \cY(n)\sqcup \cZ(n)\right\}$.

\item[(c)] The orthogonal Grothendieck polynomial of $z$ is
\[\fkGO_{z} = \sum_{w^{-1} \in \cX(n)} 2 \cdot \beta^{w(2)-w(n)-1}\cdot  \fkG_w +  \sum_{w^{-1} \in \cY(n)} \beta^{w(2)-w(n)-1}\cdot \fkG_w
\]
 and consequently
 $\supp(\GC_z) =\left\{w^{-1} : w \in \cX(n)\cup \cY(n)\right\}$.

\een
         \end{proposition}

By the results in Section~\ref{shift-sect},  this proposition 
determines $\cB(z)$, $\cB^+(z)$, and $\fkGO_z$ for all transpositions $z=(i,j) \in I_\infty$.
For these elements,
 Conjecture~\ref{b+conj} holds but $ \cB^+(z)\setminus \supp(\GC_z)$ contains one extra element if $1<i<i+1<j$.

\begin{proof}
Part (a) is straightforward to derive from \eqref{sim-eq}, and part (b) then follows as an exercise from the definitions. We omit the details.

  Let $\Sh_1(n)$ denote the set of shuffles of
 $n1$ and $ 2345\cdots (n+1)$, viewed as permutations in one-line notation. Then
define $\cX_1(n)$ to be the subset of words in $\Sh_1(n)$ for which
the last letter is $n+1$ and at most one letter appears between $n$ and $1$.
Finally let $\cY_1(n)$ be the subset of words in $\Sh_1(n)$ for which exactly one or exactly two letters appear between $n$ and $1$.

Now set $y = (1,n) \in I_\infty$. 
We shall first use Theorem~\ref{dom-thm} to prove that
\be\label{t1n-eq}
        \fkGO_{y}
        = \sum_{w^{-1}\in \cX_1(n)}2\cdot \beta^{w(1)-w(n)-1}\cdot \fkG_w + \sum_{w^{-1}\in \cY_1(n)}\beta^{w(1)-w(n)-1}\cdot \fkG_w.
        \ee     
          
It follows from \eqref{sim-eq} that $\cB(y) = \{ w^{-1} : w \in \cX_1(n)\}$.
Therefore, if $v \in \cB(y)$ and $1<i_1<i_2<i_3<\dots$ are the indices $i$ with $v(1) < v(i)$, then we must have $v(1)<v(i_1)<v(i_2)<v(i_3)<\cdots$
and also $v(i_1) = v(1)+1$. 

Set $j := v(1)$
and $k:=k(y)=1$. Given the previous paragraph, we see that if $v \in \cB(y)$ then a chain $v\xrightarrow{c(k)}w $ 
 is either trivial with $v=w$ or a single Bruhat cover of the form $v\xrightarrow{(1,b)}w$ where $w=v(1,b) = s_j v$.
 For chains of the first type we have $\varepsilon_k(v,w) = 1$ and $\rho_k(v,w) = 2$,
 and $w^{-1}=v^{-1}$ gives an arbitrary element of $\cX_1(n)$.
For chains of the second type we have $\varepsilon_k(v,w)=1$ and $\rho_k(v,w)=1$,
and the permutations $w$ range over exactly the set of inverses of the elements of $\cY_1(n)$.

The only thing left to show to deduce \eqref{t1n-eq} from Theorem~\ref{dom-thm}
is that if $w \in S_\infty$ has $w^{-1} \in \cX_1(n)\cup \cY_1(n)$ then $\ell(w) - \ellhat(y) = w(1)-w(n)-1$.
 This holds since 
 in such a permutation $w$
 it must hold that $w(n)<w(n+1)$, so $w^{-1}$ has $w(1)-1$ inversions involving $1$, along with $n-1 - (w(n)-1)$ inversions involving $n$,
 and a single inversion involving both $1$ and $n$, which gives
 \[
 \ell(w) = \ell(w^{-1}) = w(1) -1 + n-1 - (w(n)-1) - 1 = w(1) - w(n) + n-2.
 \]
 As we have $\ellhat(y) = n-1$ by \eqref{ellhat-eq}, the identity \eqref{t1n-eq}  follows. 

Let $z = (2,n) \in I_\infty$.
To prove part (c) we apply $\bpartial_1$ to both sides of \eqref{t1n-eq}.
On the left side we get $\bpartial_1\fkGO_y = \fkGO_z$ by \eqref{orthogonal-recursion}.
To evaluate $\bpartial_1$ applied to the right side of \eqref{t1n-eq}, fix
$w \in S_\infty$ with $w^{-1} \in \cX_1(n)\cup \cY_1(n)$. Then observe that
   $w(1) < w(2)$ if and only if $w^{-1} = n1234\cdots(n-1) \in \cX_1(n)$,
and in this case  it holds that $(ws_1)^{-1}  = n2134\cdots (n-1) \in \cX_1(n)$ and
\[ 
\bpartial_1 \( \beta^{ws_1(1) - w(n)-1} \fkG_{ws_1} + \beta^{w(1) - w(n)-1} \fkG_w\)
=
\bpartial_1\( \beta \fkG_{ws_1} +  \fkG_w\) = 0.
\]
The map $w \mapsto s_1 w$ defines a bijection $\cY_1(n)\to \cY(n)$ as well as a bijection 
\[ \cX_1(n) \setminus \{n1234\cdots(n-1), \ n2134\cdots (n-1)\} \to \cX(n),\]
so it follows from \eqref{Grothendieck-recursion} that $\bpartial_1$ transforms   \eqref{t1n-eq}
to the formula in part (c).
\end{proof}


\begin{remark}
We often draw $\cB^+(z)$ as an induced subgraph of the Hasse diagram of the \defn{left weak order} on $S_\infty$.
In other words, we give $\cB^+(z)$ the directed graph structure that has an edge $v \to w$ if and only if 
$v$ and $w$ are extended Hecke atoms of $z$ such that
$\ell(w) = \ell(v)+1$
and $ w=s_i v$ for some $ i\in \PP$.
Figures~\ref{t1n-fig}, \ref{gn-fig}, and \ref{w0-fig} show some instances of these graphs.

In all examples that we have computed, the directed graph structure on $\cB^+(z)$ just described is connected. This is not the case if we use the right weak order instead of the left weak order. It would be interesting to know if this connectedness property is a general phenomenon.
\end{remark}

\begin{figure}[h]
  \centerline{
      \begin{tikzpicture}[xscale=2, yscale=1,>=latex,baseline=(z.base)]
      \node at (0,0.0) (z) {};
        \node[draw] at (0.5,0) (11) {\tiny $4231^{-1}:1$};
        \node[draw] at (1.5,0) (12) {\tiny $24351^{-1}:1$};
        \node[draw,fill=lightskyblue] at (0,1) (21) {\tiny $4213^{-1}:3$};
        \node[draw,fill=lightskyblue] at (1,1) (22) {\tiny $2431^{-1}:3$};
        \node[draw] at (2,1) (23) {\tiny $23451^{-1}:1$};
        \node[draw,fill=lightskyblue] at (0,2) (31) {\tiny $4123^{-1}:2$};
        \node[draw,fill=lightskyblue] at (1,2) (32) {\tiny $2413^{-1}:2$};
        \node[draw,fill=lightskyblue] at (2,2) (33) {\tiny $2341^{-1}:2$};
        \draw[->,thick]  (31) -- (21);
        \draw[->,thick]  (32) -- (21);
        \draw[->,thick]  (32) -- (22);
        \draw[->,thick]  (33) -- (22);
        \draw[->,thick]  (33) -- (23);
        \draw[->,thick]  (21) -- (11);
        \draw[->,thick]  (22) -- (11);
        \draw[->,thick]  (22) -- (12);
        \draw[->,thick]  (23) -- (12);
       \end{tikzpicture}
  }
  \caption{The directed graph $ \cB^+(z)$ for $z=t_{4}=(1,4)$.
  The data in each box is $\boxed{w:\GC_z(w)}$ with $w$ in inverse one-line notation.
  The blue vertices are the elements of $\cB(z)$.
One has $\supp(\GC_z) = \cB^+(z)$ in this case.
  }\label{t1n-fig}
  \end{figure}


 \begin{example}\label{g-ex}
Suppose $z=g_n:=  (1, n+1)(2, n+2)\cdots(n, 2n) \in I_{2n}$. Then $z$ is dominant of shape $\lambda=(n^n)$ and 
 the set $ \cB(z)$ consists of the single element whose {inverse} is 
 $
 (n+1)1(n+2)2\cdots  (2n)n\in S_{2n}.$
 One can show that the larger set 
 $\cB^+(z)$ consists of the $2^n$ permutations   whose {inverses} have the form 
 $ (n+1)a_1b_1a_2b_2\cdots a_nb_n \in S_{2n+1}$
where $\{a_i,b_i\} = \{ i, n+1+i\}$ for $i \in [n]$.
 \end{example}
 
\begin{figure}[h]
  \centerline{
     \begin{tikzpicture}[xscale=2.5, yscale=1.2,>=latex,baseline=(z.base)]
      \node at (0,0.0) (z) {};
      \node[draw] at (1,0) (00) {\tiny $4516273^{-1}:1$};
      \node[draw] at (0,1) (21) {\tiny $4156273^{-1}:2$};
      \node[draw] at (1,1) (22) {\tiny $451623^{-1}:2$};
      \node[draw] at (2,1) (23) {\tiny $4512673^{-1}:2$};
      \node[draw] at (0,2) (31) {\tiny $415623^{-1}:4$};
      \node[draw] at (1,2) (32) {\tiny $4152673^{-1}:4$};
      \node[draw] at (2,2) (33) {\tiny $451263^{-1}:4$};
      \node[draw, fill=lightskyblue] at (1,3) (41) {\tiny $415263^{-1}:8$};
        \draw[->,thick]  (21) -- (00);
        \draw[->,thick]  (22) -- (00);
        \draw[->,thick]  (23) -- (00);
        \draw[->,thick]  (31) -- (21);
        \draw[->,thick]  (31) -- (22);
        \draw[->,thick]  (32) -- (21);
        \draw[->,thick]  (32) -- (23);
        \draw[->,thick]  (33) -- (22);
        \draw[->,thick]  (33) -- (23);
        \draw[->,thick]   (41) -- (31);
        \draw[->,thick]   (41) -- (32);
        \draw[->,thick]   (41) -- (33);
       \end{tikzpicture}}

  \caption{The directed graph $ \cB^+(g_3)$, presented as in Figure~\ref{t1n-fig}.
  }\label{gn-fig}
  \end{figure}

 The previous example has an extension analogous to  Proposition~\ref{t-prop},
which shows that $\supp(\GC_z) = \cB^+(z)$ whenever $z=g_n$; see Proposition~\ref{g-prop}.


\begin{example}
Finally let $w_0 = n\cdots 321\in I_n$ be the longest element of $S_n$. Then $w_0$
is dominant of shape $\lambda=(n-1,\dots,3,2,1)$.
The set $\cB(w_0)$ does not seem to have any description simpler than \eqref{sim-eq}.
One can check that the larger set
$ \cB^+(w_0)$ contains the permutations in $S_{n+1}$ whose {inverses} in one-line notation have the form 
$
u_1\cdots u_i(n+1) u_{i+1}\cdots u_n 
$
where $u=u_1u_2\dots u_n $ is the {inverse} of an element of $ \cB(w_0)$ and $i \in [n]$ has $n \geq 2 u_j$ for each $i<j\leq n$. 
\end{example}

\begin{figure}[h]
  \centerline{
\begin{tikzpicture}[xscale=1.8, yscale=1,>=latex,baseline=(z.base)]
 \node at (0,0.0) (z) {};
 \node[draw] at (2,0) (00) {\tiny $43521^{-1}:1$};
 \node[draw] at  (1,1) (11)   {\tiny $43251^{-1}:3$};
 \node[draw] at  (2,1) (12)   {\tiny $34521^{-1}:1$};
 \node[draw] at  (3,1) (13)   {\tiny $43512^{-1}:1$};
 \node[draw] at (0,2) (21) {\tiny $42351^{-1}:2$};
 \node[draw] at (1,2) (22) {\tiny $34251^{-1}:3$};
 \node[draw, fill=lightskyblue] at (2,2) (23) {\tiny $4321^{-1}:6$};
 \node[draw] at (3,2) (24) {\tiny $34512^{-1}:1$};
 \node[draw] at (4,2) (25) {\tiny $43152^{-1}:3$};
 \node[draw,fill=lightskyblue] at (0,3) (31) {\tiny $4231^{-1}:4$};
 \node[draw] at (1,3) (32) {\tiny $32451^{-1}:2$};
 \node[draw,fill=lightskyblue] at (2,3) (33) {\tiny $3421^{-1}:6$};
 \node[draw] at (3,3) (34) {\tiny $34152^{-1}:2$};
 \node[draw,fill=lightskyblue] at (4,3) (35) {\tiny $4312^{-1}:6$};
 \node[draw] at (5,3) (36) {\tiny $41352^{-1}:2$};
 \node[draw, fill=lightskyblue] at (1.5,4) (41) {\tiny $3241^{-1}:4$};
 \node[draw, fill=lightskyblue] at (3,4) (42) {\tiny $3412^{-1}:4$};
 \node[draw, fill=lightskyblue] at (4.5,4) (43) {\tiny $4132^{-1}:4$};
   \draw[->,thick]  (11) -- (00);
   \draw[->,thick]  (12) -- (00);
   \draw[->,thick]  (13) -- (00);
   \draw[->,thick]  (21) -- (11);
   \draw[->,thick]  (22) -- (11);
   \draw[->,thick]  (22) -- (12);
   \draw[->,thick]  (23) -- (11);
   \draw[->,thick]  (24) -- (12);
   \draw[->,thick]  (24) -- (13);
   \draw[->,thick]  (25) -- (13);
   \draw[->,thick]  (31) -- (21);
   \draw[->,thick]  (31) -- (23);
   \draw[->,thick]  (32) -- (22);
   \draw[->,thick]  (33) -- (22);
   \draw[->,thick]  (33) -- (23);
   \draw[->,thick]  (34) -- (24);
   \draw[->,thick]  (34) -- (25);
   \draw[->,thick]  (35) -- (23);
   \draw[->,thick]  (35) -- (25);
   \draw[->,thick]  (36) -- (25);
   \draw[->,thick]   (41) -- (32);
   \draw[->,thick]   (41) -- (33);
   \draw[->,thick]   (42) -- (33);
   \draw[->,thick]   (42) -- (34);
   \draw[->,thick]   (42) -- (35);
   \draw[->,thick]   (43) -- (35);
   \draw[->,thick]   (43) -- (36);
  \end{tikzpicture}}
  \caption{The directed graph $ \cB^+(4321)$, presented as in Figure~\ref{t1n-fig}.
  }\label{w0-fig}
  \end{figure}

Extending this example,
let 
\be w_{ij} := (i, j)(i+1,j-1)(i+2,j-2)\cdots (i+k,j-k) \in \Ivex\ee for any integers $1\leq i<j$ where $k=\lfloor\tfrac{j-i-1}{2}\rfloor$. Computer calculations support the following conjecture.
The case when $i=1$ has been checked for all
$j\leq 11$.

 \begin{conjecture} 
One has
$\supp(\GC_{w_{ij}}) = \cB^+(w_{ij})$ if $i=1$ or $j-i$ is odd.
\end{conjecture}

In the examples we have computed 
when $i>1$ and $j-i$ is even, the set $\supp(\GC_{w_{ij}})$ can be a proper subset of $\cB^+(w_{ij})$ with much smaller cardinality.
We do not have any specific conjecture for $\supp(\GC_{w_{ij}})$ in these cases
(apart from it being a subset of  $\cB^+(w_{ij})$, as asserted in Conjecture~\ref{b+conj}),
or for the values of $\GC_{w_{ij}}$.
This set of coefficients can be quite complicated,
although its maximum seems to be a central binomial coefficient; see Figure~\ref{values-fig}.
 
 \begin{figure}[h]
 \centerline{\begin{tabular}{| l | l|}
 \hline&\\[-8pt] 
 $n$ & nonzero values of $\GC_{w_0}$ for $w_0=n\cdots 321\in I_n$ \\[-8pt]&\\ \hline & \\[-8pt]
 1  & 1 \\[-8pt]&\\
2  & 1, 2 \\[-8pt]&\\
3  & 1, 2, 3 \\[-8pt]&\\
4  & 1, 2, 3, 4, 6 \\[-8pt]&\\
5  & 1, 2, 3, 4, 6, 9, 10 \\[-8pt]&\\
6  & 1, 2, 3, 4, 6, 8, 9, 10, 12, 18, 20 \\[-8pt]&\\
7  & 1, 2, 3, 4, 5, 6, 8, 9, 10, 12, 13, 14, 15, 18, 20, 27, 30, 33, 34, 35 \\[-8pt]&\\
8  & 1, 2, 3, 4, 5, 6, 8, 9, 10, 12, 13, 14, 15, 16, 18, 20, 24, 27, 30, 33, \\  & 34, 35, 36, 40, 54, 60, 66, 68, 70
\\[-8pt]&\\ \hline
\end{tabular}}
\caption{Nonzero values of $\GC_{w_0}$ for $w_0=n\cdots 321\in I_n$.}\label{values-fig}
\end{figure}



\section{Expansions for vexillary involutions}\label{vex-sect}

The first main result of this section is a generalization of Theorem~\ref{qd-thm}
that holds any vexillary involution; see Theorem~\ref{ivex-thm}. This formula is more complicated than in the quasi-dominant case,
and we use Section~\ref{setting-sect} to explain the relevant notation.
The rest of the section derives some interesting consequences.

\subsection{Setting up the main theorem}\label{setting-sect}

Choose a vexillary involution $z \in \Ivex$.
Since a permutation is vexillary if and only if it is $2143$-avoiding, 
it follows by definition that $z$
has cycle notation \[z = (a_1, b_1)(a_2, b_2)\cdots (a_q, b_q)\]
where $1\leq a_1 <a_2 < \dots <a_q < \min\{b_1,b_2,\dots,b_q\}$. 
We refer to the numbers $a_i$ as \defn{left endpoints} and to the numbers $b_i$ as \defn{right endpoints}.

The \defn{left segments} of $z$ are the maximal subsets of consecutive left endpoints,
that is, the equivalence classes in $\{a_1,a_2,\dots,a_q\}$ under the transitive closure of the relation with $a_i \sim a_j$ if $|a_j-a_i|\leq 1$.
There is at most one left segment containing $1$, which we refer to as the \defn{immobile segment}.
All other left segments are \defn{mobile}.

Suppose $L$ is a mobile left segment of $z$ and define $c_0: = \min(L) - 1$. Notice that we must have $c_0 = z(c_0) \in \PP$.
Now, for any subset $  S \subseteq \{a_1,a_2,\dots, a_q\}$ define $\sigma_{S,L} \in S_\infty$ to be the cyclic permutation
\be \sigma_{S,L}: = (c_0,c_1,c_2,\dots,c_k) \quad\text{where }S\cap L = \{c_1<c_2<\dots<c_k\}.\ee
This is the identity element when $S\cap L$ is empty. 
The conjugated permutation 
\[  (\sigma_{S,L})^{-1} \cdot z \cdot \sigma_{S,L} \] 
 is not always vexillary, but
its arc diagram is obtained from $z$ by simply replacing each edge $\{c_i < z(c_i)\}$ for $i \in [k]$ with $\{ c_{i-1} < z(c_i)\}$.
Also define \be \textstyle \sigma_{S}: = \prod_L \sigma_{S,L}\ee where the product is over all mobile left segments of $z$.
This product can be taken in any order since as $L$ varies the permutations $\sigma_{S,L}$ have disjoint supports.

\begin{example} \label{sigma-ex}
Suppose our vexillary involution is

\[z = (2, 7)(3, 8)(4, 6)(5, 9)
=
\arcstart
{
*{\bullet}     & *{\bullet}  \arc{1.2}{rrrrr}   & *{\bullet} \arc{1.2}{rrrrr}  & *{\bullet}  \arc{0.6}{rr}   & *{\bullet} \arc{1}{rrrr}  & *{\bullet}    & *{\bullet}  & *{\bullet}  & *{\bullet}  \\
1 & 2 & 3 & 4 & 5  & 6 & 7 & 8 & 9      
} 
\arcstop.
\]
This permutation has a unique left segment $L=\{2,3,4,5\}$, which is mobile.
For the subset $S = \{2,4,5\}$ we have $\sigma_{S,L} = (1,2,4,5)$ and 

\[  (\sigma_{S,L})^{-1} \cdot z\cdot \sigma_{S,L} =
(1,7)(2,6)(3,8)(4,9)
= \arcstart
{
*{\color{red}\bullet}     \arc{1.4}{rrrrrr}   & *{\color{red}\bullet}  \arc{1}{rrrr}   & *{\bullet} \arc{1.2}{rrrrr}  & *{\color{red}\bullet}  \arc{1.2}{rrrrr}   & *{\color{red}\bullet}  & *{\bullet}    & *{\bullet}  & *{\bullet}  & *{\bullet}  \\
1 & 2 & 3 & 4 & 5  & 6 & 7 & 8 & 9      
} 
\arcstop.
\]
\end{example}

Suppose $a_i$ and $a_j$  are left endpoints of $z$ in the same left segment with $i<j$.
 We say that $a_j$ is a \defn{crossing bound} of $a_i$ if $\{i\} = \{ t : i\leq t <j\text{ and }b_t<b_j\}$.
 Let $\crb_z(a_i)$ denote the set of all crossing bounds of $a_i$.
While this set may have multiple elements,
notice that a given left endpoint $a_j$ can belong to $\crb_z(a_i)$ for at most one choice of $a_i$.

\begin{definition}
A set $  S \subseteq \{a_1,a_2,\dots, a_q\}$   is \defn{shiftable} relative to $z$ if both
\bei
\item[(1)] no element of $S$ is in a left segment of $z$ containing $1$, and
\item[(2)] if some $a_i \notin S$ then $S$ does not contain any crossing bound of $a_i$.
\eei
Let $\shiftable(z)$ be the set of all such subsets, and for each $S \in \shiftable(z)$
define  
\[\iG_{z,S} := \iG_v\quad\text{where }v = (\sigma_S)^{-1} \cdot z\cdot \sigma_S.\]
\end{definition}
The empty set is always shiftable, as is union of all mobile left segments of $z$.
A singleton set $S = \{b\}$ is in $\shiftable(z)$
if and only if  $b$ belongs to a mobile left segment $L$ of $z$ with the property that every $a \in L$ with $a<b$ has   $z(a)>z(b)$.

\begin{example}\label{ivex-ex}
When $z=(2, 7)(3, 8)(4, 6)(5, 9)$ there are 9 shiftable subsets of left endpoints,
given by $\emptyset$, $\{2\}$, $\{4\}$, $\{2,3\}$, $\{2,4\}$, $\{4,5\}$,
$\{2,3,4\}$, $\{2,4,5\}$, and $\{2,3,4,5\}$. The polynomial $\fkGO_z $ can be expressed as a sum of 9 terms
\[
\ba \fkGO_z &= 
   (2 + \beta x_2)(2 + \beta x_3)(2 + \beta x_4)(2 + \beta x_5) \cdot \iG_{(2, 7)(3, 8)(4, 6)(5, 9)}  \\&\quad +
   \beta\cdot (1 + \beta x_2)(2 + \beta x_3)(2 + \beta x_4)(2 + \beta x_5)\cdot \iG_{(1, 7)(3, 8)(4, 6)(5, 9)}   \\&\quad +
   \beta\cdot(2 + \beta x_2)(2 + \beta x_3)(1 + \beta x_4)(2 + \beta x_5)\cdot \iG_{(1, 6)(2, 7)(3, 8)(5, 9)} \\&\quad -
    \beta^2\cdot(1 + \beta x_3)(2 + \beta x_4)(2 + \beta x_5)\cdot \iG_{(1, 7)(2, 8)(4, 6)(5, 9)} \\&\quad +
   \beta^2\cdot(1 + \beta x_2)(2 + \beta x_3)(1 + \beta x_4)(2 + \beta x_5)\cdot \iG_{(1, 7)(2, 6)(3, 8)(5, 9)} \\&\quad -
   \beta^2\cdot(2 + \beta x_2)(2 + \beta x_3)(1 + \beta x_5)\cdot \iG_{(1, 6)(2, 7)(3, 8)(4, 9)} \\&\quad -
   \beta^3\cdot(1 + \beta x_3)(1 + \beta x_4)(2 + \beta x_5)\cdot \iG_{(1, 7)(2, 8)(3, 6)(5, 9)} \\&\quad -
   \beta^3\cdot(1 + \beta x_2)(2 + \beta x_3)(1 + \beta x_5)\cdot \iG_{(1, 7)(2, 6)(3, 8)(4, 9)}        \\&\quad +
   \beta^4\cdot(1 + \beta x_3)(1 + \beta x_5)\cdot\iG_{(1, 7)(2, 8)(3, 6)(4, 9)}.
\ea
\]
\end{example}

Each term in the expression for $\fkGO_{(2, 7)(3, 8)(4, 6)(5, 9)}$ just given has the form 
\[\beta^{|S|} \cdot \customvarpi_{z,S}  \cdot \iG_{z,S}
\]
where $S$ ranges over all shiftable subsets for $z=(2, 7)(3, 8)(4, 6)(5, 9)$ and $\customvarpi_{z,S} $ is the following polynomial.
Given  any subset $  S \subseteq \{a_1,a_2,\dots, a_q\}$ define    
\be
\label{varpi-eq}
\textstyle\customvarpi_{z,S} := \prod_{i \in [q]}\customvarpi_{z,S}^{(a_i)}  \in \ZZ[\beta][x_{a_1},x_{a_2},\dots,x_{a_q}]
\ee where
\be\label{customvarpi-eq}
\customvarpi_{z,S}^{(a)} := \begin{cases}
-1&\text{if $a \in S$ and $S$ contains any crossing bound of $a$} \\
2+\beta x_{a} &\text{if }a \notin S \\
1 + \beta x_{a} &\text{otherwise}.
\end{cases}
\ee
If $z$ has an immobile left segment $L$ (that is, containing $1$) and $S \in \shiftable(z)$,
then  $\customvarpi_{z,S} $ is divisible by the product $\prod_{a \in L} (2+\beta x_a)$
 since $S\cap L =\varnothing$.

\begin{example} \label{sigma-ex2}
If $z=(2, 7)(3, 8)(4, 6)(5, 9)$ is as in Example~\ref{sigma-ex} then 
\[\ba \customvarpi_{z,\varnothing} &= (2+\beta x_2)(2+\beta x_3)(2+\beta x_4)(2+\beta x_5),\\
 \customvarpi_{z,\{2,4,5\} } &= -(1+\beta x_2)(2+\beta x_3)(1+\beta x_5).\ea\]
\end{example}

The formula in Example~\ref{ivex-ex} generalizes to the following theorem.

\begin{theorem}\label{ivex-thm}
If $z \in \Ivex$ then $\fkGO_z = \sum_{S\in\shiftable(z)}   \beta^{|S|}\cdot \customvarpi_{z,S}  \cdot \iG_{z , S}$.
\end{theorem}

We give the proof of this theorem in the next section.
This result should be compared with \cite[Thm.~3.26]{MP2020},
which presents a different formula for $\fkGO_z$ in the vexillary case in terms of certain Pfaffians.
Our new formula in Theorem~\ref{ivex-thm} is more explicit and easier to compute.

For a general $z \in \Ivex$, 
the coefficients $\customvarpi_{z,S}$ may have negative signs although
we know that $\fkGO_z\in \NN[\beta][x_1,x_2,\dots]$.
Thus, there are typically many cancellations in the formula in Theorem~\ref{ivex-thm}.

Define
$z \in \Ivex$ to be \defn{locally noncrossing} if
any numbers $a_i<a_j$ in the same mobile left segment of $z$ have $z(a_i)>z(a_j)$.
As each $\iG_y\in \NN[\beta][x_1,x_2,\dots]$, 
the following proposition shows that
 Theorem~\ref{ivex-thm} gives a (cancellation-free) monomial-positive formula for $\fkGO_z$ precisely when  $z $ 
 is
 locally noncrossing.
  
\begin{proposition}\label{locally-nc-prop}
Let $z \in \Ivex$. Then $\customvarpi_{z,S} \in \NN[\beta][x_1,x_2,\dots]$
for all shiftable sets $S \in \shiftable(z)$ if and only if $z$ is locally noncrossing.
\end{proposition}

\begin{proof}
If $z$ is locally noncrossing then the first case in \eqref{customvarpi-eq} never occurs,
so each $\customvarpi_{z,S}^{(a)} \in \{ 2+\beta x_a, 1+\beta x_a\}$ and $\customvarpi_{z,S} \in \NN[\beta][x_1,x_2,\dots]$.

Conversely, suppose $z$ is not locally noncrossing,
so that $z$ has a mobile left segment $L$ containing numbers $a_i<a_j$ with $z(a_i)>z(a_j)$.
Assume $a_j$ is as small as possible, and then let $a_i \in L$ be as large as possible
while still satisfying $a_i < a_j$ and $z(a_i)>z(a_j)$.

The minimality of $a_j$ implies that $a_i$ is not a crossing bound for any element of $L$.
The maximality of $a_i$ implies that $a_j$ is a crossing bound for $a_i$ but not for any other element of $L$.
Hence the set $S = \{a_i,a_j\} $ belongs to $ \shiftable(z)$, but we have
$
\customvarpi_{z,S}^{(a_i)} =-1
$  
while every other $\customvarpi_{z,S}^{(a)} \in \{ 2+\beta x_a, 1+\beta x_a\}$.
For this shiftable set we have $\customvarpi_{z,S} \notin   \NN[\beta][x_1,x_2,\dots]$.
\end{proof}

 The numbers $1$, $2$, $4$, $9$, $20$, $47$, $109$, $261$, $621$, $1511$, $3664$, $9052$, $\dots$  of locally noncrossing vexillary involutions in $S_n$ 
 do not appear to given a known integer sequence.
The slightly more natural sequence of locally noncrossing vexillary involutions $z \in S_{n+1}$ with $z(1)=1$
(whose counts are $1$, $2$, $4$, $8$, $17$, $36$, $77$, $166$, $359$, $780$, $1700$, $3715$, $\dots$)
also appears to be unknown.

\subsection{Proof of Theorem~\ref{ivex-thm}}\label{proof-sect}

We continue the setup of the previous section, so $z \in \Ivex$ is a vexillary involution
with left endpoints $a_1<a_2<\dots<a_q$ and right endpoints $b_1,b_2,\dots,b_q$.
We begin with three technical lemmas about shiftable subsets.

\begin{lemma}\label{shiftable-lem}
Let $S \in \shiftable(z)$. Choose a left segment $L$ of $z$ with $i \in S\cap L$.
Write $j $ for the smallest element $j \in S$ with $i<j$, when this exists.
\ben
\item[(a)] If $\crb_z(i)\cap S = \varnothing$
and $i \neq \max (S\cap L)$ then $z(i) > z(j)$.

\item[(b)] If $\crb_z(i)\cap S \neq \varnothing$
then $i \neq \max (S\cap L)$ and $z(i) < z(j)$.

\een
\end{lemma}

\begin{proof}
First suppose $\crb_z(i)\cap S = \varnothing$
and $i \neq \max (S\cap L)$ but 
the smallest $j \in S$ with $i<j$ has  $z(i) < z(j)$. Then since $j$ cannot be a crossing bound of $i$,
the maximal $t$ with $ i \leq t < j $ and $z(t)<z(j)$ must be greater than $i$.
But then $j \in \crb_z(t)$ so we must have $t \in S$ since $j \in S$ and $S$ is shiftable.
This contradicts the minimality of $j$ and proves part (a).

For part (b), suppose $\crb_z(i)\cap S \neq \varnothing$. Then some $k \in S\cap L$
has $i<k$ and $k \in \crb_z(i)$, so $i\neq\max(S\cap L)$. 
If $j=k$ then $z(i) < z(j)$ by the definition of a crossing bound.
If $j<k$ then since $\{i\} = \{ t: i\leq t < k \text{ and }z(t)<z(k)\}$,
 we must have $z(i) < z(k) < z(j)$.
\end{proof}

\begin{lemma}\label{shiftable-lem2}
Suppose $L$ is a left segment of $z$ and $S \in \shiftable(z)$. Assume $i<j$ are left endpoints with 
$i \in S\cap L$ and $j \in L\setminus S$.
If every integer $t$ with $i\leq t <j$ has $z(t)>z(j)$, 
then $S  \sqcup\{j\} \in \shiftable(z)$.
\end{lemma}

\begin{proof}
We argue by induction on $i$. 
The union $S \sqcup\{j\}$ can only fail to be shiftable if some left endpoint $h\in L$  has $h \notin S$ and $j \in \crb_z(h)$. Such an element must have $z(h)<z(j)$ so $h<i$, and then every $t$ with $h<t<j$ must have $z(t)>z(j)$. This is impossible if $i = \min (L)$, so the lemma holds in that base case.

Assume $i> \min(L)$ and such as $h$ exists.
As we cannot have $ i \in \crb_z(h)$ since $h\notin S$ and $i \in S$, there is some maximal $t$ with $h<t <i$ and $z(j)<z(t)<z(i)$.
Then by construction $i \in \crb_z(t)$ so  $t \in S\cap L$ since $S$ is shiftable.
But this means the hypotheses of the lemma are satisfied 
if we replace $i$ by the smaller left endpoint $t$, so we conclude by induction that $S \sqcup\{j\} \in \shiftable(z)$.
\end{proof}

\begin{lemma}\label{shiftable-lem3}
Suppose $L$ is a left segment of $z$ that intersects $S \in \shiftable(z)$. 
Let $m:= \min(S\cap L)$ and assume $i\in L\setminus S$ has $i< m$. 
Then $z(m) < z(i)$.
\end{lemma}

\begin{proof}
Suppose instead that $z(i)<z(m)$. Since $S$ is shiftable but does not contain $i$, the left endpoint
$m$ cannot be a crossing bound for $i$.
Hence the maximal left endpoint $t\in L$ with $i\leq t<m$ and $z(t) <z(m)$ is greater than $i$.
But then $m $ is a crossing bound for $t$, so we must have $t \in S$
since $S$ is shiftable and contains $m$.
This contradicts the minimality of $m=\min(S\cap L)$.
\end{proof}

We now begin the main proof of this section.

\begin{proof}[Proof of Theorem~\ref{ivex-thm}]
The basic strategy in this long proof is fairly straightforward.
We will argue by reverse induction on the length of $z$, the base case being when $z$ is dominant. Then
$z$ only has one left segment, which is immobile, so only the empty set is shiftable,
and the desired formula for $\fkGO_z$ reduces to the one given in Theorem~\ref{qd-thm}.

Now assume that $z \in \Ivex$ is such that $\fkGO_z = \sum_{S\in \shiftable(z)}  \customvarpi_{z,S} \cdot \beta^{|S|}\cdot \iG_{z , S}$. 
Suppose $i \in \PP$ has $z(i)>z(i+1)$ and $z \neq s_i zs_i  \in \Ivex$. Set \[y := s_i zs_i \in \Ivex.\]
Then $\fkGO_y = \bpartial_i \fkGO_z$ by \eqref{orthogonal-recursion},
so by Theorem~\ref{vexweak-prop} it suffices to show that
\be\label{goal-eq}
 \sum_{S\in \shiftable(y)}\beta^{|S|}\cdot  \customvarpi_{y,S} \cdot  \iG_{y , S}
 =
 \sum_{S\in \shiftable(z)}  \beta^{|S|}\cdot \bpartial_i(\customvarpi_{z,S} \cdot \iG_{z , S}).
 \ee
 The argument to check this splits into several cases according to how $i$ and $i+1$ are
 distributed among the left endpoints, right endpoints, and fixed points of $z$.
 
 Since $i$ is a descent of $z$, 
 it cannot occur that $i$ and $i+1$ are both fixed points of $z$, or
 that $i$ is a fixed point and $i+1$ is a left endpoint,
 or
 that $i$ is a right endpoint and $i+1$ is a fixed point
  Also, since $y$ and $z$ are distinct and both vexillary, it cannot occur that
 $i$ is a right endpoint and $i+1$ is a left endpoint,
 or that $i$ is a left endpoint and $i+1$ is a right endpoint.
 These considerations exclude 5 of 9 cases, leaving 4 possibilities:
  \ben
 \item[(a)] The number $i$ is a fixed point of $z$ while $i+1$ is a right endpoint.
  \item[(b)] Both $i$ and $i+1$ are right endpoints of $z$.
 \item[(c)] Both $i$ and $i+1$ are left endpoints of $z$.
  \item[(d)] The number $i$ is a left endpoint of $z$ while $i+1$ is a fixed point. 
\een
 We explain why \eqref{goal-eq} holds in each of these cases below.
 
   \ben
  
 \item[(a)] Suppose $i$ is a fixed point of $z$ and $i+1$ is a right endpoint.
 \een
 
 In this case $y$ and $z$ have the same left segments, crossing bounds,
 and shiftable subsets.
 Suppose $S \in \shiftable(y)= \shiftable(z)$.
 As $i$ must be greater than the largest left endpoint of  $z\in \Ivex$,
  it follows that $\bpartial_i$ commutes with $\customvarpi_{y,S}=\customvarpi_{z,S}$ and $s_i$ commutes with $\sigma_S$. Moreover,  $i$ is still a fixed point and $i+1$ is still a right endpoint of $(\sigma_S)^{-1} \cdot z \cdot \sigma_S$, so we have $\bpartial_i \iG_{z,S} = \iG_{y,S}$ and
   \be\label{case-a-eq}
  \bpartial_i(\customvarpi_{z,S} \cdot \iG_{z , S}) 
  = \customvarpi_{y,S}  \cdot \iG_{y,S}
 .\ee
Substituting this formula turns the left side of \eqref{goal-eq} into the right.
 
 \ben
 \item[(b)] Suppose $i$ and $i+1$ are both right endpoints of $z$.
 \een
 
In this case $i$ is again greater than the largest left endpoint of  $z$,
so $\bpartial_i$ commutes with $\customvarpi_{z,S}$ for all $S \in \shiftable(z)$.
 Since $i$ is a descent of $z$, for some indices $1\leq j <k \leq q$
 it holds that $z(i+1) = a_j$ and $z(i) = a_k$. If $a_j$ and $a_k$ are not in the same left segment of $z$,
 or are both in the unique immobile left segment (when this exists),
 then the argument to deduce \eqref{goal-eq} is the same as in case (a),
 since then $y$ and $z$ have the same left segments, crossing bounds,
 and shiftable subsets, so \eqref{case-a-eq} again holds.

Instead assume that $a_j$ and $a_k$ belong to the same mobile left segment of $z$.
 Then $y=s_izs_i$ has the same left segments as $z$,
 and almost the same nontrivial cycles: 
 the only difference is that $(a_j,i+1)$ and $(a_k,i)$ have been replaced by $(a_j,i)$ and $(a_k,i+1)$.
 This means that the crossing bounds 
 for the left endpoints of $y$ are the same as for $z$, with one exception:
 if $i+1<b_t$ for all $j<t<k$ 
then $a_k$ will be a crossing bound for $a_j$ in $y$ but not in $z$. 
 If this exception does not occur then \eqref{case-a-eq} again holds for all
 $S \in \shiftable(y)=\shiftable(z)$, which suffices to prove \eqref{goal-eq}.
For the exceptional case, we make a claim:

 \begin{claim} Assume we are in the case when $i+1<b_t$ for all $j<t<k$. Let 
 \[
 \ba
 \cE &:=\{ S \in \shiftable(z) : a_j \notin S \text{ and }a_k \in S\}, \\
 \cF &:= \{ S \in \shiftable(z) : a_j,a_k \in S \text{ and }a_t \notin S\text{ for all }j<t<k\} .
 \ea
 \]
Then 
$ \shiftable(y) =  \shiftable(z)  \setminus \cE$
and $S \mapsto S\sqcup\{a_j\}$ is a bijection $\cE \to \cF$.
\end{claim}

\begin{claimproof}
We can always remove $a_j$ from $T \in \cF$ to obtain a shiftable set  $T\setminus\{a_j\}\in\cE$,
since such a set $T$ contains no crossing bounds for $a_j$, as these would have to be of the form $a_t$ for some $j<t<k$.
On the other hand, 
if $S \in \cE$ then  any left endpoint $a_t$ that has $a_j$ as a crossing bound in $z$ must also have $a_k$ 
 as a crossing bound, and therefore is already in $S$, so we can add $a_j$ to get a shiftable subset 
 $S\sqcup\{a_j\} \in \shiftable(z)$. 
 
Now we check that $a_t \notin S$ for all $j<t<k$ if $a_j \notin S \in \shiftable(z)$.
 This follows by induction on $t$ as each $a_t$ is a crossing bound for one of $a_j,a_{j+1},\dots,a_{t-1}$ as
  $b_j =i+1< b_t$.
Hence 
 if $a_j,a_{j+1},\dots,a_{t-1} \notin S$ then also  $a_t \notin S$.
 We conclude that if $S \in \cE$ then the union $S\sqcup\{a_j\} $ is shiftable and contained in $\cF$.
 \end{claimproof}
 
 Continue to assume $i+1<b_t$ for all $j<t<k$, so that the claim applies.
If $S \in \shiftable(y) \setminus \cF$ then $S$ contains neither $a_j$ nor $a_k$ so $\customvarpi_{y,S} = \customvarpi_{z,S}$.
If instead $S \in \cF$ then 
 we obtain $\customvarpi_{y,S}$ from $\customvarpi_{z,S} =\prod_{t\in[q]} \customvarpi_{z,S}^{(a_t)}$ by replacing the single 
 factor $\customvarpi_{z,S}^{(a_j)} = 1+\beta x_{a_j}$ by $\customvarpi_{y,S}^{(a_j)} =-1$. This means that if $S \in \cF$ 
 then
\be
\label{case-b-eq1}
\customvarpi_{y,S} = \customvarpi_{z,S} - \customvarpi_{z,T}\quad\text{for $T := S\setminus \{a_j\} \in \cE$.}
\ee
 Finally, we observe that if $S \in \shiftable(y)=\shiftable(z)\setminus\cE$ then 
 \be
\label{case-b-eq2}
  \beta^{|S|} \cdot \bpartial_i(\customvarpi_{z,S} \cdot \iG_{z , S}) = \customvarpi_{z,S}  \cdot  \beta^{|S|} \cdot \bpartial_i \iG_{z,S} =\customvarpi_{z,S}  \cdot \beta^{|S|} \cdot \iG_{y,S}
  \ee
  while if $T \in \cE$ and $S := T \sqcup\{a_j\} \in \cF$ then, as $a_t \notin T$ for $j<t<k$, we have
  \be
\label{case-b-eq3}
   \beta^{|T|}\cdot  \bpartial_i(\customvarpi_{z,T} \cdot \iG_{z ,T}) = 
      \customvarpi_{z,T}  \cdot  \beta^{|T|}\cdot \bpartial_i \iG_{z,T} =
       -    \customvarpi_{z,T}  \cdot   \beta^{|S|} \cdot \iG_{y,S}.
 \ee
 Since $T\mapsto T\sqcup\{a_j\}$ is a bijection $\cE\to\cF$,
 substituting \eqref{case-b-eq2} and  \eqref{case-b-eq3} into the right side of \eqref{goal-eq}
 and then using \eqref{case-b-eq1} gives the left side  \eqref{goal-eq}, as needed.

   \ben
 \item[(c)] Suppose $i$ and $i+1$ are both left endpoints of $z$.
 \een
 
 In this case $y$ and $z$ have the same left endpoints and left segments,
and the numbers $i$ and $i+1$ belong to the same left segment.
Notice that the set $\crb_z(i) $ of crossing bounds  for $i$ in $z$ is empty since whenever $i <t $ and $z(i) < z(t)$ we also have $z(i+1)<z(t)$.

 Our proof of \eqref{goal-eq} in this case reduces to two somewhat technical claims,
  which we state and prove below.
 
 \begin{claim}
 \label{c-claim1}
  Continue to assume $i$ and $i+1$ are both left endpoints of $z$. Then
   \[ 
 \sum_{\substack{S \in \shiftable(z) \\ i \in S}}  \beta^{|S|}\cdot  \bpartial_i(\customvarpi_{z,S} \cdot \iG_{z , S})
  =
  0.
  \]
  \end{claim}
  
  \begin{claimproof}
  If $i$ is part of a left segment of $z$ containing $1$, then the equality holds trivially because the sum is empty.
  We may therefore 
  assume that $i$ and $i+1$ belong to a left segment that does not contain $1$. (This assumption is not actually used in any part of the following argument, however.)
  To prove the claim, it suffices to verify the following three statements:
   \ben
 
 \item[(1)] If  $ S \in \shiftable(z)$ has $i,i+1 \in S$ and $S \cap \crb_z(i+1)=\varnothing$ then
  \[
\beta^{|S|}\cdot  \bpartial_i(\customvarpi_{z,S} \cdot \iG_{z , S}) =  \customvarpi_{z,S} \cdot \beta^{|S|}\cdot\iG_{z,T}
 \]
for the set  $T:=S\setminus\{i+1\} $, which belongs to $ \shiftable(z)$ and satisfies
  \[
 \beta^{|T|}\cdot  \bpartial_i(\customvarpi_{z,T} \cdot \iG_{z , T}) = - \customvarpi_{z,S} \cdot \beta^{|S|}\cdot    \iG_{z , T}.
 \]
 
 \item[(2)] If $T \in \shiftable(z)$ has $i \in T$ and $i+1\notin T$ then $ S:= T\sqcup\{i+1\} $ satisfies
 \[
 i,i+1\in S \in \shiftable(z)
 \quand
 S  \cap \crb_z(i+1)=\varnothing.\]

  \item[(3)] If $i,i+1 \in S\in \shiftable(z)$ and $S \cap \crb_z(i+1)\neq \varnothing$ then
$
\bpartial_i(\customvarpi_{z,S}  \iG_{z , S}) = 0.
$
\een

To show (1), choose  $ S \in \shiftable(z)$ with $i,i+1 \in S$ and $S \cap \crb_z(i+1)=\varnothing$  and define  $T:=S\setminus\{i+1\} $. Since $S$ does not contain any crossing bounds for $i+1$,
it is evident that $T \in \shiftable(z)$.
Since $\crb_z(i)=\varnothing$, we have 
\[\customvarpi_{z,S}^{(i)}\customvarpi_{z,S}^{(i+1)} = (1+\beta x_i)(1+\beta x_{i+1})\] so $\bpartial_i$ commutes with $\customvarpi_{z,S}$.
It follows from Lemma~\ref{shiftable-lem}(a) that  
\[\bpartial_i( \iG_{z , S}) = \iG_{z,T}
\quand 
\bpartial_i( \iG_{z , T}) =-\beta \iG_{z,T}.
\]
Therefore, we conclude from \eqref{fg-eq} that
\[
\beta^{|S|} \cdot  \bpartial_i(\customvarpi_{z,S} \cdot \iG_{z , S}) 
=
 \customvarpi_{z,S} \cdot  \beta^{|S|} \cdot\iG_{z,T}
\]
and
\[
\beta^{|T|} \cdot  \bpartial_i(\customvarpi_{z,T} \cdot \iG_{z , T}) =  (\bpartial_i \customvarpi_{z,T})\cdot \beta^{|S|-1} \cdot  \iG_{z , T}.
\]
It remains to check that $\bpartial_i \customvarpi_{z,T} = -  \customvarpi_{z,S} \cdot \beta$. As we have \[
\customvarpi_{z,T}^{(i)}\customvarpi_{z,T}^{(i+1)} = (1+\beta x_i)(2+\beta x_{i+1})\]
the desired identity follows by checking that 
\[
\bpartial_i\( (1+\beta x_i)(2+\beta x_{i+1})\) =  -\beta(1+\beta x_i)(1+\beta x_{i+1}).
\]
This proves item (1).

We proceed to item (2). Assume $T \in \shiftable(z)$ has $i \in T$ and $i+1\notin T$. Define $ S:= T\sqcup\{i+1\} $. Since $T$ is shiftable, we must have $T\cap \crb_z(i+1) = \varnothing$, so it also holds that $S \cap \crb_z(i+1) = \varnothing$. Finally, Lemma~\ref{shiftable-lem2} shows that $S \in \shiftable(z)$
as needed.

For item (3), assume $S\in \shiftable(z)$ has $i,i+1 \in S$ and $S \cap \crb_z(i+1)\neq \varnothing$. Then 
$
\customvarpi_{z,S}^{(i)}\customvarpi_{z,S}^{(i+1)}  = -1-\beta x_i
$ so  $\bpartial_i \customvarpi_{z,S}=0$.
It suffices in view of \eqref{fg-eq} to check that $\bpartial_i \iG_{z , S} = -\beta \iG_{z,S}$, 
and this follows by Lemma~\ref{shiftable-lem}(b).
  \end{claimproof}
 
 Our second claim involves the following definitions.
 Given $S\in\shiftable(z)$, let $m_S$ be the smallest element of $ S$ with $i+1<m_S$, if this exists. Then let
 \[
 \ba
 \cA &:= \{ S \in \shiftable(z) : i,i+1 \notin S\}, \\
 \mathcal{B} &:= \{ S \in \shiftable(z) : i\notin S,\ i+1\in S,\ S \cap \crb_z(i+1)= \varnothing \}, \\
\cC &:= \{ S \in \shiftable(z) : i\notin S,\ i+1\in S,\ S \cap \crb_z(i+1)\neq \varnothing,\ z(m_S)<z(i) \}, \\
\cD &:= \{ S \in \shiftable(z) : i\notin S,\ i+1\in S,\ S \cap \crb_z(i+1)\neq \varnothing,\ z(i)<z(m_S) \}.
 \ea
 \]
 Lemma~\ref{shiftable-lem} ensures that $m_S$ is defined in the formulas for $\cC$ and $\cD$.
If $X$ and $Y$ are two sets then we let $X\vartriangle Y := (X \setminus Y) \sqcup (Y\setminus X)$.
To establish \eqref{goal-eq}, it suffices to combine our first claim above with the following claim below:
 
 \begin{claim}\label{c-claim2}
Continue to assume $i$ and $i+1$ are both left endpoints of $z$.   Then 
\[ \shiftable(y) = \cA \sqcup \{ S \vartriangle\{i,i+1\} : S \in \mathcal{B} \sqcup \cC\}
 \sqcup \{ S \sqcup \{i\} : S \in \mathcal{B} \sqcup \cC \sqcup \cD\}.\]
 Moreover, if $S \in \cA\sqcup\mathcal{B}\sqcup\cC\sqcup\cD$ and
$T := S\vartriangle \{i,i+1\}$ and $U := S\sqcup\{i\}$,  then
\[
\beta^{|S|}\cdot  \bpartial_i(\customvarpi_{z,S} \cdot \iG_{z , S})=
\begin{cases}
   \beta^{|S|}\cdot\customvarpi_{y,S} \cdot \iG_{y,S}&\text{if }S \in \cA \\
  \beta^{|T|}\cdot  \customvarpi_{y,T}\cdot  \iG_{y , T} +\beta^{|U|}
  \cdot\customvarpi_{y,U} \cdot  \iG_{y , U} &\text{if }S \in \mathcal{B}\sqcup\cC \\
  \beta^{|U|}  \cdot \customvarpi_{y,U} \cdot \iG_{y , U} &\text{if }S \in \cD.
  \end{cases}
\]
 \end{claim}
 
 \begin{claimproof}
 We can derive the formula for $\beta^{|S|}\cdot  \bpartial_i(\customvarpi_{z,S} \cdot \iG_{z , S})$
directly from \eqref{fg-eq},
 since the setup determines the value of $\customvarpi_{z,S}^{(i)} \customvarpi_{z,S}^{(i+1)}$ in each case
 and Lemma~\ref{shiftable-lem} tells us how to evaluate $\bpartial_i \iG_{z,S}$ via Theorem~\ref{iG-thm}.
 The argument is similar to the proof of Claim~\ref{c-claim1} (and also to the proofs of Claims~\ref{d-claim2} and \ref{d-claim3} below) so we omit the details.

 The nontrivial statement left to prove is the desired expression for $\shiftable(y)$.
 To show this, we start by observing that while
 $\crb_z(i) = \varnothing$, it holds that 
\[\ba 
\crb_y(i) &= \{i+1\}\sqcup \{ t \in \crb_z(i+1) : z(i+1) < z(t)< z(i)\},
\\
\crb_y(i+1) &= \{ t \in \crb_z(i+1) : z(i) < z(t)\}.\ea
\]
Meanwhile, for any other left endpoint $t$, the set of crossing bounds $\crb_y(t)$ is either equal to 
$\crb_z(t)$, or formed from $\crb_z(t)$ by removing $i$ or $i+1$, and the second case only occurs when $\crb_z(t)$ already contains at least $i$. Specifically, for any left endpoint $t \notin\{i,i+1\}$, we are in one of the following cases:
\ben

\item[(i)] when $t<i$ and $z(t) < z(i+1)$, it can happen that 
\[ i,i+1\in \crb_z(t)
\quand   \crb_y(t) = \crb_z(t) \setminus\{i+1\}.\]

  \item[(ii)] when $t<i$ and $z(i+1)<z(t)<z(i)$,  it can happen that
  \[i \in \crb_z(t),\ i+1\notin\crb_z(t)
  \quand
\crb_y(t) = \crb_z(t) \setminus\{i\}.\]

 \item[(iii)] if cases (i) or (ii) do not apply, then $i,i+1 \notin  \crb_y(t)=\crb_z(t)$.
 \een
 It follows that if $S$ is a set containing neither $i$ nor $i+1$, then $S \in \shiftable(z)$ if and only if $S \in \shiftable(y)$. 
 
Alternatively, if $S \in \shiftable(z)$ and $i \notin S$,
 then the only way $S\sqcup\{i\}$ can fail to be in $\shiftable(y)$ is if some left endpoint $t \notin S$ 
 has $i \in \crb_y(t)$. In view of (i), (ii), and (iii), this can only happen when $i+1 \in \crb_z(t)$, which would mean that $i+1 \notin S$.
 Thus, if $S \in \shiftable(z)$ has $i \notin S$ and $i+1 \in S$, so that 
 $S \in \mathcal{B}\sqcup\cC\sqcup\cD$,
 then we have $S\sqcup\{i\}\in \shiftable(z)$ and this union contains both $i$ and $i+1$.
 
 Conversely, if $U \in \shiftable(y)$ has $i,i+1 \in U$ then since $\crb_z(i)=\varnothing$,
 the only way $U \setminus \{i\}$ can fail to be in $\shiftable(z)$ is if some left endpoint $t \notin U$ 
 has $i +1\notin \crb_y(t)$ but $i+1 \in \crb_z(t)$. Inspecting (i) and (ii) shows that this situation is impossible, 
so if $U \in \shiftable(y)$ has $i,i+1 \in U$  then $U \setminus\{i\} \in \shiftable(z)$.

Finally suppose $S$ is a set containing $i+1$ but not $i$. Form $T := S\vartriangle\{i,i+1\}$ by removing $i+1$ and adding $i$. 
It is not possible for a left endpoint $t$ to have $i \in \crb_y(t)$ and $i+1\notin \crb_z(t)$,
so if $S \in \shiftable(z)$ then we can only fail to have $T \in \shiftable(y)$
if $S \cap \crb_y(i+1) \neq \varnothing$.
This only occurs if $S$ contains some $t \in \crb_z(i+1)$ with $z(i)<z(t)$,
 but then we must have $t=m_S$ or $z(i)<z(t) < z(m_S)$.
We conclude that if $S \in \mathcal{B}\sqcup\cC$ then $T \in \shiftable(y)$.

Conversely, since $\crb_z(i)=\varnothing$, if $T \in \shiftable(y)$
then we can only fail to have $S \in \shiftable(z)$
is there is some left endpoint $t \notin T$ with  $i\notin \crb_y(t)$ and $i+1 \in \crb_z(t)$.
Inspecting (i) and (ii) shows that this situation is impossible, 
so if $T \in \shiftable(y)$ then $S \in \shiftable(z)$.
In fact, 
if $T \in \shiftable(y)$ then $S \in \mathcal{B} \sqcup \cC$,
since $i+1 \in S$ and if $m_S$ is defined
then $y(m_S)< y(i+1)$ by Lemma~\ref{shiftable-lem3}
so
$z(m_S) =y(m_S)< y(i+1)=z(i)$.

This argument shows that the given union is contained in $\shiftable(y)$,
and that each element of $\shiftable(y)$---except those containing $i+1$ but not $i$---is contained in this union. 
As no set in $\shiftable(y)$ contains $i+1$ but not $i$ since $i+1 \in \crb_y(i)$,
 this completes the proof of our claim.
 \end{claimproof}
 
The identity \eqref{goal-eq} in case (c) is immediate from Claims~\ref{c-claim1} and \ref{c-claim2}.
This leaves just one final case to consider:

 \ben
 \item[(d)] Suppose $i$ is a left endpoint of $z$ and $i+1$ is a fixed point. 
 \een

In this case $i$ is the last element of its left segment in $z$, so $\crb_z(i) = \varnothing$.
However, the left endpoints for $y$ are no longer the same as for $z$.  Compared to $z$,
the set of left endpoints for $y$ loses $i$ and gains $i+1$.

This means that the left segments of $y$ also differ from $z$.
Let $L$ be the left segment of $z$ containing $i$.
If $i+2$ is not a left endpoint of $z$, then the left segments of $y$ are formed from those of $z$
by replacing $L$ by $L\setminus\{i\}$ and adding a new singleton segment $\{i\}$.
If $i+2$ is a left endpoint of $z$, then it begins a new left segment $M$ of $z$.
In this event,  the left segments of $y$ are formed from those of $z$
by replacing $L$ by $L\setminus \{i\}$ and $M$ by $\{i+1\}\sqcup M$.

These sub-cases require slightly different arguments,
which we handle as usual through a series of claims. To start, we have this close relative of Claim~\ref{c-claim1}:

 \begin{claim}\label{d-claim1} Continue to assume that $i<z(i)$ and $i+1=z(i+1)$. Then
   \[ 
 \sum_{\substack{S \in \shiftable(z) \\ i \in S}}  \beta^{|S|}\cdot  \bpartial_i(\customvarpi_{z,S} \cdot \iG_{z , S})
  =
  0.
  \]
  \end{claim}
  
  \begin{claimproof}
This holds for a much simpler reason than Claim~\ref{c-claim1}.
Assume $S\in \shiftable(z)$ has $i \in S$. Since $i$ is the last element of its left segment, we have $\crb_z(i)=\varnothing$, so 
$
\customvarpi_{z,S}^{(i)}  = 1+\beta x_i$
 and  $\bpartial_i \customvarpi_{z,S}=0$.
At the same time, it is also clear that  $\bpartial_i \iG_{z , S} = -\beta \iG_{z,S}$
so in fact $\bpartial_i(\customvarpi_{z,S} \cdot \iG_{z , S})=0$ by \eqref{fg-eq}.
  \end{claimproof}
  
  We now have two variants of Claim~\ref{c-claim2}.
  Let $\cX = \{S \in \shiftable(z) : i \notin S\}$. 
  
\begin{claim}\label{d-claim2}
If $i$ but not $i+2$ is a left endpoint of $z$, and $i+1=z(i+1)$,
  then \[\shiftable(y) = \cX \sqcup \{ S\sqcup\{i+1\} : S \in \cX\}.\]
  Moreover, in this case, if $S \in \cX$ and $T := S\sqcup\{i+1\}$ then 
  \[
    \beta^{|S|}\cdot  \bpartial_i(\customvarpi_{z,S} \cdot \iG_{z , S})
    =
  \beta^{|S|}\cdot    \customvarpi_{y,S}\cdot    \iG_{y , S} + \beta^{|T|}
  \cdot\customvarpi_{y,T} \cdot   \iG_{y , T}.
  \]
  \end{claim}
  
  \begin{claimproof}
  In this case $\crb_y(t) = \crb_z(t) \setminus \{i\}$ for all left endpoints $t\neq i$,
  although the set difference may often satisfy $\crb_z(t) \setminus \{i\} = \crb_z(t)$.
  This property makes it clear that $\cX \subset \shiftable(y)$,
   Then, since $\{i\}$ is a left segment in $y$,
  it also follows that $ \{ S\sqcup\{i+1\} : S \in \cX\} \subset \shiftable(y)$.

  Conversely, if $S \in \shiftable(y)$ then $i\notin S$ as $i$ is not a left endpoint of $y$,
  and the same observation about crossing bounds implies that $S\setminus\{i+1\} \in \cX$.
  This proves our formula for $\shiftable(y)$. 
  
  Now let $S \in \cX$ and $T := S\sqcup\{i+1\}$.
  Let $\cQ = \prod_{a \neq i,i+1} \customvarpi_{z,S}^{(a)} = \prod_{a \neq i, i+1} \customvarpi_{y,S}^{(a)}$
  where the products are over the left endpoints $a$ shared by $y$ and $z$. Then
  \[
  \customvarpi_{z,S} = (2+\beta x_i)\cQ,
  \quad 
  \customvarpi_{y,S} = (2+\beta x_{i+1})\cQ,
  \quand
  \customvarpi_{y,T} = (1+\beta x_{i+1})\cQ.
  \]
As $\cQ$ does not involve the variables $x_i$ or $x_{i+1}$,
we have
\[ s_i(\customvarpi_{z,S}) = \customvarpi_{y,S}
\quand
  \bpartial_i(\customvarpi_{z,S}) = -\beta \cQ,\] so we can compute that
\[ s_i(\customvarpi_{z,S})\cdot \beta +  \bpartial_i(\customvarpi_{z,S}) = \beta \cdot ((2+\beta x_{i+1}) - 1) \cdot \cQ = 
\beta\cdot \customvarpi_{y,T}.\]
   Since  it is clear that
  $\iG_{z,S} = \iG_{y,T}$ and $\bpartial_i \iG_{z,S} = \iG_{y,S}$,
  the desired expression for $  \beta^{|S|}\cdot  \bpartial_i(\customvarpi_{z,S} \cdot \iG_{z , S})$ follows from \eqref{fg-eq} using the preceding identities.
  \end{claimproof}
  
Now suppose that $i$ and $i+2 $ are both left endpoints of $z$ while $i+1=z(i+1)$.
As above let $M$ be the left segment of $z$ containing $i+2$.
Then define
\[
\ba
\cV &= \{ S \in \cX:  S\cap M =\varnothing\text{ or } z(i)>z(\min(S\cap M)) \},\\
\cW &=\{ S \in \cX : S\cap M\neq\varnothing\text{ and }z(i)<z(\min(S\cap M)) \}. 
\ea
\]
Here is our second variant of Claim~\ref{c-claim2}:

 \begin{claim}\label{d-claim3}
Assume $i$ and $i+2$ are left endpoints of $z$ while $i+1=z(i+1)$.   Then 
\[ \shiftable(y) =  \cV
 \sqcup \{ S \sqcup \{i+1\} : S \in  \cV \sqcup \cW=\cX\}.\]
 Moreover, if $S \in\cV\sqcup\cW$ and $T :=   S\sqcup\{i+1\}$  then
\[
  \beta^{|S|}\cdot  \bpartial_i(\customvarpi_{z,S} \cdot \iG_{z , S})=
\begin{cases}
  \beta^{|S|}\cdot \customvarpi_{y,S}\cdot   \iG_{y , S} +\beta^{|T|}  \cdot \customvarpi_{y,T} \cdot \iG_{y , T} &\text{if }S \in \cV \\
  \beta^{|T|}\cdot \customvarpi_{y,T} \cdot  \iG_{y , T} &\text{if }S \in \cW.
  \end{cases}
  \]
 \end{claim}
 
 \begin{claimproof}
 The proof of this statement is similar to Claim~\ref{d-claim2}, and only slightly more complicated.
 As in the earlier proof, we have $\crb_y(t) = \crb_z(t) \setminus \{i\}$ for all left endpoints $t\neq i$,
 though it may often hold that $ \crb_z(t) \setminus \{i\} =\crb_z(t)$.
 Whereas $\crb_z(i)=\varnothing$, the set of crossing bound $\crb_y(i+1)$ may be a nonempty subset of the left segment $M$.
 
 Suppose $S \in \cX$.
 Then it is clear from the previous paragraph and the definition of a shiftable set that $S\sqcup\{i+1\} \in \shiftable(y)$.
 Moreover, the only way we can fail to have $S \in \shiftable(y)$ is if $S\cap \crb_y(i+1)\neq \varnothing$. 
 
 If $S\cap M$ is empty, then this is impossible so 
$S \in \shiftable(y)$.
Assume instead that $S\cap M \neq \varnothing$.
Let $m =\min(S\cap M)$ and suppose
 $z(m) <z(i)$. Then $y(m) = z(m) < z(i) = y(i+1)$, which means that $m\notin \crb_y(i+1)$,
 and 
also that no $t \in \crb_y(i+1)$ is greater than $m$, since if $m<t$ and $y(i+1) < y(t)$ then $y(m) <y(i+1)< y(t)$.  
In other words, it follows that $S\cap \crb_y(i+1)=\varnothing$,
so $S \in \shiftable(y)$.
We conclude that if $S \in \cV$ then $S\in \shiftable(y)$.

This shows that $\shiftable(y)$ contains the given union.
Conversely, suppose $T \in \shiftable(y)$ and let $S := T\setminus \{i+1\}$.
Then $i  \notin S$ since $i$ is not a left endpoint of $y$, and it follows 
from our earlier observations about the crossing bounds of $y$ versus $z$ that
 $S \in \cX$.

Assume $i+1 \notin T$ so that $S=T$. It remains only to check that $S \in \cV$.
For this, observe that if $S\cap M \neq \varnothing$
and $m=\min(S\cap M) = \min(T\cap M)$,
then we have $y(m) < y(i+1)$ by Lemma~\ref{shiftable-lem3}
and therefore $z(m) = y(m) < y(i+1) = z(i)$.
This completes our proof of the formula for $\shiftable(y)$.

We now compute $\beta^{|S|}\cdot  \bpartial_i(\customvarpi_{z,S} \cdot \iG_{z , S})$
using \eqref{fg-eq}.
The argument is almost the same as in the proof of Claim~\ref{d-claim2}.
Choose $S \in \cX$ and set $T := S\sqcup\{i+1\}$.
  Again let $\cQ = \prod_{a \neq i,i+1} \customvarpi_{z,S}^{(a)} = \prod_{a \neq i,i+1} \customvarpi_{y,S}^{(a)}$
  where the products are over the left endpoints $a$ shared by $y$ and $z$. 
Then we have
\[  \begin{cases}
  \customvarpi_{z,S} = (2+\beta x_i)\cQ
\\
  \customvarpi_{y,S} = (2+\beta x_{i+1})\cQ
  \end{cases}
\quand
  \customvarpi_{y,T} =\begin{cases}  (1+\beta x_{i+1})\cQ &\text{if }S \in \cV \\ 
  -\cQ &\text{if } S\in \cW,
  \end{cases}
  \]
  using Lemma~\ref{shiftable-lem} for the second identity.
  As $\cQ$ does not involve $x_i$ or $x_{i+1}$,
\[ s_i(\customvarpi_{z,S}) = \customvarpi_{y,S}
\quand
  \bpartial_i(\customvarpi_{z,S}) = -\beta  \cQ .\] 
Observe that $\iG_{z,S} = \iG_{y,T}$ for all $S \in \cX$ along with $\bpartial_i \iG_{z,S} = \iG_{y,S}$ for $S \in \cV$
  and 
  $\bpartial_i \iG_{z,S} = -\beta\iG_{z,S}$ for $S \in \cW$.
  Thus, if $S \in \cV$ then \eqref{fg-eq} can be written as 
  \[
  \ba
 \beta^{|S|}\cdot  \bpartial_i(\customvarpi_{z,S} \cdot \iG_{z , S}) 
 &=
  \beta^{|S|}\cdot \cP_{y,S} \cdot \( \iG_{y , S}+\beta \iG_{y , T}\) - \beta^{|S|}\cdot \beta \cQ \cdot \iG_{y , T}
   \\&=
  \beta^{|S|}\cdot\cP_{y,S} \cdot  \iG_{y , S} +  \beta^{|T|}\cdot((2+\beta x_{i+1}) - 1)\cdot \cQ  \cdot  \iG_{y , T}
  \\&
  =  \beta^{|S|}\cdot  \cP_{y,S} \cdot   \iG_{y , S}+ \beta^{|T|}\cdot\cP_{y,T}\cdot  \iG_{y,T}
 \ea
 \]
 whereas if $S \in \cW$ then 
 \eqref{fg-eq} gives
 \[
  \beta^{|S|}\cdot  \bpartial_i(\customvarpi_{z,S} \cdot \iG_{z , S}) 
  =
  \beta^{|S|}\cdot  (-\beta  \cQ) \cdot \iG_{z , S}
  =
 \beta^{|T|} \cdot\customvarpi_{y,T}\cdot    \iG_{y,T}
  \]
  as needed.
 \end{claimproof}

The identity \eqref{goal-eq} in case (d) is immediate from Claims~\ref{d-claim1}, \ref{d-claim2}, and \ref{d-claim3}:
the first claim simplifies the sum on the right, and then whichever of the other two claims applies 
(according to whether $i+2$ is a left endpoint)  transforms the remaining terms to the sum on the left.
  This completes our proof by induction.
  \end{proof}

\subsection{Shift invariance and stable limits}\label{shift-sect}

In this section we consider how $\fkG_w$, $\iG_y$, and $\fkG_z$ 
vary when we apply natural shifting operations to the indexing permutations.
This will allow us to state some additional results about the \defn{stable limits} of these polynomials.

Define $\downshift : \ZZ[\beta][x_1,x_2,\dots]\to \ZZ[\beta][x_1,x_2,\dots]$ to be the $\ZZ[\beta]$-algebra homomorphism sending
$x_1 \mapsto 0$ and $x_{i+1} \mapsto x_{i}$ for $i\in\PP$.


Recall that we view permutations $w \in S_\infty$ as maps $w : \ZZ \to \ZZ$ with $w(i)=i$ for $i\leq 0$.
Given $w \in S_\infty$ and an integer $n\geq0$, define $w\downarrow n$ to be the permutation of $\ZZ$
that sends $i \mapsto w(i+n)-n$.
Notice for $w \in S_\infty$ that $w\downarrow n\in S_\infty$ if and only if $w$ fixes each $i \in [n]$,
and in this case $\ell(w\downarrow n) = \ell(w)$.

\begin{proposition}\label{downshift-prop}
If $w \in S_\infty$  then $\downshift(\fkG_w) = 
\begin{cases}
\fkG_{w\downarrow 1} &\text{if $w(1)=1$} \\
0&\text{if }w(1) \neq 1.
\end{cases}$
\end{proposition}

\begin{proof}
This identity is possibly well-known, and in any case follows from Remark~\ref{bcs-rem},
upon noting two observations.
First, if $w(1) = 1$ then all letters of all Hecke words for $w$ are at least $2$, and the operation on words 
that decrements every letter by one is a bijection $\cH(w) \to \cH(w\downarrow 1)$.
Second,  if $w(1) \neq 1$ then every Hecke word for $w$ must have some letter equal to $1$, so in any bounded compatible sequence $(a,i)$ for $w$ it must hold that $i_1 = a_1 = 1$.
%
%
\end{proof}

If $z \in \Ivex$ fixes each $i \in [n]$ then the permutation $z\downarrow n$ is also in $\Ivex$.

\begin{proposition}\label{yz-down-prop}
Suppose $y \in I_\infty$ and $z \in \Ivex$. Then
\[
\downshift(\iG_y) = 
\begin{cases}
\iG_{y\downarrow 1} &\text{if $y(1)=1$} \\
0&\text{if }y(1) \neq 1
\end{cases}
\quand
\downshift(\fkGO_z) = 
\begin{cases}
\fkGO_{z\downarrow 1} &\text{if $z(1)=1$} \\
0&\text{if }z(1) \neq 1.
\end{cases}
\]
\end{proposition}

\begin{proof}
One can deduce from \eqref{sim-eq} that if $y(1) =1$ then
every $w \in \cB(y)$ has $w(1)=1$ and $w\mapsto w\downarrow 1$ is a length-preserving bijection $\cB(y) \to \cB(y\downarrow 1)$.
Likewise, if $y(1)\neq 1$ then 
every $w \in \cB(y)$ also has $w(1)\neq 1$. The 
identity for $\downshift(\iG_y) $ is immediate from Proposition~\ref{downshift-prop} given these observations.

For $S \in \shiftable(z)$ let $S\downarrow 1 = \{ s - 1: s \in S\}$.
As $\downshift$ is an algebra homomorphism,
Theorem~\ref{ivex-thm} implies that \[\textstyle
\downshift(\fkGO_z)=
\sum_{S\in\shiftable(z)}  \beta^{|S|} \cdot \downshift(\customvarpi_{z,S}) \cdot  \downshift(\iG_{z , S})
.\]

 Suppose $z(1) = 1$. Let $L$ be the left segment of $z$ containing $2$ or set $L =\varnothing$ if $z(2)=2$.
Then our first identity implies that
 $\downshift(\iG_{z , S})$ is zero when $S\cap L\neq \varnothing$
 and equal to $\iG_{z\downarrow 1, S\downarrow 1}$ when $S\cap L= \varnothing$, in which case we also have $\downshift(\customvarpi_{z,S})=\customvarpi_{z\downarrow 1,S\downarrow 1}$.
 As  $S\mapsto S\downarrow 1$ is a bijection \[\{S \in \shiftable(z) : S\cap L= \varnothing\} \to \shiftable(z\downarrow 1),\] we conclude that 
 $\downshift(\fkGO_{z}) = \fkGO_{z\downarrow 1}$.
 
 Alternatively, if $z(1) \neq 1$ then $y := (\sigma_S)^{-1}\cdot z \cdot \sigma_S$ has $y(1) \neq 1$ for all $S \in \shiftable(z)$, so 
 $\downshift(\iG_{z , S})$ is always zero and $\downshift(\fkGO_z)=0$.
\end{proof}

Given $w \in S_\infty$ and any positive integer $n$, define 
$ 1^n \times w \in S_\infty$
to be the permutation that fixes each integer $i\in[n]$ while sending $i+n \mapsto w(i)+n$ for all $i\in\PP$. 
We usually write $1 \times w$ in place of $1^1\times w$.
Note that 
\[
(1^n \times w) \downarrow n = w\quad\text{for all }w\in S_\infty.\]

Now define
 $\upshift : \ZZ[\beta][x_1,x_2,\dots]\to \ZZ[\beta][x_1,x_2,\dots]$ to be the $\ZZ[\beta]$-linear map
 sending $\fkG_w \mapsto \fkG_{1\times w}$ for $w \in S_\infty$.
 By Proposition~\ref{downshift-prop}, we have 
\be\label{shift-eq} \downshift \circ \upshift = 1
 \quand
 \upshift \circ \downshift(\fkG_w) = \begin{cases}
 \fkG_w &\text{if }w(1) = 1\\
 0&\text{if }w(1)\neq 1.
 \end{cases}
\ee
Note, however, that $\upshift$ is not an algebra homomorphism, and thus differs from the more familiar right inverse of $\downshift$
sending $x_i \mapsto x_{i+1}$ for all $i \in \PP$.

 \begin{proposition}
It holds that $\upshift(\iG_z) = \iG_{1\times z}$ for all $z \in I_\infty$.
 \end{proposition}
 
 \begin{proof}
The operation $w\mapsto 1\times w$ preserves $\ell$ and $\ellhat$
and given a bijection $\cB(z) \to \cB(1\times z)$, 
so the result is clear from Definition~\ref{iG-def}.
 \end{proof}
 
The following theorem gives a similar formula for $\upshift(\fkGO_z)$.
This is the one of the main results of this section, 
and is much harder to prove.

\begin{theorem}\label{supp-thm}
Let $z \in \Ivex$.
Then $\upshift(\fkGO_z) = \fkGO_{1\times z}$ if and only if $z(1)=1$.
\end{theorem}

Before proving this theorem we consider an example and state one  lemma.

\begin{example}
We have
$\fkGO_{(1, 2)} = 2 \fkG_{21} + \beta \fkG_{312}$
but
\[
\fkGO_{1\times (1,  2)}= \fkGO_{(2,  3)} = 2 \fkG_{132} + \beta \fkG_{1423} + \beta \fkG_{231} + \beta^2 \fkG_{2413}
\neq \upshift(\fkGO_{(1,2)}).
 \] 
 However, it does hold that \[ 
\fkGO_{1\times (2,  3)}= \fkGO_{(3,  4)} = 2 \fkG_{1243} + \beta \fkG_{12534} + \beta \fkG_{1342} + \beta^2 \fkG_{13524} 
=  \upshift(\fkGO_{(2,3)}).\]
\end{example}

Given   $v=v_1v_2v_3\cdots \in S_\infty$ with $v_j=1$, define 
$ v^{\natural}=v_1^{\natural}v_2^{\natural}v_3^{\natural}\cdots   \in S_\infty$ to be the permutation with
$v^{\natural}_1=2$, $v^{\natural}_{j+1}= 1$, and $v^{\natural}_{i+1} = 1+v_i$ for all $i\neq j$.
If 
\[v=3467125\quad\text{then}
\quad v^{\natural} = 24578136\quand 1\times v = 14578236.\]
Notice that $ v^{\natural} = (1,2)\cdot  (1\times v)= s_1\cdot  (1\times v)$.

\begin{lemma}\label{supp-lem}
Fix integers $1 < i_1 < i_2 <\dots <i_l$ and
let 
\[\textstyle \Pi =  \prod_{j\in[l]}(1+\beta x_{i_j}) 
\quand
\Pi^+ =  \prod_{j\in[l]}(1+\beta x_{1+i_j}) .
\]
Suppose $w \in S_\infty$ has $w(1)=1$. Let $a_{v},b_{v} \in \ZZ[\beta]$ for $v \in S_\infty$ be such that 
\[\textstyle \Pi\cdot  \fkG_w 
= \sum_{v \in S_\infty} a_{v} \cdot \fkG_v
\quand
\Pi^+\cdot  \fkG_{1\times w}    = \sum_{v \in S_\infty} b_{v} \cdot \fkG_v.
\]
Then for all $v \in S_\infty$ it holds that  
$b_{v^{\natural}} =\begin{cases} - a_{v} &\text{if }v(1) \neq 1 \\ 0 &\text{if }v(1)=1.\end{cases}$
\end{lemma}

We point out that
as $\fkG_w$ is homogeneous if we view $\beta$ as having degree $-1$,
the coefficients $a_v$ and $b_v $ in this lemma are always integers times $ \beta^{\ell(v)-\ell(w)}$.
 
\begin{proof}
Since $1\times w$ fixes $1$ and $2$ and since $\Pi^+$ does not involve 
 $x_1$ or $x_2$, it follows from \eqref{fg-eq2} and 
Theorem~\ref{iG-thm} that $\bpartial_1 (\Pi^+\cdot  \fkG_{1\times w}) = -\beta\cdot \Pi^+ \cdot  \fkG_{1\times w}$.
In view of \eqref{Grothendieck-recursion}, this is only possible if $b_u =0$ whenever $u(1)>u(2)$.
The latter condition holds whenever $u=v^\natural$ for some $v\in S_\infty$ with $v(1)=1$.

Next, let $ i = (i_1,i_2,\dots,i_l)$ and $i^+ = (1+i_1,1+i_2,\dots,1+i_l)$.
Define an \defn{$i$-chain} to be an ascending chain in Bruhat order formed by 
starting with a chain of the form \eqref{lenart-eq} with $k=i_l$,
then appending chains of the same form with $k=i_{l-1}, i_{l-2},\dots,i_1$ in succession.
We say that an $i$-chain is \defn{odd} (respectively, \defn{even}) if the sum of the values of $p$ from \eqref{lenart-eq}
corresponding to the $l$ segments of the chain  is odd (respectively, even).

Finally, fix $v \in S_\infty$ with $v(1)\neq 1$.
If we use Theorem~\ref{lenart-transition} to multiply $\fkG_{1\times w}$ first by $(1+\beta x_{i_l})$,
then by $(1+\beta x_{i_{l-1}}) $, then by $(1+\beta x_{i_{l-2}})$, and so on, to ultimately obtain the $\fkG$-expansion of $\Pi^+\cdot  \fkG_{1\times w}$, then  we see that
the coefficient $b_{v^\natural}$ is exactly the difference between the numbers of even and odd $i^+$-chains from $1\times w$ to $v^\natural$.
Similarly, the coefficient $a_v$ in the $\fkG$-expansion of $\Pi\cdot \fG_w$ is  
the difference between the numbers of even and odd $i$-chains from $w$ to $v$.

The key idea of this proof rests in the following observations.
Because $1\times w$ fixes $1$ and $2$ while $v^\natural(1)=2$, in any $i^+$-chain $\gamma$ from $1\times w$ to $v^\natural$
there must be an ascent of the form $\xrightarrow{(2,i_j)}$ (as otherwise we would have $\last(\gamma)(1)=1$) and the first such ascent must be followed 
immediately an ascent of the form $\xrightarrow{(1,i_j)}$ (as otherwise $2$ would never return to the first position in the one-line representation of $\last(\gamma)$),
and then all other ascents in the chain $\xrightarrow{(a,b)}$ must have $a>1$ (as any such ascent with $a=1$ must come later in the chain, but if there is such a later ascent then we would have $\last(\gamma)(1)>2$).
In view of these properties, we deduce that removing $\xrightarrow{(1,i_j)}$
and then replacing all other ascents $\xrightarrow{(a,b)}$ by $\xrightarrow{(a-1,b-1)}$
transforms $\gamma$ to an $i$-chain from  $w$ to $v$. Conversely, by considering the natural inverse operation we see that every  $i$-chain from  $w$ to $v$ arises in this way for a unique choice of $\gamma$. 

As the bijection just given between $i^+$-chains from $1\times w$ to $v^\natural$ and $i$-chains from $w$ to $v$
reverses parity, we conclude that $b_{v^\natural} = -a_v$ as needed.
\end{proof}

We are now ready to prove Theorem~\ref{supp-thm}.

\begin{proof}[Proof of Theorem~\ref{supp-thm}]
Assume $z(1) \neq 1$ and let $L$ be the left segment of $1\times z$ containing $2$,
so that $\{i -1 : i \in L\}$ is the unique immobile left segment of $z$.
We first check that $\upshift(\fkGO_z)\neq \fkGO_{1\times z}  $.

Sending $S \mapsto 1 +S$ gives a bijection from $\shiftable(z)$ to the family of shiftable sets $T \in \shiftable(1\times z)$ with $T\cap L = \varnothing$.
Since $ \customvarpi_{ z,S} =  \customvarpi_{1\times z,1+S}$ and  $\upshift(\iG_{z , S}) = \iG_{1\times z,1+S}$,
it follows from Theorem~\ref{ivex-thm} that
\[  \fkGO_{1\times z} - \upshift(\fkGO_z) =  \sum_{\substack{S\in\shiftable(1\times z) \\ S\cap L \neq \varnothing}}  \beta^{|S|} \cdot  \customvarpi_{1\times z,S} \cdot \iG_{1\times z , S}.
\]
This difference is divisible by $\beta$. Moreover, first dividing by $\beta$ and then setting $\beta=0$
kills all terms in the sum indexed by sets $S$ with $|S| \neq 1$.
Consulting \eqref{varpi-eq}, we see that this operation
 transforms $\fkGO_{1\times z} - \upshift(\fkGO_z)$ to the sum
\[\textstyle
 \sum_{n \in L}   \customvarpi_{1\times z,\{n\}}\big|_{\beta =0}  \cdot \iG_{1\times z , \{n\}}\big|_{\beta =0}
=
2^{\cyc(z)-1} \sum_{n \in L}    \iS_{t_n (1\times z)t_n}
 \]
where   $t_n = (1,n) \in S_\infty$ and  $\iS_y := \iG_y\big|_{\beta=0}$ is the \defn{involution Schubert polynomial} considered in \cite{HMP1}.
This sum cannot be zero as involution Schubert polynomials are nonzero elements of $\NN[x_1,x_2,\dots]$,
so  $\fkGO_{1\times z} - \upshift(\fkGO_z) \neq 0$.

Now assume instead that $z(1)=1$. Then the formula for $\fkGO_{1\times z}$ in Theorem~\ref{ivex-thm}
expands as a sum of terms of the form $ \Pi\cdot \fkG_w$ as in Lemma~\ref{supp-lem},
and the corresponding formula for $\fkGO_{1^2\times z}$
is obtained by replacing each term $ \Pi\cdot \fkG_w$ by $\Pi^+\cdot \fkG_{1\times w}$.
It therefore follows from Lemma~\ref{supp-lem} that if $v \in S_\infty$ has $v(1)\neq 1$
then $\GC_{1^2\times z}(v^\natural) = -\GC_{1\times z}(v)$.
As $\GC_{1\times z}$ and $\GC_{1^2\times z}$ both only take nonnegative values \cite[Prop.~3.29]{MS2023},
we must have $\GC_{1\times z}(v) = 0$ if $v(1)\neq 1$.
We conclude from Proposition~\ref{yz-down-prop} and \eqref{shift-eq} that
\[\fkGO_{1\times z} = \upshift\circ\downshift(\fkGO_{1\times z}) = \upshift(\fkGO_z)\]
as claimed.
\end{proof}

Using Theorem~\ref{supp-thm}, we can now prove another nontrivial result.
 Comparing this with Proposition~\ref{supp-prop}
 gives some heuristic support for Conjecture~\ref{b+conj}.
  
\begin{proposition}
Suppose $z \in \Ivex$ has   $\supp(z)  \subseteq [a,b]$ for some integers $a\leq b$. Then $\supp(w) \subseteq [a-1,b+1]$
for all  $w \in \supp(\GC_z)$.
\end{proposition}

\begin{proof}
Define $p$ and $q$ as in \eqref{pq-eq}, so that $q \leq b$, and let $\delta = \dom_{pq}$ be the dominant involution defined in \eqref{dom-def}.
It follows from Corollary~\ref{qd-supp-cor} and Proposition~\ref{supp-prop}
that $\supp(w) \subseteq [1,q+1]$ for all $w \in \supp(\GC_{\delta})$.
Since $\fkGO_z$ can be obtained from $\fkGO_{\delta}$ by applying a sequence of divided difference 
operators $\bpartial_i$ with $i <q$,
we also have $\supp(w) \subseteq[1,q+1]\subseteq [1,b+1]$ for all $w \in \supp(\GC_z)$.

It remains to check that $\supp(w)\subseteq [a-1,\infty)$ for all $w \in \supp(\GC_z)$.
This is trivial unless $a\geq 3$,
and then the desired containment follows from
the fact that
$\fkGO_z = \upshift^{a-2}(\fkGO_y)$
 for $y := z\downarrow (a-2) \in \Ivex$
by Theorem~\ref{supp-thm}.
\end{proof}

For $z \in I_\infty$ we 
 extend the domain of the Grothendieck coefficient function $\GC_z$ by setting $ 
 \GC_z(w) = 0$ if $w \notin S_\infty$.
The following identities
are immediate from
Proposition~\ref{yz-down-prop} and
Theorem~\ref{supp-thm}.

\begin{corollary}\label{shift-cor}
 Fix $n \in \NN$ and $z \in \Ivex$.
\ben
\item[(a)] If $z \downarrow n \in S_\infty$ then $\GC_{z\downarrow n }(w) = \GC_z(1^n \times w)$ for all $w \in S_\infty$.

\item[(b)] If $z(1)=1$ then $\GC_{1^n\times z }(w) = \GC_z(w \downarrow n)$ for all $w \in S_\infty$.
\een
\end{corollary}


The properties in this corollary are consistent with Conjecture~\ref{b+conj} since 
it is easy to show that
if $n \in \NN$ and  $z \in \Ivex$ has $z(1)=1$ then
\be
\label{BB-eq}
\ba
\cB(1^n\times z ) &= \left\{ w \in S_\infty : w \downarrow n \in \cB(z)\right\},
\\
\cB^+(1^n \times z ) &= \left\{ w \in S_\infty : w \downarrow n \in \cB^+(z)\right\}.
\ea
\ee
In fact, the identity for $\cB(1^n\times z ) $ holds even if $z(1) \neq 1$.

Theorem~\ref{supp-thm} has an application concerning the \defn{stable limit} of $\fkGO_z$.
The \defn{symmetric Grothendieck function} of $w \in S_\infty$ 
is
\be\textstyle G_w := \lim_{n\to \infty} \fkG_{1^n\times w}\ee where the limit is taken in the sense of formal power series. This means the limit exists precisely when the coefficients of any fixed monomial 
in $\fkG_{1^n\times w}$ is eventually a constant sequence.

It is known \cite{Buch2002,BuchEtAl,FominKirillov} that  $G_w$ always exists and is a formal power series that is symmetric in the $x_i$ variables. Moreover, every symmetric power series over $\ZZ[\beta]$ is given by a unique, possibly infinite $\ZZ[\beta]$-linear combination
of $G_w$'s \cite[Thm.~6.12]{Buch2002}.
Consulting Remark~\ref{bcs-rem}, one obtains the explicit generating function 
$G_w = \sum_{(a,i)} \beta^{\ell(a)-\ell(w)} x^i$
where the sum is over all (not necessarily bounded) compatible sequences $(a,i)$ with $a \in \cH(w)$.

Following \cite{MP2020}, for $z \in \Ivex$ we define an analogous power series limit
\be\label{gq-eq}\textstyle \GQ_z := \lim_{n\to \infty} \fkGO_{1^n\times z}.\ee
We can immediately see that this limit converges to a symmetric function:

\begin{corollary}\label{supp-cor}   If $z \in \Ivex$  has $z(1)=1$ then 
\[\textstyle \GQ_z = \sum_{w \in S_\infty}  \beta^{\ell(w) - \ellhat(z)}  \cdot \GC_{z}(w) \cdot   G_w.\]
Moreover,   it holds for any $z \in \Ivex$ that $\GQ_z = \GQ_{1\times z}$.
\end{corollary}

\begin{proof}
If $z(1)=1$ then 
$\fkGO_{1^n\times z} = \sum_{w \in S_\infty}\beta^{\ell(w) - \ellhat(z)}  \cdot \GC_{z}(w) \cdot   \fkG_{1^n\times w}$ holds by Theorem~\ref{supp-thm}
and the corollary follows by taking the limit as $n\to \infty$.
\end{proof}

\begin{remark}
Comparing Corollary~\ref{supp-cor} with \cite[Lem.~3.7 and Prop.~3.29]{MS2023} shows that the conjecture \cite[Conj.~3.39]{MS2023} holds when $z$ is vexillary.
\end{remark}

A stronger result about $\GQ_z$ is known: if $z \in \Ivex$ then $\GQ_z$ is equal to the \defn{$K$-theoretic Schur $Q$-function} $\GQ_\mu$ defined in  \cite{Ikeda.Naruse2013} when $\mu=\sh(z)$ is the \defn{involution shape} of $z$ \cite[Thm.~4.11]{MP2020}.
The involution shape   $\sh(z)$ is the transpose of the partition that sorts the sequence of numbers 
$\hat c(z)$ 
from  \eqref{hat-d-eq}.
One can view Corollary~\ref{supp-cor} as a new $G$-positive expansion of $\GQ_\mu$.

%

 It is an open problem to show that 
\eqref{gq-eq} converges when $z \in I_\infty$ is not vexillary.
By contrast, it is clear from \eqref{BB-eq}  that for any $z \in I_\infty$ the limit
\be\label{gp-eq1}
\textstyle \GP_z  := \lim_{n\to \infty}  \iG_{1^n\times z }
\ee
converges in the sense of formal power series to the symmetric function 
\be\label{gp-eq2}
\textstyle \GP_z  =\sum_{w \in \cB(z)} \beta^{\ell(w)-\ellhat(z)} G_w.
\ee
 In fact, it is known \cite[Thm.~1.9]{Mar2019} that the power series $\GP_z$ are all $\NN[\beta]$-linear combinations of the
\defn{$K$-theoretic Schur $P$-functions} $\GP_\mu$ defined in \cite{Ikeda.Naruse2013}.

We can derive one other formula for $\GQ_z$ using this notation and Theorem~\ref{ivex-thm}.
Given $z \in \Ivex$
and $S \in \shiftable(z)$,
define $\GP_{z , S} := \GP_{y}$ where $y:=(\sigma_S)^{-1}\cdot z \cdot \sigma_S$
and form $ \Theta_{z,S} \in \ZZ$ from $ \customvarpi_{z,S}$ by setting $x_1=x_2=\dots=0$. 
More explicitly, if $L= \{ a\in \PP : a<z(a)\}$ then
\[\Theta_{z,S} = (-1)^{|\{ a \in S \hs : \hs \crb_z(a) \cap S \neq \varnothing\}|} 2^{|L|-| S|} = \pm 2^{\cyc(z)-| S|}. \]
The same argument as in the proof of Proposition~\ref{locally-nc-prop} shows 
that $\Theta_{z,S}>0$ for all $S \in \shiftable(z)$ if and only if $z\in \Ivex$ 
is locally noncrossing.

\begin{corollary}\label{ivex-cor}
If $z \in \Ivex$ has $z(1)=1$ then 
\[\textstyle\GQ_z = \sum_{S\in\shiftable(z)} \beta^{|S|} \cdot \Theta_{z,S} \cdot   \GP_{z , S}.\]
\end{corollary}

 \begin{proof}
Form $\customvarpi_{z,S}^{\uparrow n}$ from $\customvarpi_{z,S}$ by applying the $\ZZ[\beta]$-algebra homomorphism that sends $x_i \mapsto x_{i+n}$ for all $i \in \PP$.
Theorem~\ref{ivex-thm} implies that if $z(1)=1$ then 
$\fkGO_{1^n\times z} = \sum_{S\in\shiftable(z)}   \beta^{|S|} \cdot\customvarpi_{z,S}^{\uparrow n} \cdot  \iG_{1^n\times (\sigma_S)^{-1}\cdot  z \cdot \sigma_S}$.
The result follows by taking the limit  $n\to \infty$, which turns $ \customvarpi_{z,S}^{\uparrow n}$ into $ \Theta_{z,S}$ and $\iG_{1^n\times (\sigma_S)^{-1}\cdot  z \cdot \sigma_S}$ into $\GP_{z , S}$.
\end{proof}

\subsection{Grassmannian identities}\label{gr-sect}

There is an interesting special case when the formulas for $\fkGO_z$ and $\GQ_z$ in the previous sections become considerably more explicit. We explain this below after defining some additional notation.

 Recall the definition of the visible descent set $\DesV(z)$ from \eqref{vdes-def}.
 Following \cite{HMP3}, we say  $z \in I_\infty$ is \defn{$n$-I-Grassmannian} 
$\DesV(z)\subseteq \{n\}$.
An element of $I_\infty$ is \defn{I-Grassmannian} if it is $n$-I-Grassmannian for some $n \in \NN$.

For a \defn{strict} partition $\mu=(\mu_1>\mu_2>\dots>\mu_r>0)$ with $r>0$ nonzero parts and largest part $\mu_1\leq n$, define $\langle \mu \mid n\rangle \in I_\infty$
to be the permutation 
\be  (n+1-\mu_1, n+1)(n+1-\mu_2,n+2)\cdots(n+1-\mu_r,n+r) \in I_\infty.\ee
When $\mu=\emptyset$ is empty let $\langle \mu \mid n\rangle =1$ for all $n \in \PP$.
 This notation deliberately mimics our convention for writing $n$-Grassmannian permutations $[\lambda\mid n]\in S_\infty$, for reasons that shall soon be apparent.

For each strict partition $\mu$ with largest part $\mu_1\leq n$,
 the involution $\langle \mu \mid n \rangle$ turns out \cite[\S4.1]{HMP3} to be the unique
 $n$-I-Grassmannian element of $I_\infty$ with involution shape $\sh(\langle \mu \mid n \rangle)=\mu$.
From this observation and Remark~\ref{arc-vex-rem},
we deduce that every $I$-Grassmannian involution is vexillary.
Thus, we can apply
Theorem~\ref{ivex-thm}  when $z$ is $I$-Grassmannian.

Let $\SD_\mu := \{ (i,i+j-1) : (i,j) \in \D_\mu\}$ denote the \defn{shifted Young diagram}
of a strict partition $\mu$. Write $\mu\subseteq\lambda$ if $\mu$ and $\lambda$ are strict partitions
with $\SD_\mu\subseteq \SD_\lambda$.
When this containment holds we define 
\[
\SD_{\lambda/\mu} := \SD_\lambda\setminus \SD_\mu
\quand |\lambda/\mu| := |\lambda|-|\mu|
\quad\text{where $|\mu| := |\SD_\mu| $.}\]

The shifted skew shape $\SD_{\lambda/\mu}$ is a \defn{vertical strip} if it contains at most one position in each row,
or equivalently if $\lambda_i-\mu_i \in \{0,1\}$ for all $i \in \PP$.
Assume this holds and let $r$ be the number of nonzero parts of $\lambda$.
Define $\topRows(\lambda/\mu)$ to be the set of indices $i \in [r]$
such that row $i$ of $\SD_{\lambda/\mu}$ is empty or 
such that $(i+1,j) \notin \SD_{\lambda/\mu}$ for the unique column $j$ with $(i,j) \in \SD_{\lambda/\mu}$.
Then let 
\be
\textstyle
\customvarpi_{\lambda/\mu}^{(n)}  :=
\prod_{i \in\topRows(\lambda/\mu) } (2 +\mu_i- \lambda_i +\beta x_{n+1-\mu_i}).
\ee
Also write $\cols(\lambda/\mu) = |\{ j : (i,j) \in \SD_{\lambda/\mu}\text{ for some }i \in \PP\}|$
for the number of distinct columns occupied by the positions in $\SD_{\lambda/\mu}$.

\begin{example}
If $\mu =(8,7,4,3,1)$ and  $\lambda = (9,7,5,4,2)$ then 
\[
\ytableausetup{boxsize=0.3cm,aligntableaux=center}
\SD_{\lambda/\mu} = \begin{ytableau}
\none[\cdot] & \none[\cdot] & \none[\cdot] & \none[\cdot] & \none[\cdot] & \none[\cdot] & \none[\cdot] & \none[\cdot] & \ \\
\none  & \none[\cdot] & \none[\cdot] & \none[\cdot] & \none[\cdot] & \none[\cdot] & \none[\cdot] & \none[\cdot]  \\
\none  & \none  & \none[\cdot] & \none[\cdot] & \none[\cdot] & \none[\cdot] & \  \\
\none  & \none  & \none  & \none[\cdot] & \none[\cdot] & \none[\cdot] & \ \\
\none  & \none  & \none   & \none & \none[\cdot] & \ 
\end{ytableau}
\]
so $|\lambda/\mu| = 4$, $\cols(\lambda/\mu) =3$, $\topRows(\lambda/\mu) = \{1,2,4,5\}$, and 
\[
\customvarpi_{\lambda/\mu}^{(9)} = (1+ \beta x_{2})(2+\beta x_{3}) (1+\beta x_{7})(1+\beta x_9).
\]
\end{example}

\begin{corollary}\label{igrass-cor}
If $\mu$ is a strict partition with largest part at most $n$
then
\[
\fkGO_{\langle \mu \mid n \rangle} = 
 \sum_\lambda   (-\beta)^{|\lambda/\mu|} \cdot  (-1)^{\cols(\lambda/\mu)}  \cdot \customvarpi_{\lambda/\mu}^{(n)} \cdot\iG_{\langle \lambda \mid n\rangle}
\]
where the sum is over strict partitions $\lambda$ with the same number of parts as $\mu$ 
and with largest part at most $n$, such that $\lambda \supseteq\mu$ and $\SD_{\lambda/\mu}$ is a vertical strip.
\end{corollary}

\begin{proof}
Let $z = \langle\mu \mid n \rangle$. Then
as $S$ ranges over the elements of $\shiftable(z)$,
the permutations $ (\sigma_S)^{-1}\cdot z\cdot  \sigma_S$
range over exactly the involutions $\langle \lambda\mid n\rangle$ for $\lambda$
as described in the corollary. 

Additionally, when going across the correspondence $S\leftrightarrow \langle \lambda\mid n\rangle$,
one has $|S| = |\lambda/\mu|$, 
while the number of factors $\customvarpi_{z,S}^{(a)}$ equal to $-1$ is $|\lambda/\mu| - \cols(\lambda/\mu)$
and the product of the other factors $\customvarpi_{z,S}^{(a)}$ is exactly $\customvarpi_{\lambda/\mu}^{(n)}$.

Verifying these claims is a straightforward exercise, and then the desired formula
for $\fkGO_{\langle \mu \mid n\rangle}$
is immediate from Theorem~\ref{ivex-thm}.
\end{proof}

\begin{example}
If we fix $n=4$ and write $\langle \mu_1\mu_2\cdots\mu_r \rangle$ instead of $\langle \mu\mid 4\rangle $ then
\[
\fkGO_{\langle 32  \rangle} = 
(2+\beta x_{2})(2+\beta x_{3})\iG_{\langle 32  \rangle}
+
\beta (1+\beta x_{2})(2+\beta x_{3})  \iG_{\langle 42  \rangle}
-
\beta^2(1+\beta x_{3})  \iG_{\langle 43  \rangle}.\]
\end{example}

Corollary~\ref{igrass-cor} is a non-symmetric
generalization of \cite[Thm.~1.1]{Chiu.Marberg2023},
which one recovers by taking stable limits in the following way.
Let $\mu$ be a strict partition with $\mu_1\leq n$.
Observe that $1^m\times  \langle \mu \mid n \rangle  = \langle \mu \mid m+n\rangle$ for all $m \in \NN$.
Therefore 
\be
\fkGO_{1^m\times \langle \mu \mid n \rangle} = 
 \sum_\lambda (-\beta)^{|\lambda/\mu|} \cdot (-1)^{\cols(\lambda/\mu)}  \cdot \customvarpi_{\lambda/\mu}^{(m+n)} \cdot \iG_{1^m\times\langle \lambda \mid n\rangle }
\ee
where the sum is over strict partitions $\lambda$ as in Corollary~\ref{igrass-cor}.

We have $\lim_{m\to \infty} \customvarpi_{\lambda/\mu}^{(m+n)} = 2^{|\{i \in [r] : \lambda_i=\mu_i\}|} = 2^{\ell(\mu) - |\lambda/\mu|}$ with convergence in the sense of formal power series,
 where $\ell(\mu)$ denotes the number of nonzero parts of $\mu$.
Hence, taking the limit as $m\to \infty$ of both sides above gives
\be\label{almost-eq}
\GQ_{\langle \mu \mid n \rangle} = 
 \sum_\lambda  (-\beta)^{|\lambda/\mu|} \cdot(-1)^{\cols(\lambda/\mu)}  \cdot 2^{\ell(\mu)-|\lambda/\mu|}\cdot  \GP_{\langle \lambda \mid n\rangle}
\ee
where the sum is over strict partitions $\lambda$ as in Corollary~\ref{igrass-cor}.
This identity is nearly the same as \cite[Thm.~1.1]{Chiu.Marberg2023},
which asserts the same equation but with ``$\GQ_{\langle \mu \mid n \rangle}$''
and ``$\GP_{\langle \lambda \mid n \rangle}$''
replaced by 
``$\GQ_{\mu}$''
and ``$\GP_{ \lambda}$''.

The subtle difference between these statements, which is obscured by our notation,
is that $\GQ_{\langle \mu \mid n \rangle}$
and $\GP_{\langle \lambda \mid n \rangle}$
are instances of 
the stable limits \eqref{gq-eq} and \eqref{gp-eq1},
while $\GQ_\lambda$ and $\GP_\lambda$ are the symmetric functions defined in a quite different way in \cite{Ikeda.Naruse2013}. 
However, as mentioned earlier, it is actually already known \cite[Thm.~4.11]{MP2020} that $\GQ_{\langle \mu \mid n \rangle} = \GQ_{\sh(\langle \mu \mid n \rangle)} = \GQ_\mu$.
To complete the circle, it would suffice to have the analogous identity $
\GP_{\langle \mu \mid n \rangle} = \GP_\mu
$. This does not seem to have appeared in the literature,
but we can prove it here:

\begin{proposition}
If $\mu$ is a strict partition and $n\geq \mu_1$ 
then 
$
\GP_{\langle \mu \mid n \rangle} = \GP_\mu
$.
\end{proposition}

\begin{proof}
Our argument is very similar to the proofs of \cite[Prop.~4.14]{MP2021Sp}
and
\cite[Prop.~3.37]{MS2023}.
Suppose $\mu$ has $r$ nonzero parts, so that 
$\langle \mu \mid n \rangle $ is the involution $ (\phi_1 , n+1)(\phi_2 , n+2) \cdots (\phi_r , n+r)$ where $\phi_i = n+1-\mu_i$.
Then
\[
\textstyle
x^{\mu}  \prod_{i=1}^{r}  \prod_{j=i+1}^n \tfrac{x_i\oplus x_j}{x_i} = x_1^{1-\phi_1}x_2^{2-\phi_2}\cdots x_r^{r-\phi_r} \iG_w
\]
for the dominant permutation $w := (1 , n+1)(2 , n+2)\cdots (r , n+r)$,
since  by Theorem~\ref{iG-thm2} it holds that $\iG_w = \prod_{i=1}^r x_i\cdot \prod_{i=1}^r \prod_{j=i+1}^n x_i\oplus x_j$.

Define $\bpi_i f  := \bpartial(x_i f)$ for $i \in \PP$; see \cite[\S4.2]{MP2021Sp} for the basic properties of these 
\defn{isobaric divided difference operators}.
For integers $0 <a \leq b$ let \be\bpi_{b \searrow a} := \bpi_{b-1}\bpi_{b-2}\cdots \bpi_a\ee
so that $\bpi_{a\searrow a} = 1$, and define $\bpartial_{b\searrow a}$ analogously.
Exactly as in the proofs of \cite[Prop.~4.14]{MP2021Sp}
or
\cite[Prop.~3.37]{MS2023}, one can check by an easy calculation that 
\[ \bpi_{\phi_j \searrow j} \bpi_{\phi_{j+1}\searrow j+1} \cdots \bpi_{\phi_r \searrow r} \( x_1^{1-\phi_1}x_2^{2-\phi_2}\cdots x_r^{r-\phi_r} \iG_w\)\] is fixed by $\bpi_i$ for all $j -1 < i < \phi_{j-1} < \phi_j$,
and so  is equal to 
\[
\bpartial_{\phi_1 \searrow 1} \bpartial_{\phi_{2}\searrow 2} \cdots \bpartial_{\phi_r \searrow r} \iG_w = \iG_{\langle\mu\mid n\rangle}
\]
by \cite[Lem.~4.16]{MP2021Sp} and Theorem~\ref{iG-thm}.

Let 
$\mathsf{stab}_n = \bpi_{2\searrow 1} \bpi_{3\searrow 1}  \bpi_{4\searrow 1}\cdots \bpi_{n\searrow 1}.$
Since every $w \in \cB(z)$ for $z \in I_\infty$ has $\DesR(w)\subseteq \DesR(z)$ by \eqref{iwi-eq},
it follows from  \cite[Lem.~3.7]{MS2023}
that  
\[ G_w(x_1,x_2,\dots,x_n) =\mathsf{stab}_n( \fkG_w)\quad\text{for all }w \in \cB(\langle \mu \mid n\rangle).\]
Next, it is known \cite[Eq.~(3.17)]{MS2023} that
\[\GP_\mu(x_1,x_2,\dots,x_n) = \mathsf{stab}_n \( \textstyle x^{\mu}  \prod_{i=1}^{r}  \prod_{j=i+1}^n \tfrac{x_i\oplus x_j}{x_i} \).\]
Finally, one has  $\mathsf{stab}_n \pi_i = \mathsf{stab}_n$ if $i \in [n-1]$ since the $\pi_i$ operators are idempotent and
 obey the braid relations for $S_\infty$ \cite[\S4.2]{MP2021Sp}.
Combining these observations with \eqref{gp-eq2} shows that
\[\GP_\mu(x_1,x_2,\dots,x_n) =  \mathsf{stab}_n (\iG_{\langle \mu\mid n\rangle})=\GP_{\langle \mu\mid n\rangle}(x_1,x_2,\dots,x_n).\]
If $N\geq n$ then  $\GP_{\langle \mu \mid n \rangle} =\GP_{1^{N-n}\times \langle \mu \mid n \rangle} = \GP_{\langle \mu \mid N \rangle} $ so  \[\GP_\mu(x_1,x_2,\dots,x_N)  = \GP_{\langle \mu \mid n\rangle}(x_1,x_2,\dots,x_N ).\]
On taking the limit $N \to\infty$ we conclude that $\GP_\mu = \GP_{\langle \mu \mid n\rangle}$.
\end{proof}
 
There is one case where we can turn Corollary~\ref{igrass-cor} into an explicit positive formula.
Define $ g_{ij} = (i, j+1)(i+1,j+2)(i+3,j+4)\cdots(j,2j-i+1)$
for any integers $0<i<j$. This coincides with the $I$-Grassmannian involution
$\langle \mu \mid j\rangle$ for $\mu = (m,m-1,m-2,\dots,1)$ and $m=j-i+1$.
By \cite[Thm.~3.35 and Cor.~3.36]{HMP1}, the involutions $g_{ij}$ are the only elements in $I_\infty$ that are both vexillary and 321-avoiding.

Notice that $1^n \times g_{ij} = g_{i+n,j+n}$. 
Thus, to compute $\cB(g_{ij})$, $\cB^+(g_{ij})$, and $\fkGO_{g_{ij}}$ for arbitrary integers $0<i<j$, it suffices by the results in Section~\ref{shift-sect}
to consider just the case when $i=2$.
For this purpose, given any  $w\in S_\infty$ recall that $\DesL(w) = \DesR(w^{-1})$ and let $\ODes_L(w) := \{ i \in \DesL(w) : i\text{ is odd}\}$.

\begin{proposition}\label{g-prop}
Assume $n>2$ and define $z \in \Ivex$ to be the involution
 \[z=  g_{2n} =\langle(n-1,\dots,3,2,1)\mid n\rangle= (2,n+1)(3,n+2)\cdots(n,2n-1).\] Then the following properties hold:
   \ben
    \item[(a)] $\cB(z)$ contains only the inverse of
    $1(n+1)2(n+2)3\cdots (2n-1)n(2n).$
    
    \item[(b)] $\cB^+(z)$  consists of the $2^{n}$ permutations   whose inverses in one-line notation 
    have the form $a_1b_1a_2b_2\cdots a_{n} b_{n} $
    where $\{a_i,b_i\}=\{i,n+i\}$ for each $i\in [n]$.

\item[(c)] If $u_{\max} \in \cB^+(z)$ is the inverse of $(n+1)1(n+2)2\cdots (2n)n$
 then
 \[\textstyle\fkGO_{z} =   \sum_{w \in \cB^+(z)}  \beta^{|\ODes_L(w)|}  \cdot 2^{n-1-|\ODes_L(w)|} \cdot  \fkG_w +
(-\beta)^n \cdot\tfrac{1}{2} \cdot  \fkG_{u_{\max}}
\]
and consequently
 \[\GC_z(w) = \begin{cases} 2^{n-1-|\ODes_L(w)|}&\text{if }u_{\max} \neq w \in \cB^+(z) \\
1 &\text{if $u_{\max}=w$ and $n$ is even} \\ 
0 &\text{if $u_{\max}=w$ and $n$ is odd, or $w \notin \cB^+(z)$}.
\end{cases}\]

    \een

\end{proposition}

See Figure~\ref{g2n-fig} for an illustration of this proposition.

\begin{proof}
Part (a) is immediate from \eqref{sim-eq} and part (b) follows as an easy exercise from the definitions. We omit the details.

We may prove part (c) in the following way. 
First, it is straightforward using Remark~\ref{arc-vex-rem} to check that $s_izs_i$ and $z$ are vexillary for $i\in [n-1]$.
%
%
    Because $\Des_R(z) = \Des_V(z) =\{n\}$, we have $\fkGO_z\in \ZZ[\beta][x_1,x_2,\cdots,x_n]$ 
    by \eqref{vdes-eq}
    and $\bpartial_i(\fkGO_z) =-\beta \fkGO_z$ for all $i\in [n-1]$ by \eqref{orthogonal-recursion}.
 The latter identity implies that $\fkGO_z$ is symmetric in $x_1,x_2,\cdots,x_n$. 
By Remark~\ref{grass-sym-rem},
 we therefore have 
\[\textstyle
\fkGO_z= \sum_{\mu}\beta^{|\mu| - \ellhat(z)} \cdot \GC_z([\mu\mid n]) \cdot  \fkG_{[\mu\mid n]}
\]
where the sum is over partitions $\mu$ with at most $n$ nonzero parts.
On the other hand,
as $z\in \Ivex$ has $z(1)=1$, Corollary~\ref{supp-cor} implies that
\[\textstyle\GQ_z=\sum_{\mu} \beta^{|\mu| - \ellhat(z)} \cdot\GC_z([\mu\mid n]) \cdot   G_{[\mu\mid n]}.\]
Since
the polynomial 
 $G_{[\mu\mid n]}(x_1,x_2,\cdots,x_n)$ obtained from $G_{[\mu\mid n]}$ 
 by setting the variables $x_{n+1}=x_{n+2}=\dots=0$ is equal to
 $ \fkG_{[\mu| n]}$ for all partitions $\mu$ with at most $n$ parts (compare
\cite[Thm.~6.12]{Buch2002} and \cite[Thm.~2.2]{lenart2000}), we have
\be\label{GQ-suff-eq}
\textstyle
\GQ_z(x_1,x_2,\cdots,x_n) = \sum_{\mu}  \beta^{|\mu| - \ellhat(z)} \cdot \GC_z([\mu\mid n]) \cdot \fkG_{[\mu\mid n]}.
\ee Hence,
 it suffices to compute 
the $\fkG$-expansion of $\GQ_z(x_1,x_2,\cdots,x_n)$.

We can accomplish this by assembling some miscellaneous identities from \cite{Chiu.Marberg2023,Ikeda.Naruse2013,LM2019,MP2020}.
 First, \cite[Thm.~4.11]{MP2020} tells us that $\GQ_z$ is equal to the $K$-theoretic Schur $Q$-function indexed by $\rho_{n-1}=(n-1,n-2,\dots,3,2,1)$, which can be expanded in terms of $K$-theoretic Schur $P$-functions via \cite[Thm.~1.1]{Chiu.Marberg2023} as
   \[
    \textstyle
  \GQ_z=   2^{n-1}\GP_{\rho_{n-1}}-2^{n-1}\sum_{j=1}^{n-1}(-\beta/2)^j \GP_{\rho_{n-1}+1^j}.
\]
    Next, it follows from \cite[Prop.~2.3]{Ikeda.Naruse2013} that
    \[
    \GP_{\rho_{n-1}+1^j}(x_1,x_2,\cdots,x_n) = \GP_{\rho_{n-1}}(x_1,x_2,\cdots,x_n)G_{1^j}(x_1,x_2,\cdots,x_n).
    \]
Since by \cite[Thm.~6.8]{LM2019} we have $\GP_{\rho_{n-1}} = G_{\rho_{n-1}}$,  
we can rewrite this as \[
    \GP_{\rho_{n-1}+1^j}(x_1,x_2,\cdots,x_n) 
    = \fkG_{[\rho_{n-1}| n]}\fkG_{[1^j\mid n]}.
    \]
    If $I$ is a subset of $[n]$ then we let $\mu_{I}= \rho_{n-1}+\sum_{i\in I} \e_{i}$.
Theorem~\ref{lenart-pieri} says that 
\[\textstyle
 \fkG_{[\rho_{n-1}| n]}\fkG_{[1^j\mid n]}=
\sum_{p=j}^n \sum_{|I|=p}\beta^{p-j}
\cdot \tbinom{p-1}{p-j}
\cdot  \fkG_{[\mu_{I}\mid n]}.
\]
%
Thus,
if
 $|I|=p< n$ then 
 the coefficient of $\fkG_{[\mu_{I}\mid n]}$ in  $  \GQ_z(x_1,x_2,\dots,x_n)$ 
 is
\[
  -  2^{n-1}\sum_{j=1}^p \beta^{p-j} (-\tfrac{\beta}{2})^{j}   \tbinom{p-1}{p-j}=2^{n-2} \beta^p \underbrace{\sum_{j=0}^{p-1} (-\tfrac{1}{2})^{-j} \tbinom{p-1}{p-1-j}}_{=(-\frac{1}{2} + 1)^{p-1} = 2^{1-p}} = 2^{n-1-|I|}\beta^{|I|}.
\]
If $I=[n]$ then $\mu_{I} = \rho_n$ and 
the coefficient of $\fkG_{[\rho_n\mid n]}$ in $  \GQ_z(x_1,x_2,\dots,x_n)$ 
is  
\[
   - 2^{n-1}\sum_{j=1}^{n-1} \beta^{n-j} (-\tfrac{\beta}{2})^{j} \tbinom{n-1}{n-j}=2^{n-2} \beta^n\big(
       2^{1-n} -(-2)^{1-n}  
        \big)= \begin{cases}
     \beta^n &\text{ $n$ is even}\\ 
     0 &\text{ $n$ is odd.}\\ 
    \end{cases}
\]
The desired formula for $\GC_z$ follows from these computations via \eqref{GQ-suff-eq}
 since the elements of $\cB^+(z)$
described in part (b) are precisely the Grassmannian permutations $[\mu_{I}\mid n]$ for $I \subseteq[n]$, and it holds that 
 $u_{\max}=[\rho_n\mid n]$ along with
\[\ell([\mu_I \mid n]) - \ellhat(z) = |\mu_I| - |\rho_{n-1}| = |I| = \ODes([\mu_{I}\mid n]).\]
This completes our proof of the proposition.
\end{proof}

Using Corollary~\ref{shift-cor}, we see 
from this proposition that 
 Conjecture~\ref{b+conj} holds if $z=g_{ij}$ for any integers $0<i<j$.
 In fact, we have $\supp(\GC_{g_{ij}})= \cB^+(g_{ij})$
 unless $i>1$ and $j-i$ is odd, in which case $ \cB^+(g_{ij})$ contains one extra element.


\begin{figure}[h]
\centerline{
\resizebox{\textwidth}{!}{%
\begin{tikzpicture}[xscale=2, yscale=1.875,>=latex,baseline=(z.base)]
\node at (0,1.5) (z) {};
  \node[draw,fill=lightgray] at (1,0) (11) {\tiny $415263^{-1}:0$};
  \node[draw] at (0,1) (21) {\tiny $145263^{-1}:1$};
  \node[draw] at (1,1) (22) {\tiny $412563^{-1}:1$};
  \node[draw] at (2,1) (23) {\tiny $41523^{-1}:1$};
  \node[draw] at (0,2) (31) {\tiny $142563^{-1}:2$};
  \node[draw] at (1,2) (32) {\tiny $14523^{-1}:2$};
  \node[draw] at (2,2) (33) {\tiny $41253^{-1}:2$};
  \node[draw,fill=lightskyblue] at (1,3) (41) {\tiny $14253^{-1}:4$};
  \draw[->,thick]  (41) -- (31);
  \draw[->,thick]  (41) -- (32);
  \draw[->,thick]  (41) -- (33);
  \draw[->,thick]  (31) -- (21);
  \draw[->,thick]  (31) -- (22);
  \draw[->,thick]  (32) -- (21);
  \draw[->,thick]  (32) -- (23);
  \draw[->,thick]  (33) -- (22);
  \draw[->,thick]  (33) -- (23);
  \draw[->,thick]  (21) -- (11);
  \draw[->,thick]  (22) -- (11);
  \draw[->,thick]  (23) -- (11);
 \end{tikzpicture}
 \quad
\begin{tikzpicture}[xscale=2.1, yscale=1.5,>=latex,baseline=(z.base)]
\node at (0,2) (z) {};
  \node[draw] at (2.5,0) (11) {\tiny $51627384^{-1}:1$};
  \node[draw] at (1,1) (21) {\tiny $51623784^{-1}:1$};
  \node[draw] at (2,1) (22) {\tiny $51267384^{-1}:1$};
  \node[draw] at (3,1) (23) {\tiny $15627384^{-1}:1$};
  \node[draw] at (4,1) (24) {\tiny $5162734^{-1}:1$};
  \node[draw] at (0,2) (31) {\tiny $51263784^{-1}:2$};
  \node[draw] at (1,2) (32) {\tiny $15623784^{-1}:2$};
  \node[draw] at (2,2) (33) {\tiny $5162374^{-1}:2$};
  \node[draw] at (3,2) (34) {\tiny $15267384^{-1}:2$};
  \node[draw] at (4,2) (35) {\tiny $5126734^{-1}:2$};
  \node[draw] at (5,2) (36) {\tiny $1562734^{-1}:2$};
  \node[draw] at (1,3) (41) {\tiny $15263784^{-1}:4$};
  \node[draw] at (2,3) (42) {\tiny $5126374^{-1}:4$};
  \node[draw] at (3,3) (43) {\tiny $1562374^{-1}:4$};
  \node[draw] at (4,3) (44) {\tiny $1526734^{-1}:4$};
  \node[draw,fill=lightskyblue] at (2.5,4) (51) {\tiny $1526374^{-1}:8$};
  \draw[->,thick]  (51) -- (41);
  \draw[->,thick]  (51) -- (42);
  \draw[->,thick]  (51) -- (43);
  \draw[->,thick]  (51) -- (44);
  \draw[->,thick]  (41) -- (31);
  \draw[->,thick]  (41) -- (32);
  \draw[->,thick]  (41) -- (34);
  \draw[->,thick]  (42) -- (31);
  \draw[->,thick]  (42) -- (33);
  \draw[->,thick]  (42) -- (35);
  \draw[->,thick]  (43) -- (32);
  \draw[->,thick]  (43) -- (33);
  \draw[->,thick]  (43) -- (36);
  \draw[->,thick]  (44) -- (34);
  \draw[->,thick]  (44) -- (35);
  \draw[->,thick]  (44) -- (36);
  \draw[->,thick]  (31) -- (21);
  \draw[->,thick]  (31) -- (22);
  \draw[->,thick]  (32) -- (21);
  \draw[->,thick]  (32) -- (23);
  \draw[->,thick]  (33) -- (21);
  \draw[->,thick]  (33) -- (24);
  \draw[->,thick]  (34) -- (22);
  \draw[->,thick]  (34) -- (23);
  \draw[->,thick]  (35) -- (22);
  \draw[->,thick]  (35) -- (24);
  \draw[->,thick]  (36) -- (23);
  \draw[->,thick]  (36) -- (24);
  \draw[->,thick]  (21) -- (11);
  \draw[->,thick]  (22) -- (11);
  \draw[->,thick]  (23) -- (11);
  \draw[->,thick]  (24) -- (11);
 \end{tikzpicture}}}

\caption{The directed graphs $ \cB^+(z)$ when $z$ is 
$g_{23}=(2,4)(3,5)$ (left)  and $g_{24}=(2,5)(3,6)(4,7)$ (right), presented using the conventions in Figure~\ref{t1n-fig}.
Here, the grey box indicates the unique element of $\cB^+(z)$  not in $\supp(\GC_z)$.}\label{g2n-fig}
\end{figure}

\newpage

\printbibliography
\end{document}